\documentclass[11pt]{article}
\usepackage{amsmath}
\usepackage{amssymb}
\usepackage{amscd}
\usepackage{xspace}
\usepackage{verbatim}
\usepackage{color}
\newcommand{\mysection}[1]{
\section{#1}\setcounter{equation}{0}}
\date{}
\begin{document}
\begin{center}{\bf \large Singularities of fractional Emden's equations\medskip

 via Caffarelli-Silvestre extension\footnote{To appear in {\bf Jl. Diff. Equ.}}\\[2mm]
}

 \end{center}
 \begin{center}{\bf  Huyuan Chen}\\
{\small Department of Mathematics, Jiangxi Normal University}\\
 {\small Nanchang, Jiangxi 330022, PR China} \\[2mm]

 \end{center}
\begin{center}{\bf  Laurent V\'eron}\\
{\small Institut Denis Poisson, CNRS UMR 7013}\\
 {\small Universit\'e de Tours, 37200 Tours, France} \\[2mm]
 \end{center}



\newcommand{\txt}[1]{\;\text{ #1 }\;}
\newcommand{\tbf}{\textbf}
\newcommand{\tit}{\textit}
\newcommand{\tsc}{\textsc}
\newcommand{\trm}{\textrm}
\newcommand{\mbf}{\mathbf}
\newcommand{\mrm}{\mathrm}
\newcommand{\bsym}{\boldsymbol}
\newcommand{\scs}{\scriptstyle}
\newcommand{\sss}{\scriptscriptstyle}
\newcommand{\txts}{\textstyle}
\newcommand{\dsps}{\displaystyle}
\newcommand{\fnz}{\footnotesize}
\newcommand{\scz}{\scriptsize}
\newcommand{\be}{
\begin{equation}
}
\newcommand{\bel}[1]{
\begin{equation}
\label{#1}}
\newcommand{\ee}{
\end{equation}
}
\newcommand{\eqnl}[2]{
\begin{equation}
\label{#1}{#2}
\end{equation}
}
\newtheorem{subn}{\name}
\newcommand{\bsn}[1]{\def\name{#1}
\begin{subn}}
\newcommand{\esn}{
\end{subn}}
\newtheorem{sub}{\name}[section]
\newcommand{\dn}[1]{\def\name{#1}}   
\newcommand{\bs}{
\begin{sub}}
\newcommand{\es}{
\end{sub}}
\newcommand{\bsl}[1]{
\begin{sub}\label{#1}}
\newcommand{\bth}[1]{\def\name{Theorem}
\begin{sub}\label{t:#1}}
\newcommand{\blemma}[1]{\def\name{Lemma}
\begin{sub}\label{l:#1}}
\newcommand{\bcor}[1]{\def\name{Corollary}
\begin{sub}\label{c:#1}}
\newcommand{\bdef}[1]{\def\name{Definition}
\begin{sub}\label{d:#1}}
\newcommand{\bprop}[1]{\def\name{Proposition}
\begin{sub}\label{p:#1}}
\newcommand{\BBR}{\eqref}
\newcommand{\rth}[1]{Theorem~\ref{t:#1}}
\newcommand{\rlemma}[1]{Lemma~\ref{l:#1}}
\newcommand{\rcor}[1]{Corollary~\ref{c:#1}}
\newcommand{\rdef}[1]{Definition~\ref{d:#1}}
\newcommand{\rprop}[1]{Proposition~\ref{p:#1}} 
\newcommand{\BA}{
\begin{array}}
\newcommand{\EA}{
\end{array}}
\newcommand{\BAN}{\renewcommand{\arraystretch}{1.2}
\setlength{\arraycolsep}{2pt}
\begin{array}}
\newcommand{\BAV}[2]{\renewcommand{\arraystretch}{#1}
\setlength{\arraycolsep}{#2}
\begin{array}}
\newcommand{\BSA}{
\begin{subarray}}
\newcommand{\ESA}{\end{subarray}}
\newcommand{\BAL}{\begin{aligned}}
\newcommand{\EAL}{\end{aligned}}
\newcommand{\BALG}{\begin{alignat}}
\newcommand{\EALG}{\end{alignat}}
\newcommand{\BALGN}{\begin{alignat*}}
\newcommand{\EALGN}{\end{alignat*}}
\newcommand{\note}[1]{\textit{#1.}\hspace{2mm}}
\newcommand{\Proof}{\noindent{\bf Proof. } }
\newcommand{\qeda}{\hspace{10mm}\hfill $\square$}
\newcommand{\qed}{\\
${}$ \hfill $\square$}
\newcommand{\BBRemark}{{\bf Remark.} }
\newcommand{\modin}{$\,$\\
[-4mm] \indent}
\newcommand{\forevery}{\quad \forall}
\newcommand{\set}[1]{\{#1\}}
\newcommand{\setdef}[2]{\{\,#1:\,#2\,\}}
\newcommand{\setm}[2]{\{\,#1\mid #2\,\}}
\newcommand{\lra}{\longrightarrow}
\newcommand{\sgn}{\rm{sgn}}
\newcommand{\lla}{\longleftarrow}
\newcommand{\llra}{\longleftrightarrow}
\newcommand{\Lra}{\Longrightarrow}
\newcommand{\Lla}{\Longleftarrow}
\newcommand{\Llra}{\Longleftrightarrow}
\newcommand{\warrow}{\rightharpoonup}
\newcommand{
\paran}[1]{\left (#1 \right )}
\newcommand{\sqbr}[1]{\left [#1 \right ]}
\newcommand{\curlybr}[1]{\left \{#1 \right \}}
\newcommand{\abs}[1]{\left |#1\right |}
\newcommand{\norm}[1]{\left \|#1\right \|}
\newcommand{
\paranb}[1]{\big (#1 \big )}
\newcommand{\lsqbrb}[1]{\big [#1 \big ]}
\newcommand{\lcurlybrb}[1]{\big \{#1 \big \}}
\newcommand{\absb}[1]{\big |#1\big |}
\newcommand{\normb}[1]{\big \|#1\big \|}
\newcommand{
\paranB}[1]{\Big (#1 \Big )}
\newcommand{\absB}[1]{\Big |#1\Big |}
\newcommand{\normB}[1]{\Big \|#1\Big \|}

\newcommand{\thkl}{\rule[-.5mm]{.3mm}{3mm}}
\newcommand{\thknorm}[1]{\thkl #1 \thkl\,}
\newcommand{\trinorm}[1]{|\!|\!| #1 |\!|\!|\,}
\newcommand{\bang}[1]{\langle #1 \rangle}
\def\angb<#1>{\langle #1 \rangle}
\newcommand{\vstrut}[1]{\rule{0mm}{#1}}
\newcommand{\rec}[1]{\frac{1}{#1}}
\newcommand{\opname}[1]{\mbox{\rm #1}\,}
\newcommand{\supp}{\opname{supp}}
\newcommand{\dist}{\opname{dist}}
\newcommand{\myfrac}[2]{{\displaystyle \frac{#1}{#2} }}
\newcommand{\myint}[2]{{\displaystyle \int_{#1}^{#2}}}
\newcommand{\mysum}[2]{{\displaystyle \sum_{#1}^{#2}}}
\newcommand {\dint}{{\displaystyle \int\!\!\int}}
\newcommand{\q}{\quad}
\newcommand{\qq}{\qquad}
\newcommand{\hsp}[1]{\hspace{#1mm}}
\newcommand{\vsp}[1]{\vspace{#1mm}}
\newcommand{\ity}{\infty}
\newcommand{\prt}{
\partial}
\newcommand{\sms}{\setminus}
\newcommand{\ems}{\emptyset}
\newcommand{\ti}{\times}
\newcommand{\pr}{^\prime}
\newcommand{\ppr}{^{\prime\prime}}
\newcommand{\tl}{\tilde}
\newcommand{\sbs}{\subset}
\newcommand{\sbeq}{\subseteq}
\newcommand{\nind}{\noindent}
\newcommand{\ind}{\indent}
\newcommand{\ovl}{\overline}
\newcommand{\unl}{\underline}
\newcommand{\nin}{\not\in}
\newcommand{\pfrac}[2]{\genfrac{(}{)}{}{}{#1}{#2}}

\def\ga{\alpha}     \def\gb{\beta}       \def\gg{\gamma}
\def\gc{\chi}       \def\gd{\delta}      \def\ge{\epsilon}
\def\gth{\theta}                         \def\vge{\varepsilon}
\def\gf{\phi}       \def\vgf{\varphi}    \def\gh{\eta}
\def\gi{\iota}      \def\gk{\kappa}      \def\gl{\lambda}
\def\gm{\mu}        \def\gn{\nu}         \def\gp{\pi}
\def\vgp{\varpi}    \def\gr{\rho}        \def\vgr{\varrho}
\def\gs{\sigma}     \def\vgs{\varsigma}  \def\gt{\tau}
\def\gu{\upsilon}   \def\gv{\vartheta}   \def\gw{\omega}
\def\gx{\xi}        \def\gy{\psi}        \def\gz{\zeta}
\def\Gg{\Gamma}     \def\Gd{\Delta}      \def\Gf{\Phi}
\def\Gth{\Theta}
\def\Gl{\Lambda}    \def\Gs{\Sigma}      \def\Gp{\Pi}
\def\Gw{\Omega}     \def\Gx{\Xi}         \def\Gy{\Psi}

\def\CS{{\mathcal S}}   \def\CM{{\mathcal M}}   \def\CN{{\mathcal N}}
\def\CR{{\mathcal R}}   \def\CO{{\mathcal O}}   \def\CP{{\mathcal P}}
\def\CA{{\mathcal A}}   \def\CB{{\mathcal B}}   \def\CC{{\mathcal C}}
\def\CD{{\mathcal D}}   \def\CE{{\mathcal E}}   \def\CF{{\mathcal F}}
\def\CG{{\mathcal G}}   \def\CH{{\mathcal H}}   \def\CI{{\mathcal I}}
\def\CJ{{\mathcal J}}   \def\CK{{\mathcal K}}   \def\CL{{\mathcal L}}
\def\CT{{\mathcal T}}   \def\CU{{\mathcal U}}   \def\CV{{\mathcal V}}
\def\CZ{{\mathcal Z}}   \def\CX{{\mathcal X}}   \def\CY{{\mathcal Y}}
\def\CW{{\mathcal W}} \def\CQ{{\mathcal Q}} 
\def\BBA {\mathbb A}   \def\BBb {\mathbb B}    \def\BBC {\mathbb C}
\def\BBD {\mathbb D}   \def\BBE {\mathbb E}    \def\BBF {\mathbb F}
\def\BBG {\mathbb G}   \def\BBH {\mathbb H}    \def\BBI {\mathbb I}
\def\BBJ {\mathbb J}   \def\BBK {\mathbb K}    \def\BBL {\mathbb L}
\def\BBM {\mathbb M}   \def\BBN {\mathbb N}    \def\BBO {\mathbb O}
\def\BBP {\mathbb P}   \def\BBR {\mathbb R}    \def\BBS {\mathbb S}
\def\BBT {\mathbb T}   \def\BBU {\mathbb U}    \def\BBV {\mathbb V}
\def\BBW {\mathbb W}   \def\BBX {\mathbb X}    \def\BBY {\mathbb Y}
\def\BBZ {\mathbb Z}

\def\GTA {\mathfrak A}   \def\GTB {\mathfrak B}    \def\GTC {\mathfrak C}
\def\GTD {\mathfrak D}   \def\GTE {\mathfrak E}    \def\GTF {\mathfrak F}
\def\GTG {\mathfrak G}   \def\GTH {\mathfrak H}    \def\GTI {\mathfrak I}
\def\GTJ {\mathfrak J}   \def\GTK {\mathfrak K}    \def\GTL {\mathfrak L}
\def\GTM {\mathfrak M}   \def\GTN {\mathfrak N}    \def\GTO {\mathfrak O}
\def\GTP {\mathfrak P}   \def\GTR {\mathfrak R}    \def\GTS {\mathfrak S}
\def\GTT {\mathfrak T}   \def\GTU {\mathfrak U}    \def\GTV {\mathfrak V}
\def\GTW {\mathfrak W}   \def\GTX {\mathfrak X}    \def\GTY {\mathfrak Y}
\def\GTZ {\mathfrak Z}   \def\GTQ {\mathfrak Q}

\font\Sym= msam10 
\def\SYM#1{\hbox{\Sym #1}}
\newcommand{\bdw}{\prt\Gw\xspace}
\medskip

{\abstract We study the isolated singularities of functions satisfying 
$$(E) \qquad\qquad\qquad \quad  (-\Gd)^s v \pm|v|^{p-1}v=0\quad {\rm in}\ \, \Gw\setminus\{0\},\qquad v=0\quad {\rm in}\ \, \BBR^N\setminus\Omega,\qquad\qquad\qquad\qquad\qquad\qquad\qquad  $$
  where $0<s<1$,  $p>1$ and $\Gw$ is a bounded domain containing the origin. We use the Caffarelli-Silvestre extension to $\BBR_+\ti\BBR^N$. We emphasize the obtention of {\it a priori} estimates and analyse the set of self-similar solutions via energy methods to characterize the singularities.}\bigskip

\nind {\it AMS Subject Classification:} 35B40, 35J60, 35J70, 45G05\smallskip

\nind{\it  Keywords:} Fractional Laplacian; Isolated singularity; Lane-Emden-Fowler equation; Linear extension; Energy methods; Limit set.

\tableofcontents
\date{}

\mysection{Introduction}
In the past forty years a large series of articles have been devoted to study of the singular behaviour of solutions for the following two classes of semilinear equations 
\bel{IN0}\BA{lll}
-\Gd v+\ge v^p=0\quad {\rm in}\ \, B_1\setminus\{0\},
\EA\ee
where $p>1$, $\ge=\pm 1$ and $v^p=|v|^{p-1}v$. The first empirical studies of radial solutions to the Lane-Emden equation ($\ge=-1$) are due to J. Lane and R. Emden when analyzing the structure of polytropic objects submitted to their gravitational fields. A comprehensive presentation can be found in Chandrasekhar's book \cite[pp 84-182]{Ch}. The so-called Emden-Fowler equation 
($\ge=1$), was treated in details by R. Fowler in the radial case. For some specific values of $p$ ($\frac 32$ and $\frac 52$) this equation appeared in the study of the density of the electronic cloud in the Thomas-Fermi-Dirac model of atoms (see e.g. \cite{So}).

The study of non-radial solutions started in the heighties in connection with conformal deformation of a Riemannian metric to a metric with constant scalar curvature ($K=1$ or $K=-1$ and $p=\frac{N+2}{N-2}$). The first important results concerning the Lane-Emden equation are due to Lions \cite{Li} in the case $1<p<\frac{N}{N-2}$ who proved that any positive solution of 
\bel{IN1}\BA{lll}
-\Gd v-v^p=0\quad {\rm in}\ \,   B_1\setminus\{0\}
\EA\ee
satisfies 
\bel{IN2}\BA{lll}
v(x)=k|x|^{2-N}+O(1)\quad\text{as }x\to 0
\EA\ee
when $N\geq 3$, (with standard modification if $N=2$) and the following holds
$$
-\Gd v-v^p=c_Nk\gd_0\in \CD'(B_1).
$$
In the case $\frac{N}{N-2}<p<\frac{N+2}{N-2}$, Gidas and Spruck \cite{GS} introduced new powerful methods and proved that any positive solution of $(\ref{IN1})$ satisfies
$$
v(x)\leq c_{p}|x|^{-\frac 2{p-1}}\quad\text{for\,}\  x\in B_{\frac12}\setminus\{0\}. 
$$
As a consequence, a positive solution is either regular at $x=0$ or has the following behaviour
\bel{IN5}\BA{lll}
v(x)=\ga_{{\scriptscriptstyle N},p}|x|^{-\frac 2{p-1}}\big(1+o(1)\big)\quad\text{as }x\to 0.
\EA\ee
When $p=\frac{N+2}{N-2}$, the energy method used in \cite{GS} failed because of the conformal invariance of the equation. In \cite{CGS} a new method of asymptotic symmetry allowed to treat this case showing that the positive solutions behaved like the radial ones, which in turn were easy to describe because of this invariance. Finally, when $p=\frac{N}{N-2}$ the result of \cite{GS} was sharpened by Aviles \cite{Av} who proved that either the positive solutions are smooth or they satisfy
$$
\lim_{x\to 0}|x|^{N-2}(-\ln |x| )^{\frac{N-2}{2}}v(x)=\Big(\frac{N-2}{\sqrt 2}\Big)^{N-2}.
$$

The non-radial study of the Emden-Fowler equations started with the works of Brezis and Lieb who obtained in \cite{BLb} the universal {\it a priori} estimate (also called Keller-Osserman estimate), valid for signed solutions of 
\bel{IN7}\BA{lll}
-\Gd v+ v^p=0\quad {\rm in}\ \,   B_1\setminus\{0\},
\EA\ee
and is expressed as follows:
$$
|v(x)|\leq c'_{p}|x|^{-\frac 2{p-1}}\quad\text{for\ any}\   x\in B_{\frac12}\setminus\{0\}.  
$$
With this estimate, Brezis and V\'eron proved in \cite{BV} that if $p\geq\frac N{N-2}$ any solution is at least $C^2$ in $B_1$. When $1<p<\frac{N}{N-2}$ the precise description of singular solutions was obtained by V\'eron who proved that any positive solution $v$ is either smooth, either satisfies $(\ref{IN2})$  or 
$(\ref{IN5})$ (with a new constant positive $\tilde\ga_{N,p}$), all these behaviours being effective. Furthermore when $\frac{N+1}{N-1}\leq p<\frac{N}{N-2}$ the assumption of positivity could be removed, the sign of $k$ not prescribed in $(\ref{IN2})$, and 
$(\ref{IN5})$ is replaced by there exists $\ell\in\BBR$ such that
$$
\displaystyle
\lim_{x\to 0}|x|^{\frac{2}{p-1}}v(x)=\ell\in\{\tilde\ga_{N,p},-\tilde\ga_{N,p},0\}\quad\text{as }x\to 0.
$$
In the case $N=2$, the requirement $p\geq \frac{N+1}{N-1}$ (i.e. $q\geq 3$) for signed solutions was removed by Chen, Matano and V\'eron \cite{CMV} who proved that for any solution $v$ of $(\ref{IN7})$  there exists a $2\gp$-periodic solution of 
$$
\displaystyle
-\gw_{\gth\gth}-\left(\frac{2}{p-1}\right)^2\gw+ \gw^p=0
$$
such that 
$$
\displaystyle
\lim_{x\to 0}|x|^{\frac{2}{p-1}}v(x)=\lim_{r\to 0}r^{\frac{2}{p-1}}v(r,\gth)=\gw(\gth)\quad\text{as }r=|x|\to 0.
$$

This panorama of the classical Emden equations intended to show what could be expected when the ordinary Laplacian is replaced by the fractional Laplacian in these equations. If $s\in (0,1)$ and $(-\Gd)^s$  the fractional Laplacian in $\BBR^N\setminus\{0\}$ is defined on functions $u\in C^2(\BBR^N\setminus\{0\})\cap L^1_{\gm_s}(\BBR^N)$ with $\gm_s(x)=(1+|x|)^{-N-2s}$ by the expression
\bel{I-1}\displaystyle
(-\Gd)^su(x)=c_{{\scriptscriptstyle N},s}\lim_{\vge\to 0}\myint{|x-y|>\vge}{}\myfrac{u(x)-u(y)}{|x-y|^{N+2s}}dy\quad \text{for }\, x\in\BBR^N\setminus\{0\},
\ee
where $c_{{\scriptscriptstyle N},s}=2^{2s}\gp^{-\frac N2}\frac{\Gg(\frac{N+2s}{2})}{\Gg(1-s)}$. The singularity problem for the  fractional Emden  equations  in a punctured domain $\Gw$ containing $\overline{B}_1$ is 
\bel{IN12}\left.\BA{lll}
(-\Gd)^sv+\ge v^p=0\quad \text{in }\,\Gw\setminus\{0\}\\[2mm]
\phantom{(-\Gd)^s+\ge v^p}
v=0\quad \text{in }\, \Gw^c,
\EA\right.\ee
where $\ge=\pm 1$ and $\Gw^c=\BBR^N\setminus\Gw$. Some results concerning this equation in the case $\ge=-1$ (the fractional Lane-Emden equation) have already been obtained in particular in the case $p=\frac{N+2s}{N-2s}$, see \cite{CJSX}. Therein it is proved that either the solution is smooth at $0$, or there exists some constant $c\geq1$ such that 
$$
\myfrac{1}{c}|x|^{-\frac{N-2s}{2}}\leq v(x)\leq c|x|^{-\frac{N-2s}{2}}\quad\text{for all }x\in B_{\frac12}\setminus\{0\}.
$$
Furthermore, there exists a radial solution $\tilde v$ of $(\ref{IN12})$ such that 
$$
v(x)=\tilde v(|x|)(1+o(1))\quad\text{as }\, x\to 0.
$$
In the general case $\frac{N}{N-2s}<p<\frac{N+2s}{N-2s}$, Yang and Zou obtained in \cite{YZ} (also in \cite{YZ1}) the analogous of Gidas and Spruck estimates, namely 
any positive solution of 
\bel{IN15}\left.\BA{lll}
(-\Gd)^sv-v^p=0\quad \text{in }\, \Gw\setminus\{0\}\\[2mm]
\phantom{(-\Gd)^s+  v^p}
v=0\quad \text{in }\,\Gw^c,
\EA\right.\ee
either is smooth or satisfies for some constant $c>0$ 
$$
\myfrac{1}{c}|x|^{-\frac{2s}{p-1}}\leq v(x)\leq c|x|^{-\frac{2s}{p-1}}\quad\text{for all }x\in B_{\frac12}\setminus\{0\}.
$$
Concerning the fractional Emden-Fowler equation
 \bel{IN17}\left.\BA{lll}
(-\Gd)^sv+ v^p=0\quad \text{in }\,\Gw\setminus\{0\}\\[2mm]
\phantom{(-\Gd)^s+  v^p}
v=0\quad \text{in }\, \Gw^c
\EA\right.\ee
not so many results are known. The problem 
$$
\left.\BA{lll}
(-\Gd)^sv+v^p=\gm\quad \text{in }\,\CD'(\Gw)\\[2mm]
\phantom{(-\Gd)^s+  v^p}
v=0\quad \text{in }\,\Gw^c,
\EA\right.$$
where $\gm$ is a positive bounded Radon measure in $\Gw$, is studied in \cite{CV} in a more general framework, that is by replacing $u^p$ by a nondecreasing nonlinearity 
$g(u)$. The particular case where $\gm=k\gd_0$ and $g(u)=u^p$ is analyzed in \cite{CV1}. In particular, the asymptotics of the solutions $u_k$ (they are unique) when $k\to\infty$ are thoroughly studied in \cite{CV1}. \medskip

The aim of this article is to present a unified analysis of the isolated singularities  for problem $(\ref{IN12})$. This analysis is based upon the Caffarelli-Silvestre lifting 
of the equation which associates to it the following degenerate elliptic equation in $\BBR^{N+1}_+=\big\{\xi=(x,z):x\in\BBR^N,z>0\big\},$ 
\bel{I-2-0}\left.\BA{lll}\displaystyle
{\rm div}(z^{1-2s}\nabla u)=0\quad&\text{in }\,\BBR^{N+1}_+\\[2mm]
\phantom{div(^{1,2s}}
u(\cdot,0)=v \quad\quad&\text{in }\BBR^N\sim\prt \BBR^{N+1}_+,
\EA\right.
\ee
with the property that 
$$
(-\Gd)^sv(x)=-\lim_{z\to 0}z^{1-2s}u_z(x,z):=\displaystyle\prt_{\gn^s}u(x,0).
$$
The study of $(\ref{IN12})$ is replaced by 
\bel{IN19}
 \BA{lll}\displaystyle
\qquad\ \ \, {\rm div}(z^{1-2s}\nabla u)=0\quad&\text{in }\,\BBR^{N+1}_+\\[2mm]
\prt_{\gn^s}u(\cdot,0) +\ge u(\cdot,0)^p=0\quad &\text{in } \, \Gw\setminus\{0\}\\[2mm]
\phantom{(-\Gd)^s+\ge v^p}
\ \  u(\cdot,0)=0\quad &\text{in } \, \BBR^N\setminus \Gw,
\EA
\ee
and
 $$v=u(\cdot,0)\quad \text{in }\, \BBR^N.$$

We prove that any positive solution of $(\ref{I-2-0})$ admits a trace on $\BBR^N$ which is a nonnegative Radon measure $\gm$ satisfying that  
$$
\int_{\BBR^N}\frac{d\gm}{(1+|x|)^{N+2s}}<+\infty. 
$$

The first step for describing the behaviour of solutions of $(\ref{IN19})$ near $0$  is to obtain an {\it  a priori} estimate of $u(x,z)$.  \smallskip

  In the case $\ge=1$, we use a blow-up technique combined with the trace theorem to  prove that when $p\geq 1+\frac{2s}{N}$, any positive solution $u(x,z)$  of $(\ref{IN19})$ satisfies 
\bel{IN21}\BA{lll}
u(x,z)\leq c\gr^{-\frac{2s}{p-1}}\quad\text{for all }(x,z)\in \tilde {\bf B}^+_{\frac12}\setminus\{0\},
\EA\ee
where $\gr=\sqrt{|x|^2+z^2}$ and  $\tilde {\bf B}^+_a:=\big\{(x,z)\in\overline{\BBR^{N+1}_+}:\gr<a\big\}$. This estimate implies that any positive solution of $(\ref{IN17})$ verifies
\bel{IN22}\BA{lll}
v(x)\leq c|x|^{-\frac{2s}{p-1}}\quad\text{for all }x\in \tilde B^+_{\frac12}\setminus\{0\}. 
\EA\ee
Note that this estimate has no interest if $p\geq \frac{N}{N-2s}$ since the solution $v$ is smooth. \smallskip

   In the case $\ge=-1$, we use the estimate $(\ref{IN22})$ proved  \cite[Proposition 3.1]{YZ} in the case $1<p<\frac{N+2s}{N-2s}$ combined with the expression of the Poisson kernel of the operator $u\mapsto div(z^{1-2s}\nabla u)$ in $\BBR^{N+1}_+$ obtained by \cite{CS} to prove that positive solutions of 
$(\ref{IN19})$ satisfy also $(\ref{IN21})$. Furthermore, when $v$ is a radially symmetric decreasing solution of $(\ref{IN15})$,  we prove that $(\ref{IN22})$ holds for any 
$p>1$. \smallskip 

The second step in our study is to consider  the self-similar solutions of $(\ref{IN19})$ in $\overline{\BBR^{N+1}_+}\setminus\{0\}$, with $\Gw=\BBR^N$. In spherical coordinates in $\BBR^{N+1}_+$, these solutions have the following expression 
$$
u(x,z)=u(\gr,\gs)=\gr^{-\frac{2s}{p-1}}\gw(\gs)\quad\text{for all }(\gr,\gs)\in\BBR_+\ti \BBS^{N}_+,
$$
and, up to a rotation and a good choice of spherical variables on $\BBS^{N}$, $\gw$ satisfies 
\bel{IN23}\left.\BA{lll}\displaystyle
\phantom{}\CA_s[\gw] +\Gl_{s,p,{\scriptscriptstyle N}}\gw=0\ \ &\text{in }\,\BBS^N_+\\[2mm]
\displaystyle\phantom{}
\frac{\prt \gw}{\prt \gn^s}+\ge|\gw|^{p-1}\gw=0\ \ &\text{in }\,\BBS^{N-1}.
\EA\right.\ee
Therein, $\CA_s$ is a degenerate elliptic operator on the $N$-sphere $\BBS^N$, $\frac{\prt }{\prt \gn^s}$ the corresponding conormal outward derivative on 
$\prt \BBS^N_+=\BBS^{N-1}$ and 
$$
\Gl_{s,p,{\scriptscriptstyle N}}=\frac{2s}{p-1}\Big(\frac{2s}{p-1}+2s-N\Big).
$$
The sign of $\Gl_{s,p,{\scriptscriptstyle N}}$ which is fundamental in the study of $(\ref{IN23})$ depends on the value of $p$ with respect to $\frac{2s}{p-1}$. The structure of the set $\CE_\ge$ (resp. $\CE^+_\ge$) of solutions (resp. positive solutions) plays a key role in our study.\medskip

\nind{\bf Theorem A} {\it Let $s\in (0,1)$, $\ge=1$ and $p>1$.\smallskip

\nind 1- If $p\geq \frac{N}{N-2s}$, then $\CE_1=\{0\}$. \smallskip

\nind 2- If $1<p\leq 1+\frac{2s}{N}$, then $\CE^+_1=\{0\}$. \smallskip

\nind 3- If $1+\frac{2s}{N}<p<\frac{N}{N-2s}$, then $\CE^+_1=\{0,\gw_1\}$, where $\gw_1$ is a positive solution of $(\ref{IN23})$. 
}\medskip

Besides the case $p\geq\frac{N}{N-2s}$ we can describe the set $\CE_1$ in another case.\medskip

\nind{\bf Theorem B} {\it There exists $p^*=p^*(s,N)\in (1,\frac{N}{N-2s})$ such that if $p^*\leq p<\frac{N}{N-2s}$, then $\CE_1=\{0,\gw_1,-\gw_1\}$}.\medskip

The exponent $p^*$ is explicitly given by $(\ref{XX1})$ below  and its  algebraic expression is heavy. However if $N\geq 3$ or $N=2$ and $\frac12\leq s<1$ one has $p^*\leq 1+\frac{2s}{N}$. The set of positive self-similar solutions of the Lane-Emden is characterized by the following statement.\medskip

\nind{\bf Theorem C} {\it Let $s\in (0,1)$, $\ge=-1$ and $p>1$.\smallskip

\nind 1- If $p\leq \frac{N}{N-2s}$, then $\CE^+_{-1}=\{0\}$. \smallskip

\nind 2- If $p>\frac{N}{N-2s}$, then $\CE^+_{-1}=\{0,\gw_2\}$, where $\gw_2$ is a positive solution of   $(\ref{IN23})$ depending only on one variable. }
\medskip

We use the classical $\ln\gr$ variable in order to transform  $(\ref{IN23})$ into an autonomous equation, writing any $u$ under the form 
\bel{IN24}\BA{lll}\displaystyle
u(x,z)=\gr^{-\frac{2s}{p-1}}w(t,\gs)\quad\ \text{with }\, t=\ln\gr.
\EA\ee
We introduce standard tools of the dynamical systems theory such as energy functional and limit sets and prove the following general theorem. \medskip

\nind{\bf Theorem D} {\it Assume $s\in (0,1)$, $\ge=\pm1$, $p\in (1,+\infty)\setminus\{\frac{N+2s}{N-2s}\}$. Let $u$ be a solution of $(\ref{IN19})$ satisfying $(\ref{IN21})$, even with $u$ replaced by $|u|$ and $w$ be defined by $(\ref{IN24})$. Then the limit set at $-\infty$ of the negative trajectory $\CT_-[w]:=\bigcup \big\{w(t,\cdot):t<0\big\}$ of $w$ 
is a nonempty, compact and connected subset of $\CE^+_{\ge}$ in the $C^{2s}(\overline {\BBS^N_+})$-topology.
}\medskip
 
 Thanks to this result and the properties of $\CE^+_{\ge}$, we can give the precise behaviour of positive solutions of $(\ref{IN12})$ near the origin. This happens
 when $\CE_\ge$ is disconnected, the limit set is reduced to a single element usually $\gw_j$ (j=1,2) or $0$. When this limit set is zero, the solutions have either a weak singularity or regular, according to the value of $\ge$.\medskip
 
\nind{\bf Theorem E} {\it Let $s\in (0,1)$, $\ge=1$ and $1+\frac{2s}{N}<p<\frac {N}{N-2s}$. If $u$ is a positive solution of $(\ref{IN19})$ satisfying
 \bel{IN25}\displaystyle\lim_{(x,z)\to 0}\gr^{\frac{2s}{p-1}}u(x,z)=0,
\ee
then there exists $k>0$ such that 
$$
u(x,z)=C_{{\scriptscriptstyle N},s}k\frac{z^{2s}}{(|x|^2+z^2)^{\frac N2+s}}+O(1)\quad\ \text{as }\,x\to 0,
$$
and $v_k =u(\cdot,0)$ satisfies
\bel{IN27}\BA{lll}\displaystyle
(-\Gd)^sv+v^p=k\gd_0\quad\ \text{in }\,\CD'(\Gw),
\EA\ee
where 
$C_{{\scriptscriptstyle N},s}=\frac14 \pi^{\frac{N+2-2s}{2}} \Gamma(\frac{N-2s}{2}).$

Furthermore, if $k=0$, $v_0$ is identically zero and $u$ is smooth. 
}\medskip

\nind{\bf Theorem F} {\it Let $s\in (0,1)$, $\ge=-1$ and $p\in \big(\frac{N}{N-2s},+\infty\big)\setminus\{\frac{N+2s}{N-2s}\}$. If $u$ is a positive solution of $(\ref{IN19})$ satisfying
$(\ref{IN25})$, then $u$ is smooth in  $\overline{\BBR^{N+1}_+}$.}\medskip

 Finally, we can give precise classification of nonnegative solutions to the fractional Emden-Fowler equation. \smallskip

\nind{\bf Theorem G} {\it Let $s\in (0,1)$,  $p>1$ and $v_k$ be the solution of $(\ref{IN27})$ for $k\geq0$ vanishing in $\Gw^c$. \smallskip

\nind 1- For $p\in \big(1, 1+\frac{2s}{N} \big]$,  the set of nonnegative solutions of $(\ref{IN17})$ is the set of $\{v_k\}_{k\geq0}$. \smallskip

\nind 2- For $p\in \big(1+\frac{2s}{N}, \frac{N}{N-2s}\big)$, the set of nonnegative solutions of $(\ref{IN17})$ is the set of $\{v_k\}_{k\geq0}\cup \{v_\infty\}$, where  $\displaystyle v_\infty=\lim_{k\to+\infty} v_k$. 

\nind 3- For $p  \geq \frac{N}{N-2s}$, the set of nonnegative solutions to $(\ref{IN17})$ is  reduced to $\{0\}$. 

 }\medskip

 The remainder of this paper is organized as follows. In Section 2, we provide  preliminary tools: Caffarelli-Silvestre extension, the related Poisson kernel and basic setting.      In Section 3, we classify     singular self-similar solutions of $(\ref{IN12})$.
  Section 4 is devoted to obtain some piori estimates the upper bounds of singular solutions $(\ref{IN12})$ near the origin.   Finally, we give the proof of the classification of  isolated singularities of $(\ref{IN12})$. 


\mysection{The  lifting}

\subsection{Caffarelli-Silvestre extension and its Poisson's kernel}
Caffarelli and Silvestre \cite{CS} introduced a general method of lifting from $\BBR^N$ to $\BBR^N\ti(0,+\infty)$ allowing to define the fractional Laplacian as the weighted trace of the normal derivative. Setting 
$$\BBR^{N+1}_+=\big\{\xi=(x,z):x\in\BBR^N,z>0\big\},$$ 
they proved that if $v\in C^2(\BBR^N)\cap L^1_{\gm_s}(\BBR^N)$ and $u\in C^2(\BBR^{N+1}_+)$ satisfy
\bel{I-2}\left.\BA{lll}\displaystyle
\CD_su:={\rm div}(z^{1-2s}\nabla u)=0\quad&\text{in }\,\BBR^{N+1}_+\\[2mm]
\phantom{\CL_su:=div(z^{1 } \ \ }
u(\cdot,0)=v\quad &\text{in }\, \BBR^N,
\EA\right.
\ee
then 
$$
(-\Gd)^sv(x)=-\lim_{z\to 0}z^{1-2s}u_z(x,z).
$$
In the sequel we denote
$$
\displaystyle\prt_{\gn^s}u(x,0):=-\lim_{z\to 0}z^{1-2s}u_z(x,z).
$$
 %
Set 
\bel{I-5}\BA{lll}\displaystyle
\gm_s(x)=(1+|x|)^{-N-2s}\qquad\text{for }\, x\in\BBR^N.
\EA
\ee
The Poisson's kernel for the operator  $u\mapsto \CD_su$
 has been obtained in $\cite[(2.1)]{CS}$
\bel{I-6}\BA{lll}\displaystyle
\CP_s(x,z)=C_{{\scriptscriptstyle N},s}\frac{z^{2s}}{\left(|x|^2+z^2\right)^{\frac{N}{2}+s}},
\EA
\ee
where $C_{{\scriptscriptstyle N},s}=\frac14 \pi^{\frac{N+2-2s}{2}} \Gamma(\frac{N-2s}{2})$. 
If $v\in L_{\gm_s}^1(\BBR^N)$, then 
\bel{I-7}\BA{lll}\displaystyle
u(x,z):=\BBP_s[v](x,z)=C_{{\scriptscriptstyle N},s}z^{2s}\int_{\BBR^N}\frac{v(y)dy}{\left(|x-y|^2+z^2\right)^{\frac{N}{2}+s}}
\EA
\ee
satisfies $(\ref{I-2})$. 
\subsection{The trace theorem}
The counter part of Caffarelli-Silvestre theorem is the following.
\bth{trace} Let $s\in (0,1)$ and $u\in L^1_{loc}(\BBR_+^{N+1})$ satisfies $u\geq 0$ and
\bel{I-T-1}\BA{lll}\displaystyle
\CD_s u=0\quad\text{in }\, \BBR_+^{N+1}.
\EA
\ee
Then $u\in C^{\infty}(\BBR_+^{N+1})$ and $u$ admits a boundary trace on $\prt\BBR^{N+1}_+$, which is a nonnegative Radon measure $\gm$ in the sense that for any $\gz\in C^\infty_0(\BBR^N)$, 
$$
\lim_{z\to 0}\myint{\BBR^N}{}u(x,z)\gz(x) dx=\myint{\BBR^N}{}\gz d\gm(x),
$$
and $\gm$ satisfies
\bel{I-T-3}\BA{lll}\displaystyle
\myint{\BBR^N}{}\frac{d\gm(x)}{(|x|+1)^{N+2s}}<+\infty.
\EA
\ee
Furthermore, if $\norm{u(z,\cdot)}_{L^1(\BBR^N)}=o(z^{2s})$ when $z\to\infty$, then
\bel{I-T-4}\BA{lll}\displaystyle
u(x,z)=C_{{\scriptscriptstyle N},s} z^{2s}\myint{\BBR^N}{}\frac{d\gm(y)}{(|x-y|^2+z^2)^{\frac{N}{2}+s}}\qquad\text{for all }(x,z)\in \BBR_+^{N+1}.
\EA
\ee
\es
\Proof The fact that $u\in C^{\infty}(\BBR^{N+1})$ follows from the standard theory of elliptic operators. Furthermore, the equation $(\ref{I-T-1})$ not only holds in the sense of distributions but also in the strong sense. The core of the proof is an adaptation of the Brezis-Lions wellknown note on isolated singularities of linear elliptic equations \cite{BL}.

Let $R>0$,  $\gl_{\scriptscriptstyle R,1}>0$ be first eigenvalue  of $-\Gd $ in $B_R\subset \BBR^N$
under the zero Dirichlet boundary condition and  $\gf_{\scriptscriptstyle R,1}$ be the related first eigenfunction. We set
$$X(z)=\int_{B_R}u(x,z)\gf_{\scriptscriptstyle R,1}(x)dx.
$$
Then $X$ satisfies the ODE
$$
(z^{1-2s}X')'(z)-\gl_{\scriptscriptstyle R,1} z^{1-2s}X(z)\leq 0\quad{\rm for}\ \,z\in(0,\infty).
$$
Integrating the equation on $(z,T)$, we obtain that 
$$
z^{1-2s}X'(z)\geq T^{1-2s}X'(T)-\gl_{\scriptscriptstyle R,1}\int_z^Tt^{1-2s}X(t)dt=T^{1-2s}X'(T)-\gl_{\scriptscriptstyle R,1} F(z),
$$
where we have denoted
$$F(z)=\int_z^Tt^{1-2s}X(t)dt.
$$
The function $F$ is decreasing and 
\bel{BL3}\BA{lll}\displaystyle
X(z)\leq X(T)-\frac{1}{2s}(T-T^{1-2s}z^{2s})X'(T)+\gl_{\scriptscriptstyle R,1}\int_z^Tt^{2s-1}F(t)dt\\[4mm]\displaystyle
\phantom{X(y)}
\leq X(T)-\frac{1}{2s}(T-T^{1-2s}z^{2s})X'(T)+\frac{\gl_{\scriptscriptstyle R,1}}{2s}(T^{2s}-z^{2s})F(z).
\EA\ee
Hence $F'(z)=-z^{1-2s}X(z)$ and $(\ref{BL3})$ becomes
\bel{BL4}\BA{lll}\displaystyle
-F'(z)\leq z^{1-2s}X(T)-\frac{1}{2s}(Tz^{1-2s}-T^{1-2s}z)X'(T)+\frac{\gl_{\scriptscriptstyle R,1}}{2s}\left(T^{2s}z^{1-2s}-z\right)F(z)\\[4mm]
\phantom{-F'(z)}
=: A(z)+B(z)F(z)
\EA\ee
for $z\in (0,\infty)$. Since $A,B\in L^1(0,T)$ and $F(T)=0$, we deduce by integration,
\bel{BL5}\BA{lll}\displaystyle
F(z)\leq e^{-\int_z^TB(t)dt}\int_z^T e^{\int_t^TB(\gt)d\gt}A(t)dt.
\EA\ee
Therefore $F$ is uniformly bounded on $(0,T)$, and from $(\ref{BL3})$ the function $X$ shares the same property.  This implies that $u\in L_{loc}^1(\BBR^N\ti [0,+\infty))$.

Next, let $\gz\in C^\infty_0(\BBR^N)$ with support in $B_R$ for some R. We set 
$$Y(z)=\int_{\BBR^N}u(x,z)\gz(x) dx\quad\text{and }\, \ \Gf(z)=\int_{\BBR^N}u(x,z)\Delta_x \gz(x) dx.
$$
Then 
$$\left(z^{1-2s}Y'(z)\right)'+z^{1-2s}\Gf(z)=0.
$$
Consequently, we have that 
$$Y(z)=Y(1)-\frac{1-z^{2s}}{2s}Y'(1)-\frac{1}{2s}\int_{z}^{1}(t^{2s}-z^{2s})\Gf(t)dt.
$$
Since 
$$\int_{z}^{1}|(t^{2s}-z^{2s})\Gf(t)|dt\leq\norm{\Gd\gz}_{L^\infty}\myint{0}{1}\int_{B_R}u(x,t)dxdt<+\infty,
$$
it follows that $Y(z)$ admits a limit when $z\to 0$, which defines a positive linear functional on $ C^\infty_0(\BBR^N)$, hence a Radon measure denoted by 
$\gm$. 

Let $R>0 $ and $\vge>0$. For $k>0$ we set 
$$Q^T_{R,R+k,\vge}=B_{R+k}\ti (\vge,T), $$
and we consider the problem 
$$
\left.\BA{lll}\displaystyle
 \ \ \, \CD_s w=0\qquad&\text{in }\  Q^T_{R,R+k,\vge}\\[1.5mm]
\phantom{ \CD_s- }
w=0&\text{in }\  \prt B_{R+k}\ti [\vge,T]\cup \overline B_{R+k}\ti\{T\}\\[1.5mm]
\phantom{ }
w(\cdot,\vge)=u(\cdot,\vge)\chi_{_{B_R}}&\text{in }\  B_{R+k}\ti\{\vge\}.
\EA\right.
$$
The solution $w=w_{\vge,{\scriptscriptstyle R},k,{\scriptscriptstyle T}}$ is unique and satisfies 
$$
w_{\vge,R,k,T}\leq u\quad\text{in }\, Q^T_{R,R+k,\vge}.
$$
The correspondence $(R,k,T)\mapsto w_{\vge,\scriptscriptstyle R,k,\scriptscriptstyle T}$  is increasing and we have that
$$\displaystyle w_{\vge,\scriptscriptstyle R}:=\lim_{k,\scriptscriptstyle T\to\infty}w_{\vge,\scriptscriptstyle R,k,\scriptscriptstyle T}\leq u\quad\text{in }\, \BBR^N\ti(\vge,+\infty).
$$
When $\vge\to 0$, $w_{\vge,\scriptscriptstyle R}(\cdot,\vge)$ converges to $\chi_{_{B_{\scriptscriptstyle R}}}\gm$. Hence its limit $w_{\scriptscriptstyle R}$ is expressed by 
 the Poisson kernel $(\ref{I-6})$, then 
$$
 w_{\scriptscriptstyle R}(x,z)=C_{{\scriptscriptstyle N},s}z^{2s}\myint{\BBR^N}{}\frac{\chi_{_{B_{\scriptscriptstyle R}}}(y)d\gm(y)}{(|x-y|^2+z^2)^{\frac{N}{2}+s}}\qquad\text{in }\BBR^N\ti(0,+\infty).
$$
As a consequence, letting $R\to+\infty$, we get
 $$
u(0,1)\geq C_{{\scriptscriptstyle N},s}\myint{B_R}{}\frac{d\gm(y)}{(|y|^2+1)^{\frac{N}{2}+s}},
$$
which implies $(\ref{I-T-3})$.

Finally, since $R\mapsto w_{\scriptscriptstyle R}$ is increasing and $w_{\scriptscriptstyle R}$  is upper bounded by $u$, we obtain that $w_{\scriptscriptstyle R}$ admits a limit $w$ when $R\to+\infty$. The function $w$ is a positive solution of $(\ref{I-T-1})$ in $\BBR^{N+1}_+$ which satisfies 
 \bel{I-T-9}\BA{lll}\displaystyle
 w(x,z)=C_{{\scriptscriptstyle N},s}z^{2s}\myint{\BBR^N}{}\frac{d\gm(y)}{(|x-y|^2+z^2)^{\frac{N}{2}+s}}\leq u(x,z)\quad\text{for } \, (x,z)\in\BBR^N\ti(0,+\infty).
\EA
\ee
Thus $w$ is the minimal positive solution of $(\ref{I-T-1})$ in $\BBR^{N+1}_+$ with trace $\gm$ on $\prt\BBR^{N+1}_+$. Set $\psi=u-w$, then 
$\psi$ is a nonnegative solution of $(\ref{I-T-1})$ in $\BBR^{N+1}_+$ with trace $0$ on $\prt\BBR^{N+1}_+$. The even extension of $\psi$ to whole 
$\BBR^{N+1}$ defined by \cite[Lemma 4.1]{CS} is a nonnegative solution of $(\ref{I-T-1})$  in $\BBR^{N+1}$. It is therefore continuous and it satisfies Harnack inequality \cite{FKS}. For proving uniqueness, we set  $t=\left(\frac{z}{2s}\right)^{2s}$ and $\tilde \psi(x,t)=\psi(x,z)$. Then $\tilde \phi$ satisfies (see\cite[1.8]{CS})
$$
t^{\frac{2s-1}{s}}\tilde \psi_{tt}+\Gd_x\tilde \psi=0.
$$
Since the operator $-\Gd_x$ is accretive in $L^1(\BBR^N)$, the function $t\mapsto \norm{\tilde \psi(\cdot,t)}_{L^1(\BBR^N)}$ is convex on $(0,\infty)$. 
If $\norm{\tilde \psi(\cdot,0)}_{L^1(\BBR^N)}=0$ and $\norm{\tilde \psi(\cdot,t)}_{L^1(\BBR^N)}=o(t)$ when $t\to\infty$, then $\norm{\tilde \psi(\cdot,t)}_{L^1(\BBR^N)}=0$ for all 
$t>0$. The condition $\norm{\tilde \psi(\cdot,t)}_{L^1(\BBR^N)}=o(t)$ is equivalent to $\norm{ \psi(\cdot,z)}_{L^1(\BBR^N)}=o(z^{2s})$. 

Therefore,  the uniqueness and $(\ref{I-T-9})$ implies $(\ref{I-T-4})$, which ends the proof.\qeda
\subsection{Spherical coordinates}

We recall spherical coordinates $(r,\gs)\in (0,+\infty)\ti \BBS^N$. The parametrization of $\BBR^{N+1}_+$ is the following
$$
\BBR^{N+1}_+=\big\{\xi=(r,\gs):r>0,\, \gs\in \BBS^N_+\big\},
$$
where
$$
\BBS^N_+=\left\{\gs=(x,z)=(\gs'\cos\gf ,\sin\gf):\gs'\in \BBS^{N-1},\gf\in \left[0,\frac\gp 2\right]\right\}.
$$
Without confusion, we also use the next notation 
$$\BBS^{N-1}= \partial \BBS^N_+=\big\{(0,\gs')\in\BBR^{N+1}: |\gs'|=1\big\}.$$%
With this parametrization 
$$
\Gd_{\BBS^N}u=u_{\phi\phi}-(N-1)(\tan\phi)\, u_\phi+\myfrac{1}{\cos^2\gf}\Gd_{\BBS^{N-1}}u.
$$
We define the operator $\CA_s$ in $C^2(\BBS^N_+)$ by 
$$
\CA_s[w]=\myfrac{1}{\gl_s(\phi)(\cos\phi)^{N-1}}\left(\gl_s(\phi)(\cos\phi)^{N-1}\;w_\phi\right)_\phi+\myfrac{1}{\cos^2\gf}\Gd_{\BBS^{N-1}}w,
$$
where 
$$
\gl_s(\phi)=(\sin\gf)^{1-2s}.
$$

We denote by $dS$  the invariant measure obtained by the isometric imbedding of $\BBS^N$ into $\BBR^{N+1}$, then
$$
dS(\gs)=(\cos\phi)^{N-1}dS'(\gs')d\gf,
$$
and by $\nabla'$ the covariant gradient in the canonical metric on $\BBS^{N}$ identified with the tangential gradient on the unit sphere of $\BBR^{N+1}$.
The bilinear form associated to $\CA_s$ is 
$$\BA{lll}\displaystyle
\CB[w,\tilde w]=\myint{\BBS^N_+}{}\left(w_\gf\tilde w_\gf+\frac{1}{\cos^2\gf}\nabla'w.\nabla'\tilde w\right)\gl_s(\gf) dS=-\myint{\BBS^N_+}{}\tilde w\CA_s(w)\gl_s(\gf)  dS,
\EA$$
and we have the corresponding Green's formula
$$
\myint{\BBS^N_+}{}\tilde w\CA_s(w)\gl_sdS=\int_{\BBS^{N-1}}\tilde w\prt_{\phi^s}\gw dS'-\CB[w,\tilde w].
$$
If $u=\BBP_s[v]$, then
$$
(-\Gd)^sv(r,\gs')=-\lim_{\gf\to 0}(\sin\gf)^{1-2s}\, u_\gf(r,\gs',\phi):=\prt_{\gf^s}u(r,\gs',0).
$$
\mysection{Self-Similar solutions }
Let $\Gw\subset\BBR^N$ be a bounded domain containing the origin. A function $v\in C^{2}(\Gw\setminus\{0\})\cap L^1_{\gm_s}(\BBR^N)$ satisfies the fractional Emden equations in $\Gw\setminus\{0\}$ with zero exterior value if
\bel{X1-3.1}\left.\BA{lll}\displaystyle
(-\Gd)^sv+\ge  v^p=0\ \  \text{in }\Gw\setminus\{0\},\\[2mm]
\phantom{(-\Gd)^s+\ge v^p}
 v=0\ \  \text{in }\Gw^c,
\EA\right.\ee
where  $p>1$ and $\ge=\pm 1$. If $\ge=1$ the equation considered is called {\it the fractional Emden-Fowler equation} while if $\ge=-1$ it is called {\it the fractional  Lane-Emden equation}. We denote 
$u=\BBP_s[v]$, which is admissible since $v\in L^1(\BBR^N)$, then 
 \bel{X-2}\left.\BA{lll}\displaystyle
\phantom{\ge z^{1-2s}----\,}\CD_s u=0\ \ &\text{in }\,\BBR^{N+1}_+\\[2mm]
\prt_{\gn^s}u(\cdot,0)+\ge  u(\cdot,0)^p=0\ \ &\text{in }\,\Gw\setminus\{0\}\\[2mm]
\phantom{\prt_{\gn^s}u(x ) |u|^{p-1}\,}
u(\cdot,0)=0\ \ &\text{in }\,\Gw^c.
\EA\right.
\ee
Hence $(\ref{X-2})$ implies that 
\bel{X-3}
\left.\BA{lll}\displaystyle
\phantom{\gth,\gth,\gth^pu \, }\;u_{rr}+\myfrac{N+1-2s}{r}u_r+\frac{1}{r^2}\CA_s[u]=0\ \
&\text{in }(0,\infty)\ti \BBS^N_+
\\[3mm]
\displaystyle
-\lim_{\gf\to 0}(\sin\gf)^{1-2s} u_\gf(r,\gth,\phi)+\ge  u(r,\gth,0)^p=0\ \ &\text{in }(0,\infty)\ti \BBS^{N-1}.
\EA\right.
\ee
Equation $(\ref{X-2})$ is equivariant under the family of transformations $\CS_\ell$, $\ell>0$, defined by
$$
\CS_\ell[u](x,z)=\ell^{\frac{2s}{p-1}}u(\ell x,\ell z)=\ell^{\frac{2s}{p-1}}u(\ell \xi).
$$
Therefore self-similar solutions of $(\ref{X-2})$ have the form 
$$
u(r,\gs)=r^{-\frac{2s}{p-1}} \gw(\gs)\quad\text{with } (r,\gs)\in (0,\infty)\ti \BBS^N_+,
$$
and $\gw$ satisfies 
\bel{X-5}\left.\BA{lll}\displaystyle
\CA_s[\gw]+ \Gl_{s,p,{\scriptscriptstyle N}}\gw=0\ \ &\text{in }\, \BBS^N_+\\[2mm]
\displaystyle\phantom{--\ \ }
\frac{\prt \gw}{\prt \gn^s}+\ge \gw^p=0\ \ &\text{in }\, \BBS^{N-1},
\EA\right. \ee
where 
$$
\Gl_{s,p,{\scriptscriptstyle N}}=\frac{2s}{p-1}\Big(\frac{2s}{p-1}+2s-N\Big).
$$
We denote by {\it $\CE_\ge$ (resp. $\CE_\ge^+$) the set of functions (resp. positive functions) $\gw\in C(\overline{\BBS^N_+})\cap C^2(\BBS^N_+)$  satisfying $(\ref{X-5})$.} There are several critical exponents will play some role later on:
\bel{X-7}
\BA{lll}\displaystyle
(i)\qquad\qquad &1<p\leq \frac{N}{N-2s}&\Longleftrightarrow\ \ \Gl_{s,p,{\scriptscriptstyle N}}\geq 0;\qquad\qquad\qquad\qquad\qquad\qquad\qquad
\\[4mm]\displaystyle
(ii)\qquad\qquad &\phantom{1<}p\geq 1+\frac{2s}{N}&\Longleftrightarrow \ \ \Gl_{s,p,{\scriptscriptstyle N}}\leq 2sN;
\\[4mm]\displaystyle
(iii)\qquad\qquad &\phantom{1<}p=\frac{N+2s}{N-2s}&\Longleftrightarrow\ \ \Gth_{s,p,{\scriptscriptstyle N}}:=N-2s\frac{p+1}{p-1}=0.
\EA\ee

\subsection{ Fractional Emden-Fowler equation }
\nind The structure of $\CE_1$ is described as follows.

\bth{Th1} Assume $s\in (0,1)$, $\ge=1$ and $p>1$.\smallskip

\nind 1- If $p\geq \frac{N}{N-2s}$, then $\CE_1=\{0\}$.  \\[1mm]
\nind 2- If $1<p\leq 1+\frac{2s}{N}$, then $\CE^+_1=\{0\}$.\\[1mm]
\nind 3- If $1+\frac{2s}{N}<p< \frac{N}{N-2s}$, then $\CE^+_1=\{0,\gw_1\}$, where $\gw_1$ is a positive solution of $(\ref{X-5})$ depending only on the variable $\gf$.  
\es
\Proof \nind {\it 1}- Let $\gw$ be a solution, then 

$$\int_{\BBS^N_+}\Big(\gw^2_\gf+\frac{1}{\cos^2\gf}|\nabla'w|^2-\Gl_{s,p,{\scriptscriptstyle N}}\gw^2\Big)\gl_s(\gf) dS=\int_{\BBS^{N-1}}\gw\prt_{\gf^s}\gw dS'=-\int_{\BBS^{N-1}}|\gw|^{p+1} dS'.
$$
If $p\geq \frac{N}{N-2s}$, then $\Gl_{{\scriptscriptstyle N},s,p}\leq 0$. This implies that $\gw=0$.\smallskip

\nind {\it 2}- It is easy to check that the first eigenfunction $\psi_1$ of $\CA_s$ in $W^{1,2}_0(\BBS^{N}_+)$ is $\gf\mapsto\psi_1(\gth,\gf)=(\sin\gf)^{2s}$ with corresponding eigenvalue $\ell_1=2sN$. Assume now that $\gw$ is a positive solution of $(\ref{X-5})$, then, multiplying the equation by $\psi_1$ and integrating yields
$$\left(\Gl_{s,p,{\scriptscriptstyle N}}-2sN\right)\int_{\BBS^N_+}\gw\psi_1\gl_s(\gf)dS=\int_{\BBS^{N-1}}\frac{\prt\psi_1}{\prt\gn^s}\gw dS'=-2s\int_{\BBS^{N-1}}\gw dS'.
$$
The claim follows by $(\ref{X-7})$-(ii).
\smallskip

\nind {\it 3- Existence.} We denote by $W(\BBS^{N}_+)$ the space of functions $w$ such that $\CB[w,w]<+\infty$, where 
$$\CB[w,v]=\int_{\BBS^N_+}\Big(w_\gf v_\gf  +\frac{1}{\cos^2\gf} \nabla' w\cdot \nabla' v \Big)\gl_s(\gf) dS.
$$
Then $w\mapsto \sqrt{\CB[w,w]}$ is a semi-norm on $W(\BBS^{N}_+)$ and it is a norm on the subspace $C^\infty_c(\BBS^N_+)$. 
Indeed 
$$\CB[w,w]\geq 2sN\int_{\BBS^N_+}w^2\gl_s(\gf) dS\qquad\text{for all }\;w\in C^\infty_c(\BBS^N_+).
$$
The closure of $C^\infty_c(\BBS^N_+)$ into
$W(\BBS^{N}_+)$ is denoted by $W_0(\BBS^{N}_+)$. We define the functional 
$$
J(w)=\frac 12\int_{\BBS^N_+}\Big(w_\gf^2+\frac{1}{\cos^2\gf}|\nabla' w|^2-\Gl_{s,p,{\scriptscriptstyle N}}w^2\Big)\gl_s(\gf) dS
+\frac{1}{p+1}\int_{\BBS^{N-1}}|\gg_0(w)|^{p+1}dS',
$$
where $\gg_0$ denotes the trace operator from $W(\BBS^{N}_+)$ onto $L^2(\BBS^{N-1})$ identified with $w(0,\gth)$. Then $J$ is a quadratic perturbation of a convex lower semicontinous functional defined in $W(\BBS^{N}_+)$, hence it is weakly lower continuous with domain of definition
$$D(J)=\{w\in W(\BBS^{N}_+): \gg_0(w)\in L^{p+1}(\BBS^{N-1})\}.$$ 
Let $w\in W(\BBS^{N}_+)$ and
$$
w=w_1+w_2,
$$
where $w_1\in W_0(\BBS^{N}_+)$ is defined by 
$$
\CB[w_1,\gz]=\CB[w,\gz]\qquad\text{for all }\;\gz\in W_0(\BBS^{N}_+).
$$
Since $\CB$ is a continuous quadratic form in $W(\BBS^{N}_+)$, $\gz\mapsto \CB[w,\gz]$ is a continuous linear form in $W(\BBS^{N}_+)$ and there holds
$$
\CB[\gz,\gz]\geq \left(2sN-\Gl_{s,p,{\scriptscriptstyle N}}\right)\int_{\BBS^N_+}\gz^2\gl_s(\gf)  dS\qquad\text{for all }\;\gz\in W_0(\BBS^{N}_+),
$$
 the function $w_1$ is uniquelly defined. Thus, one has 
 $$J(w)=J(w_1+w_2)=\frac 12\CB[w_1,w_1]+\frac 12\CB[w_2,w_2]+2\CB[w_1,w_2]+\frac{1}{p+1}\int_{\BBS^{N-1}}|\gg_0(w)|^{p+1}dS'.
 $$
We have $\CB[w_1,w_1]=\CB[w,w_1]=\CB[w_1,w_1]+\CB[w_2,w_1]$, therefore
$$\CB[w_2,w_1]=\CB[w_1,w_2]=0.
$$
Furthermore $w_2(\gf,\gth)=w_2(\gf,\gth)-w_2(0,\gth)+w_2(0,\gth)$. Since for any $\vge>0$, 
$$w^2_2(\gf,\gth)\leq (1+\vge)\left(w_2(\gf,\gth)-w_2(0,\gth)\right)^2+\left(1+\frac 1\vge\right)w^2_2(0,\gth),
$$
and $(\gf,\gth)\mapsto w_2(\gf,\gth)-w_2(0,\gth)$ belongs to $W^{1,2}_0(\BBS^N_+)$, we have the inequality

$$\BA{lll}\displaystyle\CB[w_2,w_2]\geq \int_{\BBS^N_+}\frac{1}{\cos^2\gf}|\nabla w_2|^2\gl_s(\gf) dS+\int_{\BBS^N_+}\left(\left(w_2(\gf,\gth)-w_2(0,\gth)\right)_\gf\right)^2\gl_s(\gf) dS\\[4mm]
\phantom{\CB[w_2,w_2]--\quad }\displaystyle-(1+\vge)\Gl_{s,p,{\scriptscriptstyle N}}\int_{\BBS^N_+}\left(w_2(\gf,\gth)-w_2(0,\gth)\right)^2\gl_s(\gf) dS
\\[4mm]
\phantom{\CB[w_2,w_2]\quad ----}\displaystyle-\Big(1+\frac 1\vge\Big)\Gl_{s,p,{\scriptscriptstyle N}}\int_{\BBS^N_+}w^2_2(0,\gth)\gl_s(\gf) dS. 
\EA$$
We choose $\vge$ such that 
$$(1+\vge)\Gl_{s,p,{\scriptscriptstyle N}}<2sN
$$
and we get
$$\BA{lll}\displaystyle\CB[w_2,w_2]\geq \Big(1-\frac{(1+\vge)\Gl_{s,p,{\scriptscriptstyle N}}}{2sN}\Big)\int_{\BBS^N_+}\Big(w_{2\,\gf}^2+\frac{1}{\cos^2\gf}|\nabla' w_2|^2\Big)\gl_s(\gf) dS\\[4mm]\phantom{\CB[w_2,w_2]--}\displaystyle
-\Big(1+\frac 1\vge\Big)\Gl_{s,p,{\scriptscriptstyle N}}\int_{\BBS^N_+}w^2_2(0,\gth)\gl_s(\gf) dS.
\EA$$
Notice that 
$$\BA{lll}\displaystyle\int_{\BBS^N_+}w^2_2(0,\gth)\gl_sdS=\int_0^{\frac\gp 2}(\sin\gf)^{1-2s}(\cos\gf)^{N-1}d\gf\int_{\BBS^{N-1}}w^2_2(0,\gth)dS'\\[4mm]
\phantom{\displaystyle\int_{\BBS^N_+}w^2_2(0,\gth)\gl_sdS}\displaystyle
\leq c_1\Big(\int_{\BBS^{N-1}}|w_2(0,\gth)|^{p+1}dS'\Big)^{\frac{2}{p+1}},
\EA$$
where $c_1=c_1(s,p,N)>0$. Because $w_1\in W^{1,2}_0(\BBS^N_+)$, there holds  also
$$\CB[w_1,w_1]\geq \Big(1-\frac{\Gl_{s,p,{\scriptscriptstyle N}}}{2sN}\Big)\int_{\BBS^N_+}\Big(w_{1\,\gf}^2+\frac{1}{\cos^2\gf}|\nabla' w_1|^2\Big)\gl_s(\gf) dS.
$$
Combining the previous inequalities, we finally obtain
$$
\BA{lll}\displaystyle J(w)\geq \Big(1-\frac{(1+\vge)\Gl_{s,p,{\scriptscriptstyle N}}}{2sN}\Big)\int_{\BBS^N_+}\Big(w_{\gf}^2+\frac{1}{\cos^2\gf}|\nabla' w|^2\Big)\gl_s(\gf) dS\\
\phantom{---\,}\displaystyle -\Big(1+\frac 1\vge\Big)M\Gl_{s,p,{\scriptscriptstyle N}}
\Big(\int_{\BBS^{N-1}}|\gg_0w|^{p+1}dS'\Big)^{\frac{2}{p+1}}+\frac{2}{p+1}\int_{\BBS^{N-1}}|\gg_0w|^{p+1}dS'.
\EA$$
This implies that $J(w)$ tends to $\infty$ when $\min\Big\{\sqrt{\CB[w,w]},\norm {\gg_0(w)}_{L^{p+1}(\BBS^{N-1})}\Big\}$ tends to infinity. Consequently the functional 
$J$ admits a minimum $\gw$ in $W(\BBS^N_+)$. Moreover, since $J(|\gw|)=J(\gw)$ the minimum is achieved by a nonnegative function and it 
 is therefore classical that $\gw$ is a solution of $(\ref{X-5})$, and it is positive by the strong maximum principle. \medskip
 
 Finally, if we consider the restriction $J_{rad}$ of $J$ to the space $W_{rad}(\BBS^N_+)$ of function of $W(\BBS^N_+)$ depending only on the variable $\gf$
(they are called radial), then $\gg_0(w^{p+1})$ is a real number and  $J_{rad}(w)$ has the following form provided we write $w=w(\gf)$,
$$
J_{rad}(w)=\frac {|\BBS^{N-1}|}2\int_{0}^\frac{\gp}{2}\left(w_\gf^2-\Gl_{s,p,{\scriptscriptstyle N}}w^2\right)\gl_s(\gf)(\cos\gf)^{N-1}d\gf
+\frac{|\BBS^{N-1}|}{p+1}|w(0)|^{p+1}.
$$
 The functional $J_{rad}$ is also a quadratic perturbation of a convex lower semicontinuous functional defined in $W_{rad}(\BBS^N_+)$ and it tends to $\infty$ at $\infty$. Thus it admits a minimum which is achieved by a nonnegative function due to $J_{rad}(w)=J_{rad}(|w|)$. {\it Hence there exists a minimizing nonnegative solution of $(\ref{X-5})$ which depends only on the variable $\gf$}.
\smallskip

\nind {\it 3- Uniqueness. } Let $\gw$ and $\tilde\gw$ be two positive solutions of $(\ref{X-5})$, then

$$\BA{lll}\displaystyle 0=\myint{\BBS^N_+}{}\Big(\frac{\CA_s[\gw]}{\gw}-\frac{\CA_s[\tilde\gw]}{\tilde\gw}\Big)(\gw^2-\tilde\gw^2)\gl_s dS\\[4mm]
\displaystyle \phantom{0}=\int_{\BBS^{N-1}}\left(\prt_{\gf^s}\gw\Big(\gw-\frac{\tilde\gw^2}{\gw}\Big)-\prt_{\gf^s}\tilde\gw\Big(\tilde\gw-\frac{\gw^2}{\tilde\gw}\Big) \right) dS'-\left(\CB\Big[\gw,\gw-\frac{\tilde\gw^2}{\gw}\Big]-\CB\Big[\tilde\gw,\tilde\gw-\frac{\gw^2}{\tilde\gw}\Big]\right)\\[4mm]
\displaystyle \phantom{0}=A-B.
\EA$$
Furthermore
$$\BA{lll}\displaystyle A=\int_{\BBS^{N-1}}\left(-\gw^p\Big(\gw-\frac{\tilde\gw^2}{\gw}\Big)+\tilde\gw^p\Big(\tilde\gw-\frac{\gw^2}{\tilde\gw}\Big) \right) dS'
\\[4mm]
\displaystyle \phantom{A}=-\int_{\BBS^{N-1}}\left( \gw^2-\tilde\gw^2\right) \left(\gw^{p-1}-\tilde\gw^{p-1} \right) dS'\leq 0 
\EA$$
and
$$\BA{lll}\displaystyle B=
\int_{\BBS^{N-1}}\Big(\frac{1}{\gw^2}+\frac{1}{\tilde\gw^2}\Big)\Big(\big(\gw\tilde\gw_\gf-\tilde\gw\gw_\gf\big)^2+\frac{1}{\cos^2\gf}\left|\gw\nabla'\tilde\gw-\tilde\gw\nabla'\gw\right|^2\Big)\gl_s(\gf) dS'\geq 0.
\EA$$
Therefore, $A=B=0$ and $\gw=\tilde\gw$. \qeda
\medskip

\nind\BBRemark To the unique positive element $\gw_1$ of $\CE^+_1$ corresponds a unique positive singular self-similar solution $U_p$ of 
$$
(-\Gd)^s v+v^p=0\quad\text{in }\,\BBR^N\setminus\{0\},
$$
where $U_p(x)=c_p|x|^{-\frac{2s}{p-1}}$ with 
$$
 c_p=\gw_1(1)=\left(-2^{2s}
\myfrac{\Gg\big(\frac N2-\frac{s}{p-1}\big)\Gg\big(s+\frac{s}{p-1}\big)}
{\Gg\big(\frac{s}{p-1}\big)\Gg\big(\frac N2-s-\frac{s}{p-1}\big)}\right)^\frac{1}{p-1}.
$$

In order to prove \rth{Th2}, we need the following intermediate result.
\bprop{rad} Assume $\Gl\neq 0$, then for any $a\neq 0$, there exists a unique function $\gw_a$   satisfying 
 \bel{V3-7}\BA{lll}\displaystyle
\gw_a(\gf)=a-\Gl\myint{\gf}{\frac\gp2}(\sin\gs)^{2s-1}(\cos\gs)^{1-N}\myint{\gs}{\frac\gp2}\gw_a(\gth)(\sin\gth)^{1-2s}(\cos\gth)^{N-1} d\gth d\gs\ \ \text{in }\big(0,\frac\gp 2\big).
 \EA\ee
Furthermore $a=\gw_a(\frac\gp2)$ and $\gw_a=a\gw_1$. If $a>0$, then $\gw_a$ is positive and increasing (resp. decreasing) on $\left(0,\frac\gp 2\right)$ if $\Gl<0$ (resp. $\Gl>0$).
\es
\Proof  Let $a>0$ and $\tilde\gw_a(\gf)=\gw_a(\frac\gp 2-\gf)$. To find a function $\gw_a$ satisfying $(\ref{V3-7})$ is equivalent to finding  $\tilde \gw_a$ satisfying
$$
\tilde \gw_a(\gf)=a-\Gl\myint{0}{\gf}(\cos\gs)^{2s-1}(\sin\gs)^{1-N}\myint{0}{\gs}\tilde\gw_a(\gth)(\cos\gth)^{1-2s}(\sin\gth)^{N-1}  d\gth d\gs\quad\text{in }\big(0,\frac\gp 2\big).
$$
We define the operator $\CT$ on $C([0,\frac\gp 2])$ by 
$$
\CT[w](\gf)=a-\Gl\myint{0}{\gf}(\cos\gs)^{2s-1}(\sin\gs)^{1-N}\myint{0}{\gs}w(\gth)(\cos\gth)^{1-2s}(\sin\gth)^{N-1}   d\gth d\gs\quad\text{in }\big(0,\frac\gp 2\big).
$$
We have that
$$\BA{lll}\displaystyle|\CT[w-w'](\gf)|\\[4mm]\phantom{--}
\displaystyle\leq |\Gl|\left(\myint{0}{\gf}(\cos\gs)^{2s-1}(\sin\gs)^{1-N}\myint{0}{\gs}(\cos\gth)^{1-2s}(\sin\gth)^{N-1}    d\gth d\gs\right)\sup_{[0,\gf]}|(w-w')(\gth)|.
\EA$$
1- We first assume that $1-2s\geq 0$, then
$$\BA{lll}\displaystyle
\quad \myint{0}{\gf}(\cos\gs)^{2s-1}(\sin\gs)^{1-N}\myint{0}{\gs}(\cos\gth)^{1-2s}(\sin\gth)^{N-1} d\gth d\gs 
\leq 
\myint{0}{\gf}(\cos\gs)^{2s-1}  \gs d\gs\\[4mm]
\phantom{-------------------- }
\displaystyle\leq  
\frac{\pi}2\myint{0}{\gf}(\cos\gs)^{2s-1} \sin \gs d\gs
\leq   c_2,
\EA$$
where $c_2=\frac{\pi}{4s}$,  we used here the monotonicity $\gth\in(0,\frac{\pi}{2})\mapsto  (\sin\gth)^{N-1}$ and $(\cos\gth)^{1-2s}\leq 1$.
Hence
 $$
|\CT[w-w'](\gf)|\leq c_2\sup_{[0,\gf]}|(w-w')(\gth).
 $$
It is therefore straightforward to check by induction that for $k\in\BBN^*$,
 \bel{V3-11}\BA{lll}\displaystyle
|\CT^k[w-w'](\gf)|\leq \myfrac{c_2^k}{k!}\sup_{[0,\gf]}|(w-w')(\gth),
\EA\ee
then $\CT^k$ admits a unique fixed point $\tilde\gw_a$, and uniqueness implies that $\tilde\gw_a$ is also a fixed point of $\CT$ and it is unique.\smallskip

\nind 2- Assume that $1-2s< 0$. Then
$$\BA{lll}\displaystyle
  \myint{0}{\gf}(\cos\gs)^{2s-1}(\sin\gs)^{1-N}\myint{0}{\gs}(\cos\gth)^{1-2s}(\sin\gth)^{N-1} d\gth d\gs 
\\[4mm]
\phantom{ ----------- }
\displaystyle \leq 
\myint{0}{\gf} \myint{0}{\gs}  d\gth d\gs=\frac{\gf^2}{2}\leq \frac{\pi^2}{8}
\EA$$
by there the monotonicity $\gth\in(0,\frac{\pi}{2})\mapsto (\cos\gth)^{1-2s}(\sin\gth)^{N-1}$. 
Then $(\ref{V3-11})$ holds and the conclusion follows.\smallskip

If $\Gl<0$ (resp. $\Gl<0$), the function $\gw_a$ is increasing (resp. decreasing).
The fact that $\gw_a=a\gw_1$ is a consequence of the linearity of  $\CT$ and uniqueness.
\qeda\medskip

The following critical exponent corresponds to the sequel of isotropy  of solutions of  $(\ref{X-5})$ in the case $\ge=1$,
 \bel{XX1}
p^*:=\frac{N+2s+\sqrt{N^2+4(N-s)+4s^2-4}}{N-2s+\sqrt{N^2+4(N-s)+4s^2-4}}.
 \ee
This exponent corresponds to the fact that $\frac{2s}{p^*-1}$ is the positive root of the equation
\bel{XX2}
Q(X):=X^2+(2s-N)X+1-N=0.
\ee

\bth{Th1*} Assume $s\in (0,1)$, $\ge=1$ and $p^*\leq p<\frac{N}{N-2s}$. Then $\CE_1=\{0,\gw_1,-\gw_1\}$, where $\gw_1$ is a positive solution of $(\ref{X-5})$ depending only on the variable $\gf$.
\es
\Proof {\it Step 1.} We first prove that under the condition $p\geq p^*$ any element of $\CE_1$ does not depend on the variable $\gs'\in \BBS^{N-1}$. If $\psi$ is any function defined on $\BBS^N$, we set 
$$\overline\psi(\gf)=\frac{1}{|\BBS^{N-1}|}\int_{\BBS^{N-1}}\psi(\gf,\gs')dS'.
$$
By averaging $(\ref{X-5})$ we obtained that $\bar\gw$ satisfies
$$
\left.\BA{lll}\displaystyle
\CA_s[\overline\gw]+\Gl_{s,p,{\scriptscriptstyle N}}\overline\gw=0\ \ &\text{in }\, \BBS^N_+\\[2mm]
\displaystyle\phantom{--\ \ \,\, }
\frac{\prt \overline\gw}{\prt \gn^s}+\overline{\gw^p}=0\ \ &\text{in }\, \BBS^{N-1},
\EA\right.$$
where  
$$\CA_s[\overline\gw]=\frac{1}{\gl_s(\gf)(\cos\gf)^{N-1}}\left(\gl_s(\gf)(\cos\gf)^{N-1}\overline\gw\right)_\gf.
$$
Multiplying the equation by $(\gw-\overline\gw)\gl_s$ and integrating over $\BBS^N_+$, we obtain
 \bel{XXX6}\BA{lll}\displaystyle\displaystyle\int_{\BBS^{N-1}}(\gw-\overline\gw)(|\gw|^{p-1}\gw-\overline{|\gw|^{p-1}\gw})dS'+\int_0^{\frac\gp2}\int_{\BBS^{N-1}}
\Big((\gw_\gf-\overline\gw_\gf)^2+\frac{1}{\cos^2\gf}|\nabla' (\gw-\overline\gw)|^2\Big)\gl_s(\gf)dS\\[4mm]
\phantom{---------------------}\displaystyle=\Gl_{s,p,{\scriptscriptstyle N}}\int_{\BBS^N_+}(\gw-\overline\gw)^2\gl_sdS.
\EA\ee
There holds
$$\BA{lll}\displaystyle
\int_{\BBS^{N-1}}(\gw-\overline\gw)\big( \gw^p-\overline{\gw^p}\big)dS'=\int_{\BBS^{N-1}}(\gw-\overline\gw)(\gw^p- \overline{\gw}^p)dS'
  +(\overline{\gw}^p-\overline{\gw^p} ) \int_{\BBS^{N-1}}(\gw-\overline\gw)dS'\\[4mm]\displaystyle
\phantom{\int_{\BBS^{N-1}}(\gw-\overline\gw)\big( \gw^p-\overline{\gw^p}\big)dS'}
= \int_{\BBS^{N-1}}(\gw-\overline\gw)\big( \gw^p-\overline{\gw}^p\big)dS'\\[4mm]\displaystyle
\phantom{\int_{\BBS^{N-1}}(\gw-\overline\gw)\big( \gw^p-\overline{\gw^p}\big)dS'}
\geq 2^{-p}\int_{\BBS^{N-1}}|\gw-\overline\gw|^{p+1}dS'.
\EA$$
 Furthermore, we see that 
 \bel{XXX7}\BA{lll}\displaystyle
 \int_0^{\frac\gp2}\int_{\BBS^{N-1}}
\Big((\gw_\gf-\overline\gw_\gf)^2+\frac{1}{\cos^2\gf}|\nabla' (\gw-\overline\gw)|^2\Big)\gl_s(\gf)dS\\[4mm]
\phantom{------}
\displaystyle\geq 
 \int_0^{\frac\gp2}\Big(\int_{\BBS^{N-1}}|\nabla' (\gw-\overline\gw)|^2dS'\Big)(\cos\gf)^{N-3}\gl_s(\gf)d\gf
 \EA\ee
 and the inequality is strict unless $\gw_\gf=\overline\gw_\gf$. By Wirtinger's inequality,  there holds
 $$\int_{\BBS^{N-1}}|\nabla' (\gw-\overline\gw)|^2dS'\geq (N-1)\int_{\BBS^{N-1}}(\gw-\overline\gw)^2dS'.
 $$
 Since $(\cos\gf)^{N-3}\geq (\cos\gf)^{N-1}$, we obtain
 $$\int_0^{\frac\gp2}\int_{\BBS^{N-1}}
\Big((\gw_\gf-\overline\gw_\gf)^2+\frac{1}{\cos^2\gf}|\nabla' (\gw-\overline\gw)|^2\Big)\gl_s(\gf)dS\geq (N-1)\int_{\BBS^N_+}(\gw-\overline\gw)^2\gl_s(\gf)dS. 
 $$
 Finally we have that 
 $$
\left( \Gl_{s,p,{\scriptscriptstyle N}}+1-N\right)\int_{\BBS^N_+}(\gw-\overline\gw)^2\gl_s(\gf)dS\geq 2^{-p}\int_{\BBS^{N-1}}|\gw-\overline\gw|^{p+1}dS'.
 $$
Introducing the polynomial $Q$ defined in $(\ref{XX2})$ and setting $X=\frac{2s}{p-1}$, then  $(\ref{XXX6})$ becomes
 $$Q(X)\int_{\BBS^N_+}(\gw-\overline\gw)^2\gl_s(\gf)dS\geq 2^{-p}\int_{\BBS^{N-1}}|\gw-\overline\gw|^{p+1}dS'.
 $$
 Hence if $X< X^*$, where $X^*$ is the positive root of $Q$, then $\gw-\overline\gw=0$ in  $\BBS^N_+$. If $X= X^*$ we deduce that $\gw(0,\gth)=\overline\gw(0)=0$. Then we use that fact that the inequality $(\ref{XXX7})$ is strict unless $\gw_\gf=\overline\gw_\gf$. If this holds, combining with  $\gw(\gth)-\overline\gw(0)=0$
 we deduce that $\gw=\overline\gw$ in  $\BBS^N_+$.\smallskip
 
 \nind {\it Step 2.}  We prove that $\CE_1=\left\{0,\gw_1,-\gw_1\right\}$, where $\gw_1$ is a positive solution of $(\ref{X-5})$ with $\ge=1$, depending only on $\gf$. For this sake, it   suffices to prove uniqueness among the radial solutions. If $\gf\mapsto\gw(\gf)$ is a solution of $(\ref{X-5})$ with $\ge=1$, we have
  \bel{V3-1*}\left.\BA{lll}\displaystyle
-\frac{d}{d\gf}\Big((\sin\gf)^{1-2s}(\cos\gf)^{N-1}\frac{d\gw}{d\gf}\Big)=\Gl_{s,p,{\scriptscriptstyle N}}(\sin\gf)^{1-2s}(\cos\gf)^{N-1}\gw\quad\text{in }\big(0,\frac\gp 2\big)\\[4mm]
\phantom{--------   }
\myfrac{d\gw(0)}{d\phi^s}+ \gw(0)^p=0,
\EA\right.\ee
where
$$
\myfrac{d\gw(0)}{d\phi^s}=-\lim_{\gf\to 0}(\sin\gf)^{1-2s}\myfrac{d\gw(\gf)}{d\phi}.
$$
Since $\gw'(\tfrac\gp2)=0$,  by integrating $(\ref{V3-1*})$ we obtain that 
\bel{V3-5*}\BA{lll}\displaystyle
\frac{d\gw(\gf)}{d\gf}=\Gl_{s,p,{\scriptscriptstyle N}}(\sin\gf)^{2s-1}(\cos\gf)^{1-N}\myint{\gf}{\frac\gp2}\gw(\gth)(\sin\gth)^{1-2s}(\cos\gth)^{N-1} d\gth
\EA\ee
on $\left(0,\frac\gp 2\right)$. Therefore for $\gf\in \left(0,\frac\gp 2\right)$,
\bel{V3-6*}\BA{lll}\displaystyle
\gw(\gf)=\gw(\tfrac{\gp}{2})-\Gl_{s,p,{\scriptscriptstyle N}}\myint{\gf}{\frac\gp2}(\sin\gs)^{2s-1}(\cos\gs)^{1-N}\myint{\gs}{\frac\gp2}\gw(\gth)(\sin\gth)^{1-2s}(\cos\gth)^{N-1} d\gth d\gs.
\EA\ee 
From $(\ref{V3-5*})$, 
\bel{V3-7*}
\myfrac{d\gw}{d\phi^s}(0)=-\Gl_{s,p,{\scriptscriptstyle N}}\myint{0}{\frac\gp2}\gw(\gth)(\sin\gth)^{1-2s}(\cos\gth)^{N-1} d\gth.
\ee
Note  that for $p<\frac{N}{N-2s}$, $\Gl_{s,p,{\scriptscriptstyle N}}>0$. As a parameter we take the value $a=\gw(\tfrac{\gp}{2})$, and by linearity $\gw_a=a\gw_1$. We define the 
mapping $F^*:\BBR_+\mapsto\BBR$ by 
$$F^*(a)=\frac {d\gw_a}{d\gf^s}(0)+(\gw_a(0))^{p}=a\Big(\frac {d\gw_1}{d\gf^s}(0)+a^{p-1}(\gw_1(0))^{p}\Big).
$$
If  $\frac {d\gw_1}{d\gf^s}(0)\geq 0$, then for any $a>0$ we have $\frac 1aF^*(a)>0$. Since $F^*$  is increasing, there is no solution of $(\ref{V3-1*})$, then no radial solution of $(\ref{X-5})$ with $\ge=1$, which contradicts \rth{Th1}, which states that there exists a radial solution, even a positive one. 
Therefore $\frac {d\gw_1}{d\gf^s}(0)< 0$. The mapping $a\mapsto \frac 1aF^*(a)$ is increasing and negative at $a=0$. Since  $F^*$ tends to infinity as $a\to+\infty$, it follows that there exists a unique 
$a^*>0$ such that $F^*(a^*)=0$, equivalently a unique solution $\gw:=\gw_{a^*}$ of $(\ref{V3-1*})$ verifying $a^*=\gw_{a^*}(\tfrac{\gp}{2})$. \smallskip

We end the proof as follows. From {\it Step 1} the solutions $\gw$ depends only on $\gf$. Such a solution satisfying $\gw(\frac{\gp}{2})>0$ is unique, thus it coincides  with the radial positive solution obtained by minimization among radial functions.\qeda

\subsection{ Fractional Lane-Emden equation}
\bth{Th2} Assume $s\in (0,1)$, $\ge=-1$ and $p>1$.\smallskip

\nind 1- If $p\leq \frac{N}{N-2s}$, then $\CE^+_{-1}=\{0\}$.\\[1mm]
\nind 2- If $p> \frac{N}{N-2s}$, then $\CE^+_{-1}\!=\{0,\gw_2\}$,  where $\gw_2$ is a positive solution of $(\ref{X-5})$ depending only on $\phi$.
\es
\Proof {\it 1-}  A direct computation implies that  
$$\Gl_{s,p,{\scriptscriptstyle N}}\int_{\BBS^N_+}\gw\gl_sdS+\int_{\BBS^N_+}\prt_{\gf^s}\gw dS'=\Gl_{s,p,{\scriptscriptstyle N}}\int_{\BBS^N_+}\gw\gl_sdS+\int_{\BBS^N_+}\gw^p dS'=0.
$$
Since $1<p\leq \frac{N}{N-2s}$, $\Gl_{s,p,{\scriptscriptstyle N}}\geq 0$ and the assertion 1 follows.\smallskip

\nind {\it 2- Existence of a unique positive solution depending only on $\gf$}. If $\gw$  is a solution of $(\ref{X-5})$ with $\ge=-1$ depending only on $\gf$, then
\bel{V3-1}
\left.\BA{lll}\displaystyle
-\frac{d}{d\gf}\big((\sin\gf)^{1-2s}(\cos\gf)^{N-1} \gw_{\gf} \big)=\Gl_{s,p,{\scriptscriptstyle N}}(\sin\gf)^{1-2s}(\cos\gf)^{N-1}\gw\quad\text{in }\Big(0,\frac\gp 2\Big)\\[4mm]
\phantom{------(\cos\gf)^{N-}\,}
\myfrac{d\gw(0)}{d\phi^s}= \gw(0)^{p}.
\EA\right.\ee
Since $\gw'(\tfrac\gp2)=0$, the identities $(\ref{V3-5*})$, $(\ref{V3-6*})$, and $(\ref{V3-7*})$ obtained by  integrating $(\ref{V3-1})$ are still valid. Because $p>\frac{N}{N-2s}$, we have that $\Gl_{s,p,{\scriptscriptstyle N}}<0$, hence $\frac{d\gw}{d\phi^s}(0)>0$ from  identity  $(\ref{V3-7*})$.
By linearity, $\gw_a=a\gw_1$. We define the mapping $F$ from $[0,\infty)$ to $(-\infty,\infty)$ by 
$$F(a)=\frac {d\gw_a(0)}{d\gf^s}- \gw_a(0)^{p}=a\Big(\frac {d\gw_1(0)}{d\gf^s}-a^{p-1} \gw_1(0)^{p}\Big).
$$
 This implies that $a\mapsto \frac{1}{a}F(a)$ is decreasing, with a positive limit at $a=0$, and tends to $-\infty$ when $a\to+\infty$. Therefore, $F$ admits a unique 
positive zero point for some $a_0$, and the function $\gw_{a_0}$ is the unique positive solution of  $(\ref{X-5})$ which depends only on $\gf$.\smallskip

\nind {\it 3- Any positive solution depends only on $\gf$}.\smallskip

\nind {\it Symmetry in the case $s=\frac 12$}. In that case $s=\frac 12$, the problem $(\ref{X-5})$ reduces to
\bel{X-5-1}\left.\BA{lll}\displaystyle
\Gd_{\BBS^N}\gw+\Gl_{\frac 12,p,{\scriptscriptstyle N}}\gw=0\ \ &\text{in }\BBS^N_+\\[2mm]
\displaystyle\phantom{---\  \, }
\frac{\prt \gw}{\prt \gn}-\gw^{p}=0\ \ &\text{in }\BBS^{N-1},
\EA\right.\ee
where $\Gd_{\BBS^N}$ is the Laplace-Beltrami operator on $\BBS^N$ and 
$$\Gl_{\frac 12,p,{\scriptscriptstyle N}}=\frac{1}{p-1}\Big(\frac{1}{p-1}+1-N\Big).
$$
We denote the coordinates in $\BBR^N$ by $(x_1,\cdots,x_{\scriptscriptstyle N})$ and by $(x_1,\cdots,x_{\scriptscriptstyle N},z)$, $z>0$, the coordinates in $\BBR^{N+1}_+$.  If $\gf_0\in (0,\frac\gp2)$, let $\BBH_{\gf_0}$ be the hyperplane  of $\BBR^{N+1}$ passing through the (N-1)-plane $D_{x_1}:=\big\{(0,x_2,\cdots,x_{\scriptscriptstyle N},0)\in\BBR^{N+1}:x_j\in\BBR,\,j=2,\cdots,N\big\}$ with angle $\gf_0$ with the hyperplane $z=0$. We denote by $\Gg_{\phi_0}$ the domain of $\BBS^N$ with boundary $\BBH_{\gf_0}\cap \prt\BBR^{N+1}_+\cap\big\{(x_1,\cdots,x_{\scriptscriptstyle N},z):x_1\geq 0\big\}$ and by $\tilde \Gg_{\phi_0}$ the symmetric of this domain with respect to the plane $\BBH_{\gf_0}$. If $\gw$ is a positive solution of  $(\ref{X-5})$, we define the function
\bel{V-5-2}
(\gf,\gs')\mapsto \varpi(\gf,\gs'):=\gw(\gf,\gs')-\gw(2\gf_0-\gf,\gs')\qquad\text{for }(\gf,\gs')\in (0,\gf_0)\ti \BBS^{N-1}.
\ee
Since $\gf$ can be larger than $\frac{\gp}{2}$, we make the convention that $\gw(\gf,\gs')=\gw(\gp-\gf,-\gs')$ if $\frac\gp2<\gf<\gp$. 
The symmetry with respect to $\BBH_{\gf_0}$ is an isometry of $\BBR^{N+1}$, hence $\varpi$ satisfies 
$$
\left.\BA{lll}\displaystyle
 \Gd_{\BBS^N}\varpi+\Gl_{\frac12,p,{\scriptscriptstyle N}}\varpi=0\quad&\text{in }\,\Gg_{\phi_0}\\[2mm]
 \phantom{\frac{\prt v}{\prt \gn}+\ge|v|^{p-1}{\BBS^N} }
\varpi=0\quad&\text{in }\, \BBH_{\gf_0}\cap \BBS^{N}_+\\[2mm]
 \phantom{-,--{\BBS^N}}
\!\varpi(0,\cdot)=\gw(0,\cdot)-\gw(2\gf_0,\cdot)\ \ &\text{in }\, \BBS^{{N-1}}\cap\{ (x_1,x_2,...,x_N,0):x_1>0\}.
\EA\right.$$
Since $\frac{\prt\gw}{\prt\gn^s}=\gw^p$ on $\prt \BBS^N_+$, the function $\gf\mapsto\gw(\gf,\gs')$ is strictly decreasing if $(\gf,\gs')\in (0,\gf_0)\ti \BBS^{N-1}$ for $\gf_0$ small enough. Hence $\varpi(0,\gs')>0$ in $\BBS^{N-1}\cap\{x=(x_1,x_2,...,x_N,0):x_1>0\}$, and by the maximum principle we have that $\varpi>0$ in $(0,\gf_0)\ti \BBS^{N-1}$. We proceed as in the Gidas-Ni-Nirenberg paper \cite{GNN}: we denote by $\gf^*$ the maximum of the $\gf_0\in (0,\frac\gp2]$ such that 
\bel{V-5-3-0}\gw(\gf,\gs')-\gw(2\gf_0-\gf,\gs')>0\qquad\text{for }\, (\gf,\gs')\in (0,\gf_0)\ti \BBS^{N-1}.
\ee
We assume by contradiction that $\gf^*<\frac\gp 2$. Then there exists $\gs^*\in \BBS^{N-1}$ such that $v(\gf^*,\gs^*)=0$. By the Hopf's Lemma, $(\gf^*,\gs^*)$ cannot belong to $\BBS^{N-1}_+$, hence it belongs to $\overline{\BBS^{N}_+}\cap \{z=0\}$. Since Hopf's lemma is also valid at this point, we get a contradiction. Therefore 
$\gf^*=\frac\gp2$. As a consequence, 
$$\gw(0,\gs')-\gw(2\gf^*-\gf,\gs')\geq 0\quad\text{for all }\,\gs'\in \BBS^{N-1}\cap\{(0,x_2,...,x_N,z):z>0\}. 
$$
Similarly, starting the reflexion from $x_1<0$, we get
$$\gw(0,\gs')-\gw(2\gf^*-\gf,\gs')\leq 0\quad\text{for all }\,\gs'\in \BBS^{N-1}\cap\{(0,x_2,...,x_N,z):z>0\}. 
$$
This implies that the gradient of $\gw$ alongside the ``great circle" $\BBS^N_+\cap \{(0,x_2,...,x_N,z):z>0\}$ is zero, which can be written as 
$$\nabla'\gw\cdot {\bf e}_1=0\quad\text{on }\,\BBS^N_+\cap \big\{(0,x_2,...,x_N,z):z>0\big\}.
$$
Performing a rotation we conclude that $\nabla'(\gw\gf,\gs')=0$ for all $\gs'\in \BBS^{N-1}$. This implies that $\gw$ depends only on $\gf$.

\medskip

\nind   {\it Sketch of the proof for $s\in(0,1)$}. We perform the same reflection method as in the case $s=\frac12$, defining 
$\varpi$ as in $(\ref{V-5-2})$, it satisfies 
$$
\left.\BA{lll}\displaystyle
\CA_s[\varpi]+\Gl_{s,p,{\scriptscriptstyle N}}\varpi=0\quad&\text{in }\, \Gg_{\phi_0}\\[0mm]
\displaystyle\phantom{\frac{\prt v}{\prt \gn^s}+\ge|v|^{p-1}\, }
\varpi=0\quad&\text{in }\, \BBH_{\gf_0}\cap \BBS^{N}_+\\[0mm]
\displaystyle\phantom{-,,,.-\ \ \, }
\!\varpi(0,\cdot)=\gw(0,\cdot)-\gw(2\gf_0,\cdot)\ \ &\text{in }\, \BBS^{{N-1}}\cap\{(x_1,x_2,...,x_N,0):x_1>0\},
\EA\right.$$
where
$$\CA_s[\varpi]=\Gd_{\BBS^N}\varpi +(1-2s)\cot(\phi)\, \varpi_{\phi}. $$

Since $\frac{\prt\gw}{\prt\gn^s}=\gw^p$ on $\prt \BBS^{N}_+$, the function $\gf\mapsto\gw(\gf,\gs')$ should be decreasing if $(\gf,\gs')\in (0,\gf_0)\ti \BBS^{N-1}$ for $\gf_0$ small enough. Let $\gf^*$ be the maximum of   $\gf_0\in (0,\frac\gp2]$
such that  (\ref{V-5-3-0}) holds. 

If $\gf^*=\frac{\pi}{2}$, we are done. If not, we can assume that  $\gf^*\in(0,\frac{\pi}{2})$ and let $(\gf^*,\gs^*)$ be the point verifying $\varpi(\gf^*,\gs^*)=0$.
Note that $(1-2s)\cot(\cdot)$ is always bounded in a neighborhood of $(\gf^*,\gs^*)$.
Thus, the Hopf's lemma could be applied and the remaining of the  proof  is the same as in the case $s=\frac12$. 
 \qeda\medskip

\nind{\bf Proofs of Theorem A--C. }  Theorem A--C  are exact \rth{Th1}, \rth{Th1*} and \rth{Th2} respectively.  \qeda\medskip

\nind\BBRemark  From the unique positive solution $\gw_2$ obtained in {\it 2} in above theorem, the fractional Lane-Emden equation 
$$
(-\Gd)^sv-v^p=0\qquad\text{in }\BBR^N\setminus\{0\} 
$$
has an explicit self-similar solution $U_p$.  From $(\ref{Li2-u})$ and $(\ref{Li2-c})$ below (also see \cite[Lemma 3.1]{Fa}), we have that
$$
U_p(x)=c_p|x|^{-\frac{2s}{p-1}},\qquad\text{where }\;c_p=\gw_2(1)=\left(2^{2s}
\myfrac{\Gg\big(\frac N2-\frac{s}{p-1}\big)\Gg\big(s+\frac{s}{p-1}\big)}
{\Gg\big(\frac{s}{p-1}\big)\Gg\big(\frac N2-s-\frac{s}{p-1}\big)}\right)^\frac{1}{p-1}.
$$

\nind \BBRemark Another approach could be to transform the problem on the half-sphere $\BBS^N_+$ into a problem on the unit ball $B_1$ of $\BBR^N$ using a stereographic projection from the antipodal point  $\bf S$ of the North pole $\bf N$ on $\BBS^N$. This stereographic projection $\BBP_N$ is an isomorphism from $\BBS^{N}\setminus\{\bf S\}$ onto $\BBR^N$. Then the image of $\BBS^N_+$ is the unit ball $B_1$. It is well known that $\BBP_N$ is a conformal diffeomorphism with conformal factor $D^{\frac{4}{N-2}}$, where
$$D(x)=\Big(\frac{2}{1+|x|^2}\Big)^{\frac{N-2}{2}}.
$$
If $\gw$ is a solution of $(\ref{X-5-1})$, and if 
$$\tilde w(x)=\gw(\gs)D(x)\quad\text{with }\gs=\BBP_N^{-1}(x).
$$
Then 
the above boundary value follows from the fact that
$$\tan\gf=\frac{1-|x|^2}{2|x|}
$$
and if $\gw$ is defined on $\overline {\BBS^N_+}$, we have that
$$\frac{\prt \gw}{\prt\gn}\Big\lfloor_{\gf=0}=\frac{\prt \gw}{\prt\gn}\Big\lfloor_{|x|=1}+\frac{N-2}{2}w\Big\lfloor_{|x|=1}.
$$ 
As a consequence, $\tilde w$ satisfies
$$
\left.\BA{lll}\displaystyle
\Gd\tilde w+\frac{4}{(1+|x|^2)^2}\Big(\Gl_{\frac 12,p,{\scriptscriptstyle N}}+\frac{N(N-2)}{4}\Big)\tilde w=0\ \ &\text{in }B_1\\[3.5mm]
\displaystyle\phantom{--------)-\   }
\frac{\prt \tilde w}{\prt \gn}-\tilde w^{p}+\frac{N-2}{2}\tilde w=0\ \ &\text{on }\prt B_1.
\EA\right.$$
  Strong Maximum Principle and Hopf's Lemma could imply that $\tilde w>\big(\frac{N-2}{2}\big)^{\frac1{p-1}}$ and the standard method of moving planes could be applied to 
obtain the radial symmetry. 



\mysection{A priori estimates}

\subsection{The lifting estimates}

Our aim in this sections is to  obtain the upper bounds for the positive solutions of 
\bel{X1}\left.\BA{lll}\displaystyle
(-\Gd)^sv+\ge v^p=0\quad&\text{in }\, \Omega\setminus\{0\}\\[2mm]
\phantom{(-\Gd)^s-\ge v^p}
v=0\quad&\text{in }\,\Omega^c,
\EA\right.\ee
 under the form 
\bel{Y-1}\BA{lll}\displaystyle
|v(x)|\leq c_3|x|^{-\frac{2s}{p-1}}\qquad\text{for all }x\in \Omega\setminus\{0\},
\EA\ee
where $\Gw\subset\BBR^N$ is a bounded regular domain containing the origin.

In order to study the asymptotics of singular solutions   near the origin, we use the Caffarelli-Silvestre lifting of $v$ and the Poisson kernel associated to the operator $\CD_s$ by formula $(\ref{I-6})$, $(\ref{I-7})$. 
Set $\gm_s(x)=(|x|+1)^{-N-2s}$. If $v\in L_{\mu_s}^1(\BBR^N):=L^1(\BBR^N,\mu_sdx)$, then, with $\xi=(x,z)$, 
\bel{Y-2}\BA{lll}\displaystyle
u(\xi)=u(x,z)=\myint{\BBR^N}{}\CP_s(x-y,z)v(y)dy
\EA\ee
satisfies 
$$
 \CD_s u=0\quad \text{in }\;\, \BBR^{N+1}_+,\qquad
u(\cdot,0)=v \quad \text{on }\;\BBR^{N}.
 $$
\bprop{est} Assume that $p>1+\frac{2s}{N}$ and $u$ is a positive solution of $\CD_s u=0$ in $\BBR^{N+1}_+$ such that $v(x)=u(x,0)$ satisfies estimate $(\ref{Y-1})$. Then there exists $c_4>0$ such that 
\bel{S17}\BA{lll}
u(x,z)\leq c_4\left(|x|^2+z^2\right)^{-\frac{s}{p-1}}\quad\text{for all }\;(x,z)\in\BBR^N\ti\BBR_+.
\EA\ee
\es

\Proof Assume that $u(\cdot,0)$ is nonnegative and satisfies $(\ref{Y-1})$. If $p>1+\frac{2s}{N}$, then $v\in L_{\mu_s}^1(\BBR^N)$ and
  \bel{S12}
\BA{lll}\displaystyle u(x,z)=\myint{\BBR^N}{}\CP_s(x-y,z)v(y)dy\\[4mm]\displaystyle
\phantom{u(x,z)}\displaystyle\leq c_5\myint{\BBR^N}{}\myfrac{z^{2s}}{\left(|x-y|^2+z^2\right)^{\frac{N+2s}{2}}}|y|^{-\frac{2s}{p-1}}dy
\\[5mm]\displaystyle\phantom{u(x,z)}=c_5z^{-\frac{2s}{p-1}} \myint{\BBR^N}{}\myfrac{|y|^{-\frac{2s}{p-1}}}{\left(|x_z-y|^2+1\right)^{\frac{N+2s}{2}}} dy,
\EA \ee
where we take $x_z=\frac{x}z$.

When $0\leq |x_z|\leq 8$, we obtain that 
$$\myint{\BBR^N}{}\myfrac{|y|^{-\frac{2s}{p-1}}}{\left(|x_z-y|^2+1\right)^{\frac{N+2s}{2}}} dy\leq \myint{\BBR^N}{} \left(|y|^2+1\right)^{-\frac{N+2s}{2}}|y|^{-\frac{2s}{p-1}}dy:=c_6 <+\infty.$$

When $|x_z|>8$, letting $r=\frac12|x_z|$, we have that 
$$\BA{lll}\myint{B_r(0)}{}\myfrac{|y|^{-\frac{2s}{p-1}}}{\left(|x_z-y|^2+1\right)^{\frac{N+2s}{2}}}dy\leq (1+|x_z|^2)^{-\frac{N+2s}{2}} \myint{B_r(0)}{} |y|^{-\frac{2s}{p-1}}dy\\[5mm]\displaystyle\phantom{------------\ }  
=c_7 |x_z|^{-N-2s}r^{N-\frac{2s}{p-1}}, 
\EA$$
$$\BA{lll}\myint{B_r(x_z)}{}\myfrac{|y|^{-\frac{2s}{p-1}}}{\left(|x_z-y|^2+1\right)^{\frac{N+2s}{2}}} dy\leq \big(\frac r2\big)^{-\frac{2s}{p-1}} \myint{B_r(x_z)}{} \left(|x_z-y|^2+1\right)^{-\frac{N+2s}{2}} dy
\\[5mm]\displaystyle\phantom{ ------------\ \ \, } 
 \leq c_8 r^{-\frac{2s}{p-1}}\big(1+ r^{-2s} \big),
\EA$$
and 
$$\BA{lll}\myint{\BBR^N\setminus \big(B_r(x_z)\cup B_r(0)\big)}{}\myfrac{|y|^{-\frac{2s}{p-1}}}{\left(|x_z-y|^2+1\right)^{\frac{N+2s}{2}}} dy\leq c_9 \myint{\BBR^N\setminus B_r(0)}{} \left(|y|^2+1\right)^{-\frac{N+2s}{2}}|y|^{-\frac{2s}{p-1}} dy\\[5mm]\displaystyle\phantom{------------------ }  \leq 
c_{10} r^{-\frac{2s}{p-1}-2s}. 
\EA$$
Thus, we conclude that for $|x_z|\geq 8$,
$$\myint{\BBR^N}{}\myfrac{|y|^{-\frac{2s}{p-1}}}{\left(|x_z-y|^2+1\right)^{\frac{N+2s}{2}}}dy\leq c_{11}  |x_z|^{-\frac{2s}{p-1}}.$$
As a consequence, we obtain that for $|x_z|>0$,
$$
\BA{lll}\displaystyle u(x,z)\leq c_{12} z^{2s}   (1+|x_z|^2)^{-\frac{s}{p-1}}\leq c_{12}  (z^2+|x|^2)^{-\frac{s}{p-1}}.
\EA $$
Which ends the proof.  \hfill \qeda\medskip

\nind{\bf Remark. } Since $(\ref{S17})$ can imply $(\ref{Y-1})$ immediately by taking $z=0$, the previous result shows that these two estimates  are equivalent when $p>1+\frac {2s}{N}$.

\subsection{Integral estimates}
\blemma{lm 4.1-s}
Let $p>0$ and $v\in C^{1,s}({\overline \Gw}\setminus\{0\})$ be a nonnegative solution of $(\ref{X1})$ in $\Gw\setminus\{0\}$, which vanishes in $\Gw^c$.  
Then  for any $\xi\in C^{1,s}({\overline \Gw } \setminus\{0\})$ with
compact support in ${\overline \Gw}\setminus\{0\}$, there holds
\bel{2.x4}\int_{\Gw}\big(v(-\Gd)^s\xi+\ge v^p\xi\big)dx=0.
\ee
Furthermore, if $v^p\in L^1(\Gw)$,which is always satisfied for $\ge=-1$, there exists $k\geq 0$ such that 
\bel{2.y4}
(-\Gd)^sv+\ge v^p=k\gd_0\quad\mbox{in }\, \CD'(\Gw).
\ee
\es
\Proof Without loss of generality, we can assume that $B_2\subset\Gw.$ Let $\eta_0:\BBR^N\to[0,1] $ be a nonnegative cut-off function in $C^\infty_c(\BBR^N)$,  equal to $1$ in $B_1$, vanishing in $B_2^c$. For $0<\vge<1$ we set
$\tilde\eta_\vge(x)=\eta_0(\vge x)(1-\eta_0(x/\vge)$ and denote $v_\vge=\tilde\eta_\vge v$. Then
$$\BA{lll}
(-\Gd)^sv_\vge(x)=\tilde\eta_\vge(x)(-\Gd)^sv(x)+v(x)(-\Gd)^s\tilde\eta_\vge(x)\\[2mm]
\phantom{--------}\displaystyle 
-c_{N,s}\int_{\BBR^N} \frac{(v(x)-v(y))(\tilde\eta_\vge(x)- \tilde\eta_\vge (y))}{|x-y|^{N+2s}}dy,\quad \forall\, x\in B_1\setminus\{0\}.
\EA$$
Furthermore,  
\begin{eqnarray*}
 (-\Delta)^s v_\vge (0)=\lim_{x\to0} (-\Delta)^s v_\vge (x)=-c_{N,s}\int_{ B_{\frac2\vge}  \setminus B_\vge }  \frac{v(y)\eta_\vge(y)}{|y|^{N+2s}}  dy.
\end{eqnarray*}
We can assume that $\xi$ vanishes in $\overline B_r$ and we choose $0<\vge<r/8$. By integration, there holds  
 $$
  \int_{B_1} \xi (-\Delta)^s v_\vge \,dx=  \int_{B_1}  v_\vge  (-\Delta)^s\xi\,dx.
 $$
Since $v=v_\vge+(1-\tilde\eta_\vge)v$ and $\xi(1-\tilde\eta_\vge)=0$, we have that
$$\BA{lll}\displaystyle
\int_{B_1}\left[ v (-\Delta)^s \xi  +\ge v^p\xi\right]\,dx=\int_{B_1}\left[\xi v(-\Gd)^s\tilde\eta_\vge+(1-\tilde\eta_\vge)v(-\Gd)^s\xi\right]dx\\[4mm]
\phantom{------------}
\displaystyle-{c_{N,s}}\int_{B_1\setminus B_r}\left(\int_{\BBR^N} \frac{(v(x)-v(y))(\tilde\eta_\vge(x)-\tilde\eta_\vge (y))}{|x-y|^{N+2\alpha}}dy\right)\xi(x)\,dx
\\[4mm]
\phantom{\displaystyle\int_{B_1}\left[ v (-\Delta)^s \xi  +\ge v^p\xi\right]\,dx}
=A_{1,\vge}+A_{2,\vge}+A_{3,\vge}.
\EA
$$
Since $v\in L^1(B_1)$ and $|(-\Delta)^s \xi(x)|\leq c_{13}(1+|x|)^{-N-2s}$ for any $ x\in\BBR^N$ and some $c_{13}>0$, $A_{2,\vge}\to 0$ when $\vge\to 0$. We have also for 
$r<|x|<1$ (where $\xi(x)$ may not vanish),
$$\BA{lll}\displaystyle
|(-\Gd)^s\tilde\eta_\vge(x)|\leq c_{N,s}\int_{B_{2\vge}}\frac{1}{|x-y|^{N+2s}}dy+c_{N,s}\int_{B_{2/\vge}}\frac{1}{|x-y|^{N+2s}}dy\\[3.5mm]
\phantom{|(-\Gd)^s\tilde\eta_\vge(x)|}
\leq c_{14}\left(\vge^N(r-2\vge)^{-N-2s}+\vge^{2s}\right),
\EA
$$
for some $c_{14}>0$ independent of $\vge$. This implies that $A_{1,\vge}\to 0$ when $\vge\to 0$. Finally, for $x\in B_1\setminus B_r$ and $\vge <\frac r4$, there holds
$$\BA{lll}\displaystyle
 |A_{3,\vge}|\leq \norm\xi_{L^\infty}\Big(\int_{B_1\setminus B_r}\int_{B_{2\vge} } \frac{|v(x)-v(y)|}{|x-y|^{N+2s}}dydx+\int_{B_1\setminus B_r}\int_{  B_{\frac1\epsilon}^c} \frac{|v(x)-v(y)|}{|x-y|^{N+2s}}dydx\Big) 
 \\[4mm]\phantom{ \displaystyle
 |A_{3,\vge}|}\displaystyle
\leq c_{15} (r-2\varepsilon)^{-N-2s} \norm{\xi}_{L^\infty } \Big(  (\varepsilon^N +\epsilon^{2s})\int_{B_1} v(x)\,dx  +\int_{B_{2\varepsilon}} v(y) dy  
 + \int_{B_{\frac1\varepsilon}} v(y)  d\mu_s(y)  \Big),
\EA
$$
where $c_{15}>0$ is independent of $\vge$. 
Hence $A_{3,\vge}\to 0$ as $\vge\to 0$.  This implies $(\ref{2.x4})$.\smallskip

If $\ge=-1$, it follows from \cite[Theorem 3]{LWX} that $v^p\in L^1(B_1)$ and there exists $k\geq 0$ such that $(\ref{2.y4})$ holds.

If $\ge=1$ and we assume now that $F:=v^p\in L^1(B_1)$, the function $w=v+\BBG_{s,1}[F]$ is $s$-harmonic and positive in $B_1\setminus\{0\}$. By the extension to 
$s$-harmonic functions of B\^ocher's theorem \cite[Theorem 4]{LWX}, there exists $k\geq 0$ such that 
$$(-\Gd)^sw=k\gd_0\quad\mbox{in }\CD'(B_1(0)),
$$
and this is $(\ref{2.y4})$.
\hfill$\Box$ \medskip

\subsection{A priori estimate for   Lane-Emden equation}

The following result is proved in \cite{YZ}.
\bprop{estLE} Assume that  $0<s<1$, $1<p<\frac{N+2s}{N-2s}$ and $u$ is a nonnegative solution of 
$$
\left.\BA{lll}
\phantom{\displaystyle \prt_{\gn^s}u(\cdot,0)-u^p\,  }\CD_s u=0\quad&\text{in }\,{\bf B}_1\cap\BBR^{N+1}_+\\[2mm]
\displaystyle \prt_{\gn^s}u(\cdot,0)-u(\cdot,0)^p=0&\text{in }\,  B_1\cap\BBR^{N}\setminus\{0\}.
\EA\right.$$
Then there exists a constant $c_{16}=c_{16}(N,p,s)>0$ such that 
\bel{WW-2}\BA{lll}
u(x,0):=v(x)\leq c_{16} |x|^{-\frac{2s}{p-1}}\quad\text{for all }\;x\in B_\frac{1}{2}\cap\BBR^{N}\setminus\{0\}.
\EA\ee
\es
The upper bound of $v$ comes from the technique used in the proof which is based upon a blow-up argument and a Liouville type theorem only valid in this range of exponents. When $
v$ is a radial function, we have a much shorter proof.  

\bprop{estLE-2} Let $0<s<1$, $p>1$, $\Omega=B_1$ and $v$ be a positive  and radially decreasing solution of $(\ref{X1})$ with $\ge=-1$.
Then there exists a constant $c_{17}=c_{17}(N,p,s)>0$ such that 
$$
v(x)\leq c_{17} |x|^{-\frac{2s}{p-1}}\quad\text{for all }\;x\in B_\frac{1}{2}\setminus\{0\}.
$$
\es
\Proof Let $\gz\in C^\infty_c(\overline B_1)$ be a radially symmetric function with value in $[0,1]$ and such that 
$$\gz(x)=\left\{\BA{lll}1\quad&\mbox{for }\, \, x\in \overline{ B_\frac12\setminus B_\frac14}\\[2mm]
0& \mbox{for }\, \, x\in B_\frac18 \cup  B_1^c  
\EA
\right.
$$
and for $n>1$, $\gz_n(x)=\gz(nx)$. By \rlemma{lm 4.1-s}, there holds
$$
 \int_{B_1}v^p\gz_n^{p'}(x)dx=\int_{B_1}v(-\Gd)^s\gz_n^{p'}(x)dx.
$$
From \cite[Lemma 2.3]{CV},  there holds
$$(-\Gd)^s\gz_n^{p'}=p'\gz_n^{p'-1}(-\Gd)^s\gz_n-K^2,
$$
where 
$$K^2:=K^2(x)=\frac{p'(p'-1)}{2}\gz_n^{p'-2}(z_x)\int_{B_1}\frac{(\gz_n(x)-\gz_n(y))^2}{|x-y|^{N+2s}}dy
$$
for some $z_x\in \overline B_1$. Hence
$$\BA{lll}\displaystyle
\int_{B_1}v^p\gz^{p'}_n(x)dx=\int_{B_1}v(-\Gd)^s\gz^{p'}_n(x)dx
\leq p'\int_{B_1}v\gz^{p'-1}_n(-\Gd)^s\gz_n(x)dx.
\EA$$
Now 
$$
\BA{lll}\displaystyle\int_{B_1}v\gz^{p'-1}_n(-\Gd)^s\gz_n(x)dx=\int_{\frac{1}{8n}<|x|<\frac{1}{n}}v\gz^{p'-1}_n(-\Gd)^s\gz_n(x)dx\\[4mm]
\phantom{\displaystyle\int_{B_1}v\gz^{p'-1}_n(-\Gd)^s\gz_n(x)dx}
\displaystyle\leq \Big(\int_{\frac{1}{8n}<|x|<\frac{1}{n}}v^p\gz^{p'}_n(x)dx\Big)^\frac{1}{p}\Big(\int_{\frac{1}{8n}<|x|<\frac{1}{n}}|(-\Gd)^s\gz_n(x)|^{p'}dx\Big)^\frac{1}{p'}.
\EA$$
Therefore 
$$
\int_{B_1}v^p\gz^{p'}_n(x)dx=\int_{\frac{1}{8n}<|x|<\frac{1}{n}}v^p\gz^{p'}_n(x)dx\leq p'^{p'}
\int_{\frac{1}{8n}<|x|<\frac{1}{n}}|(-\Gd)^s\gz_n(x)|^{p'}dx.
$$
Now 
$$\int_{\frac{1}{8n}<|x|<\frac{1}{n}}v^p\gz^{p'}_n(x)dx\geq \int_{\frac{1}{4n}<|x|<\frac{1}{2n}}v^p\gz^{p'}_n(x)dx\geq  \gw_{\scriptscriptstyle N} \frac{4^N-1}{4^N}n^{-N}v^p\Big(\frac{1}{2n}\Big),
$$
and 
$$\int_{\frac{1}{8n}<|x|<\frac{1}{n}}|(-\Gd)^s\gz_n(x)|^{p'}dx=n^{2sp'-N}\int_{\frac{1}{8}<|y|<1}|(-\Gd)^s\gz(y)|^{p'}dy.
$$
Since $\gz$ vanishes in $B_1^c$, 
$$
(-\Gd)^s\gz(x)=\left\{\BA{lll}\displaystyle
c_{{\scriptscriptstyle N},s}{\rm p.v.}\int_{B_1}\frac{\gz(x)-\gz(y)}{|x-y|^{N+2s}}dy+c_{N,s}\gz(x)\int_{B_1^c}\frac{dy}{|x-y|^{N+2s}}\ \, &\mbox{for }\, x\in B_1\\[4mm]\displaystyle
-c_{{\scriptscriptstyle N},s}\int_{B_1}\frac{\gz(y)}{|x-y|^{N+2s}}dy\ \, &\mbox{for }\, x\in B^c_1.
\EA\right.$$
Therefore, $\norm{(-\Gd)^s\gz}_{L^k(B_n)}$ is bounded independently of $k$ and $n$. 
Combining the previous inequalities, we obtain
$$
\gw_{\scriptscriptstyle N} \frac{4^N-1}{4^N}n^{-N}v^p\Big(\frac{1}{2n}\Big)\leq n^{2sp'-N}\int_{\frac{1}{8}<|y|<1}|(-\Gd)^s\gz(y)|^{p'}dy,
$$
which implies that there holds, for some $c_{18}>0$,
$$v\Big(\frac{1}{2n}\Big)\leq c_{18}n^{\frac{2s}{p-1}}.
$$
Taking $r=\frac {1}{2n}$, we  infer $(\ref{I-2})$. 
 \hfill$\Box$\medskip
 
 \nind\BBRemark When $\Omega=B_1$, it is easy to prove that any positive solution $v$ of $(\ref{X1})$ satisfying
 $$
 \lim_{x\to 0}v(x)=+\infty
$$
is radially symmetric and decreasing.
 \subsection{Some maximum principles}
 In this section we recall some fundamental maximum principles which will be used in the sequel. Their proof can be found in \cite{CFQ}.
 \blemma{cr 2.0} Let  $\CO\subset\BBR^N$ be any domain and $v\in C(\overline \CO)\cap L_{\mu_s}^1(\CO)$. If v  achieves the maximum at some point $x_0\in\CO$, then
$$
 (-\Delta)^s v(x_0)\geq  0.
 $$
Furthermore equality holds if and only if 
$$v(x)=v(x_0)\ \  \text{ for all}\  x\in\CO.$$
\es
The following corollaries follow from the strong maximum principle.
\bcor{teo 2.1-cr}
Let $\CO\subset\BBR^N$ be a bounded regular domain and $v\in C(\overline\CO)$ satisfy 
$$  (-\Delta)^s v \leq 0\quad \text{in }\,  \CO,\qquad
v=h\quad \text{in }\,\CO^c, $$
where $h$ is continuous and bounded in $\CO^c$. Then

$$\displaystyle\sup_{x\in\CO} v(x)\leq \sup_{x\in\CO^c}h(x).$$
\es

\bcor{teo 2.1-com}
Let $\CO\subset\BBR^N$ be a bounded regular domain, $g:\BBR\to\BBR$ be a continuous nondecreasing function  and $v_j\in C(\overline\CO)$ (j=1,2)   satisfy the two inequalities
$$\BA{lll}(-\Delta)^s  v_1+g(v_1) \geq  0\quad\text{and }\;
(-\Delta)^s  v_2+g(v_2)\leq 0 \quad \text{in }\,  \CO.
\EA
$$
Assume furthermore that $  v_1\geq  v_2$ a.e. in $\CO^c$.
Then  either $ v_1> v_2$ in $\CO$ or $v_1\equiv v_2$   in $\BBR^N.$
\es
The next statement proved in \cite[Lemma 3.1]{CV} is the analogue of the classical local regularity result for the Laplace operator, see also \cite[Theorem 2.2]{JLX} for local Schauder estimates.

\blemma{lm 2.1}
Assume that $w\in C^{2s+\varepsilon}(\overline B_1)\cap L_{\gm_s}^1(\BBR^N)$ for some $\varepsilon\in (0,1)$
 satisfies
 $$(-\Delta)^s w=f\quad\ \text{in }\;  B_1,$$
 where $f\in C^1(\overline B_1)$. Then for $\beta\in (0,2s)$, there exist  $ c_{19}, c_{20}>0$ such that
$$
\norm w_{C^\beta(\bar B_{1/2})}\leq c_{19}\Big(\norm w_{L^\infty(B_1)}+\norm f_{L^\infty(B_1)}+\norm{w}_{L_{\gm_s}^1(\BBR^N)}\Big) 
$$
and
$$
\norm w_{C^{2s+\beta}(\bar B_{1/4})}\leq  c_{20}\Big(\norm w _{C^\beta(B_\frac12)}+\norm f_{C^\beta(B_\frac12)}+\norm{w}_{L_{\gm_s}^1(\BBR^N)}\Big).
$$
\es

The next general result proved in \cite {FT} is an important tool 
for obtaining the existence of solutions of semilinear equations. 

\bprop{th 3.2.4}
Assume that $\Omega\subset\BBR^N$ is a  bounded regular domain, $g:\BBR\to\BBR$  is a continuous nondeacreasing function,  $f\in C^\beta(\Omega)$  with $\beta\in(0,1)$.
If there exist
a super-solution $\bar v$ and a sub-solution $\underline v$ of
\bel{eq 1.1-ss}\left.\BA{lll}
 (-\Delta)^s  v+ g(v)=f\ \ \ & \rm{in}\ \, \Omega\setminus\{0\}\\[2mm]
 \phantom{   (-\Delta)^s  + g(v)  }
v\geq  0\ \ & \rm{in}\ \, \Omega\setminus\{0\}\\[2mm]
 \phantom{   (-\Delta)^s  +  g(v)  }
\displaystyle v=0 &\rm{in}\ \, \Omega^c 
\EA\right.
\ee
such that  $\bar v,\underline v\in C^2(\Gw\setminus\{0\})\cap L^1(\Gw)$,
$$\bar v\geq  \underline v\geq  0\ \ \mbox{in}\ \Omega\setminus\{0\}\quad 
\mbox{and}\quad \bar v\geq  0\geq   \underline v \ \ \mbox{in}\ \Omega^c.
$$
Then
there exists at least one solution $v$ of $(\ref{eq 1.1-ss})$ which satisfies
$\underline v\leq v\leq\bar v \ \ \mbox{in}\ \ \Omega\setminus\{0\}.
$
\es  

We end this subsection with some comparison statements.

\bprop{lm 2.1-non}
 Let  $v $
be a nonnegative s-subharmonic function in $\BBR^N\setminus\{0\} $ satisfying
$$
 \lim_{x\to 0}|x|^{N-2s}v(x)=0\quad {\rm and}\quad \lim_{|x|\to+\infty}v(x)=0,
$$
then 
$$v\equiv 0\quad {\rm in}\ \ \BBR^N. $$
\es
\Proof The function $x\mapsto w(x):= |x|^{2s-N}+1$ is $s$-harmonic in $\BBR^N\setminus\{0\}$ and for any $\vge>0$ the function $v-\vge w$   is $s$-subharmonic in $\BBR^N\setminus\{0\}$  and  negative in $B_R^c$ for $R>0$ large enough and in $B_\gd$ for $\gd>0$ small enough. Taking $\CO=B_{R}\setminus {\overline B_\delta}$ and $h=0$ in   \rcor{teo 2.1-cr}, it implies  that for any $\vge>0$,
$$0\leq v\leq \vge w\quad{\rm  in}\ \,  \BBR^N\setminus\{0\}.$$  
Thus passing to the limit as $\vge\to 0$ implies the claim.  
\qeda\medskip

As a variant of this result we have

\bcor{cr 2.1-non}
 Let $\Omega\subset\BBR^N$ be a bounded regular domain containing the origin, $g:\BBR_+\to\BBR_+$ be a continuous function such that $g(0)=0$ and  $h$ be a nonnegative function bounded in $\BBR^N\setminus\Omega$. If $v$ is a  nonnegative solution of 
$$
 (-\Delta)^s  v+g(v)=0\ \   {\rm in}\ \ \Omega\setminus\{0\},\qquad
   v=h \ \    {\rm in}\ \ \Omega^c, 
$$
 satisfying 
$
\displaystyle \lim_{x\to 0}v(x)|x|^{N-2s}=0.
$
Then $v\leq \norm h_{L^\infty}$ in $\BBR^N.$
\es



 \subsection{Fractional Hardy operator}

For $\tau\in(-N,2s)$, we define  $\CC_s(\tau)$ by the expression
$$
 \CC_s(\tau) = 2^{2s} \frac{\Gamma(\frac{N+\tau}{2})\Gamma(\frac{2s-\tau}{2})}{\Gamma(-\frac{\tau}{2})\Gamma(\frac{N-2s+\tau}{2})}= -\frac{c_{N,s}}2 \int_{\BBR^N}\frac{|{\bf e}_1+z|^{\tau}+|{\bf e}_1-z|^\tau-2}{|z|^{N+2s}}\,dz,
$$
where ${\bf e}_1=(1,0,\cdots,0)\in\BBR^N$. 
Since the function $\Gg$ is infinite at $0$,  $\CC_s$ vanishes at $0$ and $2s-N$, i.e.  
$$\CC_s(0)=\CC_s(2s-N)=0.$$
The role of $C_s(\tau)$ is enlighted by the following identity \cite{CW}
$$
 (-\Delta)^s |x|^\tau = \CC_s(\tau)|x|^{\tau-2s} \quad\  \mbox{ for all }\, x\in   \BBR^N\setminus\{0\}. 
$$
The function $\CC_s$ 
is concave and achieves its maximum in $(-N,2s)$  for $\tau=\frac{2s-N}{2}$ with corresponding maximal value $2^{2s}  \frac{\Gamma^2(\frac{N+2s}4)}{\Gamma^2(\frac{N-2s}{4})}$ (see \cite[Lemma 2.3]{CW}).
Furthermore,
$$
\CC_s(\tau)= \CC_s(2s-N-\tau) \quad \text{for } \, \tau \in (-N,2s),
$$
and
$$
\lim_{\tau \to -N}\CC_s(\tau) = \lim_{\tau \to 2 s}\CC_s(\tau)=-\infty.
$$
Let $\tau_p=-\frac{2s}{p-1}$,  then  we
have that  for $1+\frac{2s}{N}<p<\frac{N}{N-2s}$,
$$\tau_p\in (-N,2s-N)\quad{\rm and}\quad \CC_s(\tau_p)>0$$ and for $p>\frac{N}{N-2s}$ 
$$\tau_p\in (2s-N, 0)\quad{\rm and}\quad\CC_s(\tau_p)<0.$$ 
 Let 
  \bel{Li2-c}
  c_p=\left\{
\BA{lll}
\CC_s(\tau_p)^{\frac{1}{p-1}}\ \ & \text{if}\ \ 1+\frac{2s}{N}<p<\frac{N}{N-2s}\\[2mm]
 \phantom{      }
\big(-\CC_s(\tau_p)\big)^{\frac{1}{p-1}}\ \ & \text{if}\ \ p>\frac{N}{N-2s},
\EA \right. \ee
then the function $U_p$ expressed by 
 \bel{Li2-u}
 U_p(x):=c_p|x|^{-\frac{2s}{p-1}}\quad{\rm in}\ \,  \BBR^N\setminus\{0\}, 
 \ee
 satisfies,  for $1+\frac{2s}{N}<p<\frac{N}{N-2s}$,
$$(-\Delta)^sU_p+U_p^p=0\quad {\rm in}\ \,  \BBR^N\setminus\{0\},$$ 
 and for $p>\frac{N}{N-2s}$
 $$(-\Delta)^sU_p=U_p^p\quad {\rm in}\ \,  \BBR^N\setminus\{0\}. $$ 
The fractional Hardy operator is defined by 
$$\CL^s_\mu v(x)=(-\Delta)^s v(x)+\frac{\mu}{|x|^{2s}}v(x), $$
under the condition 
$$
\mu\geq  \mu_0:=-2^{2s}\frac{\Gamma^2(\frac{N+2s}4)}{\Gamma^2(\frac{N-2s}{4})}.$$
It is shown in \cite{CW} that   the linear equation 
$$\mathcal{ L}_\mu^s v=0\qquad\text{in }\, \ \BBR^N\setminus \{0\}$$
 has two distinct radial solutions $\Phi_{\mu}$ and $\Gamma_{ \mu}$ defined by
$$
\Phi_{\mu}(x)=\left\{\arraycolsep=1pt
\begin{array}{lll}
 |x|^{\tau_-( \mu)}\quad
   &{\rm if}\quad \mu>\mu_0\\[1.5mm]
 \phantom{   }
|x|^{-\frac{N-2s}{2}}\ln\left(\frac{1}{|x|}\right) \quad  &{\rm   if}\quad \mu=\mu_0
 \end{array}
 \right.\qquad  {\rm and}\quad\  \Gamma_{ \mu}(x)=|x|^{\tau_+(\mu)},$$
where the exponents 
$\tau_-( \mu)  \leq  \tau_+( \mu)$ verify the relations 
 \begin{eqnarray*}
 &\tau_-(\mu)+\tau_+(\mu) =2s-N \quad\  {\rm for\ all}\ \   \mu \geq  \mu_0,\nonumber\\[2mm]
&\tau_-( \mu_0)=\tau_+( \mu_0)=\frac{2s-N}2,\quad\ 
\tau_-( 0)=2s-N, \quad\ \tau_+(0)=0,  \\[2mm]
&\displaystyle\lim_{\mu\to+\infty} \tau_-( \mu)=-N\quad {\rm and}\quad \lim_{\mu\to+\infty} \tau_+(\mu)=2s\nonumber.
  \end{eqnarray*}
In the remaining of the paper and when there is no ambiguity,  we use the notations $\tau_+=\tau_+(\mu)$, $\tau_-=\tau_-( \mu)$.  
The following inequalities are proved in  \cite{CW}, 
$$\CC_s(\tau)+\mu>0\quad {\rm for}\ \ \tau\in(\tau_-,\tau_+)$$
and
$$\CC_s(\tau)+\mu<0\quad {\rm for}\ \ \tau\in(-N,\tau_-) \cup (\tau_+,2s).$$

\nind {\bf Definition 3.1} 
We denote by
$W^s(\Omega)$ the space of functions
$v \in L^1(\BBR^N, d\mu_s)$ such that $v\big|_{\Omega} \in L^2(\Omega)$ which satisfy
$$
\CE^s_{\Omega,\Omega'}(v,v)<+\infty\quad  {\rm for\ some\ domain\ }\Omega' \subset \BBR^N\ {\rm with}\  \Omega \subset \subset \Omega',
$$
where we have set 
$$
\CE^{s}_{A,B}(v,w) = \frac{c_{N,s}}{2} \int_{A \times B} \frac{(v(x)-v(y))(w(x)-w(y))}{|x-y|^{N+2s}}dx dy.
$$
We define $W^s_*(\Omega)$ as the space of functions
$v \in L^1(\BBR^N, d\mu_s)$ such that $v|_{\Omega\setminus B_\epsilon} \in W^s(\Omega\setminus B_\epsilon)$ for every $\epsilon>0$.
\medskip

 The comparison principle for the fractional Hardy operator is the following, see \cite[Lemma 4.11]{CW},
  
\blemma {Comparison principle, Version 2}
\label{L-comp}
Let $v \in W^s_*(\Omega)$ satisfy $\CL^s_\mu v \geq  0$ in $\Omega \setminus \{0\}$ and
$$
v \geq  0 \quad  {\rm in}\  \BBR^N\setminus \Omega , \qquad \liminf_{x \to 0} \frac{v(x)}{\Phi_\mu(x)} \geq  0.
$$
Then  $v \geq  0$ in $\BBR^N$.
\es

As an application we have the following Liouville type result.
\bcor {Liouv} Let $0<s<1$ and $p>0$. If $v$ is a solution of 
  \bel{Li1}
(-\Gd)^sv+|v|^{p-1}v=0\quad\text{in }\,\BBR^N\setminus\{0\},
\ee
satisfying 
  \bel{Li2}
\displaystyle \lim_{x\to 0}|x|^{N-2s}v_+(x)=0\quad\text{and }\,\lim_{|x|\to +\infty}|x|^{\gt_0}v_+(x)=0,
\ee
for some $\gt_0\in (0,2s)\cap[2s(p+1)-Np,+\infty)$, then $v_+\equiv 0$.
\es
\Proof Let $\vge>0$ and 
$$w_\vge(x)=\vge(|x|^{2s-N}+1)+\vge^{p+1}|x|^{\gt_0}\quad\text{for }x\in\BBR^N\setminus\{0\}.
$$
Then
$$w^p_\vge(x)\geq \vge^p(|x|^{(2s-N)p}+1), 
$$
and 
  \bel{Li3}\BA{lll}(-\Gd)^sw_\vge(x)+w^p_\vge(x)= \vge^{p+1}\CC_{s}(\gt_0)|x|^{\gt_0-2s}+w^p_\vge(x)\\[2mm]
  \phantom{(-\Gd)^sw_\vge(x)+w^p_\ge(x)}\geq \vge^{p}\left(\vge\CC_{s}(\gt_0)|x|^{\gt_0-2s}+|x|^{(2s-N)p}+1\right),
\EA\ee
where $\CC_{s}(\gt_0)<0$. Since $\tau_0\in (0,2s)\cap[2s(p+1)-Np,+\infty)$, we have that $\gt_0-2s\geq (2s-N)p$. Thus, taking $\vge>0$  small enough we infer that the right-hand side of $(\ref{Li3})$ is nonnegative. Hence $w_\vge$ is a supersolution of $(\ref{Li1})$, larger than $v$ at $0$ and at 
$\infty$. Thus $v\leq  w_\vge$. Letting $\vge\to 0$ yields the claim.\qeda.\medskip


\subsection{Emden-Fowler equation with Dirac mass}
The Emden-Fowler equation with right-hand side measure is studied in \cite{CV} and \cite{CV1}. A particular case which is important in our present study is the following problem
\bel{DD1} \BA{lll}
(-\Gd)^sv+v^p=k\gd_0\ \ &\text{in}\;\CD'(\Gw)\\[2mm]
\phantom{(-\Gd)^s+   v^p}
v=0\ \ &\text{in }\;\Gw^c,
\EA \ee
where $\Gw\subset\BBR^N$ is a bounded domain containing the origin. We set 
$$\BBX_{s,\Gw}=\left\{\gz\in C(\BBR^N):\gz=0\text{ in }\Gw^c,|(-\Gd)^s\gz|\in L^\infty(\BBR^N),\,|(-\Gd)_\vge^s\gz|\leq\gf \, \text{ if }0<\vge\leq\vge_0\right\},
$$
where 
$$(-\Gd)_\vge^s\gz(x)=\int_{|x-y|>\vge}\myfrac{\gz(x)-\gz(y)}{|x-y|^{N+2s}}dy
$$
and $\gf \in L^1(\Gw, d_{\prt\Gw}dx)$ with $d_{\prt\Gw}(x)=\dist (x,\Gw^c)$.  The following result is proved in \cite{CV}.
\bprop{Dirac-1} Let  $0<p<\frac{N}{N-2s}$, then for any $k\geq 0$ there exists a unique function $v=v_{k}$ belonging to $L^1(\Gw)$ such that $v_k^p\in L^1(\Gw, d_{\prt\Gw}dx)$,  satisfying 
$$\BA{lll}\displaystyle
\int_{\Gw}\left(v_k(-\Gd)^s\gz+v_k^p\gz\right) dx=k\gz(0)\qquad&\text{for all }\gz\in\BBX_{s,\Gw}\\
\phantom{\displaystyle
\int_{\Gw}\left((-\Gd)^s\gz+v_k^p\gz\right) dx}v_k=0\qquad&\text{in }\Gw^c
\EA$$
and the following estimate holds
\bel{DD3}\BA{lll}
kG_{s,\Gw}(x,0)-c_{21}k^p\leq v_{k}(x)\leq kG_{s,\Gw}(x,0)\leq kG_{s}(x,0)\quad\text{for all }x\in\Gw\setminus\{0\},
\EA\ee
where $c_{21}>0$, $G_{s,\Gw}$ and $G_s$ denote the Green functions of $(-\Gd)^s$ in $\Gw$ and $\BBR^N$ respectively.
Furthermore, the mapping $k\mapsto v_k$ is increasing and it is stable in the sense that if $\{\gr_\ge\}$ is a sequence of nonnegative   functions with compact support in $\Gw$ converging to $k\gd_0$ in the sense of distributions in $\Gw$, then the sequence of solutions $v:=v_{\gr_\ge}$ of 
$$
(-\Gd)^sv+v^p=\gr_\ge\ \ \text{in}\;\CD'(\Gw),\qquad
v=0 \ \ \text{in }\;\Gw^c, $$
converges to $v_{\Gw,k}$ uniformly on any compact subset of $\overline\Gw\setminus\{0\}$, and $v_{\gr_\ge}$ converges to $v_{\Gw,k}$ in $L^p(\Gw,d_{\prt\Gw}dx)$. \es

Since  $k\mapsto v_k$ is increasing, we set 
$$
v_{ \infty}(x)=\lim_{k\to+\infty}v_{ k}(x)\quad\text{for all }x\in\Gw\setminus\{0\}.
$$


\bprop{Dirac-2} Let   $0<p<\frac{N}{N-2s}$, then \smallskip

\nind 1- If $0<p\leq 1+\frac{2s}{N}$, then $v_{ \infty}(x)=\infty$ for all $x\in\Gw$.\smallskip

\nind 2- If $1+\frac{2s}{N}< p<\frac{N}{N-2s}$, then $v_{ \infty}$ is a positive solution of 
\bel{DD6} 
(-\Gd)^sv+v^p=0\ \ \text{in }\;\Gw\setminus\{0\},\qquad
 v=0 \ \ \text{in }\;\Gw^c, 
 \ee
satisfying 
\bel{DD7}\BA{lll}
\displaystyle
\lim_{x\to 0}|x|^{\frac{2s}{p-1}}v_{ \infty}(x)=c_p,
\EA\ee
and
\bel{DD8}\BA{lll}
\displaystyle
U_p(x)-c^*\leq v_{ \infty}(x) \leq U_p(x)\quad\text{for all}\;x\in\Gw\setminus\{0\},
\EA\ee
where $c^*>0$,
$c_p$ and $U_p$ are defined in $(\ref{Li2-c})$ and $(\ref{Li2-u})$ respectively.
\es
\Proof When $1+\frac{2s}{N}\geq \frac{2s}{N-2s}$, the above result is proved in \cite[Theorem 1.1]{CV1}, but when $1+\frac{2s}{N}< \frac{2s}{N-2s}$ the situation is not yet completely clarified \cite[Theorem 1.2]{CV1}. In this proof, we give a proof to overcome this gap. 

To distinguish the effect of domain, we replace $v_k$ by $v_{\Gw,k}$ and for simplicity in the proof of this proposition, we set
$$v_{0,k}=v_{\BBR^N,k},\quad \quad v_{0,\infty}=v_{\BBR^N,\infty}$$
and
$$v_{1,k}=v_{B_1,k},\quad  \quad v_{1,\infty}=v_{B_1,\infty}.$$
From \cite[Theorem 1.4]{CY} for $\Omega=\BBR^N$,  there holds  
  \begin{enumerate} 
  \item[$(i)$\ ] {\it    If $p\in\big(0, 1+\frac{2s}{N}\big]$, then $v_{0,\infty}(x)=\infty$ for all $x\in\BBR^N\setminus\{0\}$;}
  \item[$(ii)$] {\it    If $ p\in\big( 1+\frac{2s}{N}, \, \frac{N}{N-2s})$, then  
$
 v_{0,\infty}(x)=U_p(x)=c_p|x|^{-\frac{2s}{p-1}},\ \ x\in\BBR^N\setminus\{0\}.
$
}
    \end{enumerate}

Assume that $B_1\subset \Omega$,
then  it follows by the comparison principle that for any $k>0$,
\bel{ABC}
v_{1,k}\leq v_{\Gw,k}\leq v_{0,k} \ \ {\rm in} \ \, \BBR^N\setminus\{0\}.
\ee
Moreover, we see that $v_{1,k}$ is radially symmetric and decreasing with respect to $|x|$, and 
$$\lim_{x\to0}v_{1,k}(x)|x|^{N-2s}=c_N k.$$

\nind {\it Step 1:  we prove that $v_{\Gw,\infty}$ is infinite everywhere in $\Omega$ when $p\in(0,1+\frac{2s}{N}]$. }\\
Since the mapping $k\mapsto v_{1,k}$ is monotone, if there is one point $x_0\not\in B_1\setminus\{0\}$ such that  $\displaystyle\lim_{k\to+\infty} v_{1,k}(x_0)$ 
is finite,, then up to changing $x_0$ we can assume that $\displaystyle\lim_{k\to+\infty} v_{1,k}(z)<\infty$ for all $z$ such that $|z|\geq |x_0|-\ge$ for some 
$\ge>0$. Hence 
$$\BA{lll}\displaystyle v^p_{1,k}(x_0)=c_{N,s}\int_{|y|<|x_0|-\ge}\frac{v_{1,k}(y)-v_{1,k}(x_0)}{|y-x_0|^{N+2s}}dy+c_{N,s}P.V.\int_{|y|>|x_0|-\ge}\frac{v_{1,k}(y)-v_{1,k}(x_0)}{|y-x_0|^{N+2s}}dy\\[3.5mm]
\phantom{v^p_{1,k}(x_0)}=A_k+B_k.
\EA$$
Since $v_{1,k}(y)$ is uniformly bounded in $B_{1}\setminus B_{|x_0|-\ge}$, the term $B_k$ remains bounded too. If there exists $y_0\in B_{|x_0|-\ge}\setminus\{0\}$ such that 
$\displaystyle\lim_{k\to+\infty} v_{1,k}(y_0)=\infty$ then $\displaystyle\lim_{k\to+\infty} v_{1,k}(y)=\infty$ for all $y\in B_{y_0}\setminus\{0\}$, hence 
$A_k\to\infty$ which contradicts the fact that  $v^p_{1,k}(x_0)$ is uniformly bounded. Consequently, \\
(i) either $\displaystyle\lim_{k\to+\infty} v_{1,k}(x)=\infty$ for all $x\in B_1\setminus\{0\}$,\\
(ii) or $\displaystyle\lim_{k\to+\infty} v_{1,k}(x):=\displaystyle v_{1,\infty}(x)<\infty$ for all $x\in B_1\setminus\{0\}$.\\
Let us assume that (ii) occurs. Then $v_{1,\infty}$ is a classical solution of 
\begin{equation}\label{eq 1.1}
 \arraycolsep=1pt
\left.\begin{array}{lll}
 (-\Delta)^s v+ v^p=0\ \ &  {\rm in}\ \, B_1\setminus\{0\}\\[2mm]
 \phantom{   (-\Delta)^s  +  u^p  }
\displaystyle v=0 & {\rm in}\ \,   B_1^c.
\end{array}\right.
\end{equation}
Let $\eta_0:\BBR^N\to[0,1]$ be a smooth, radially symmetric function such that 
$$\eta_0=1\quad {\rm in}\ \ B_\frac12\quad {\rm and}\quad \eta_0=0\quad {\rm in}\ \  B_\frac34^c.$$
Set 
$$\tilde v=v_{1,k} \eta_0\quad\text{for }\; k\leq \infty.$$
A straightforward calculation shows that there exists $c_{22}>0$ independent of $k$, such that 
  $$|(-\Delta)^s \tilde v(x)+\tilde v ^p(x)|\leq c_{22}(1+|x|)^{-N-2s}\quad {\rm for\ all}\ \ x\in\BBR^N\setminus\{0\}. $$
Therefore there exists $l_0=l_0(N,s,p)>0$ such that 
$$(-\Delta)^s(2\tilde v+l_0  )+(2\tilde v+l_0  )^p>0\quad {\rm in}\quad \BBR^N\setminus\{0\}.$$
Hence  $2\tilde v+l_0 $ is a positive super solution of
$$ (-\Delta)^s  v+ v^p=0\quad {\rm in}\quad \BBR^N\setminus\{0\}.$$
From $(\cite{CV1})$ we have that for any finite $k>0$
 $$\lim_{x\to0}\frac{v_{0,k}(x)}{v_{1,k}(x)}=1\quad{\rm and}\quad \lim_{|x|\to+\infty} v_{0,k}(x)=0.$$
The comparison principle applies and we obtain that for any $k>0$, there holds
$$v_{0,k}\leq 2\tilde v+l_0 \quad {\rm in}\ \ \BBR^N\setminus\{0\}. $$
Since $p\in(0,1+\frac{2s}{N}]$ it is proved in $(\cite{CV1})$ that $\displaystyle \lim_{k\to+\infty} v_{0,k}(x)=\infty$ for all 
$x\in\BBR^N\setminus\{0\}$, which is a contradiction. Therefore $v_{1,\infty}\equiv \infty$. \\
Since $v_{1,k}\leq v_{\Gw,k}\leq v_{0,k}$ we infer that $v_{\Gw,\infty}\equiv \infty$.\medskip

\nind {\it Step 2: Asymptotics of $v_{\Gw,\infty}$ for $1+\frac{2s}{N}<p< \frac{N}{N-2s}$.}\\
If $1+\frac{2s}{N}<p<\frac{N}{N-2s}$, it is proved in  $(\cite{CV1})$  that
 $$\lim_{k\to+\infty}v_{0,k}=U_p.$$  
 Let 
 $H_s$ be the $s$-harmonic extension of $U_p$ in $B_1$, i.e. 
 $$  
 (-\Delta)^s  H_s=0\ \   \text{in}\ \,  {B}_1,\qquad
H_s=U_p\ \   \text{in}\ \,     B_1^c.
$$
Then $H_s$ is bounded   and positive  in $B_1$.  
Note that 
$H_s\geq v_{0,k}$ in $B_1^c$ for any $k>0$.
We define $W_k$ in $\BBR^N\setminus\{0\}$ by
$$W_k=v_{0,k}-H_s.$$
Then there exists a sequence  $\{r_k\}\subset(0,1)$ such that $r_k\to1$
and $W_k(x)>0$ for $0<|x|<r_k$. Actually, this is due to the fact that for any $k>0$
$$v_{0,k}= c_Nk |x|^{2s-N}(1+o(1))\quad\text{as }x\to 0\quad{\rm and}\quad  \lim_{k\to+\infty}v_k=0.$$
 Moreover, $W_k$ is a subsolution of $(\ref{eq 1.1})$. 
 Therefore, comparison principle implies that for any $k>0$
  $$v_{\Gw,k}\geq v_{1,k}\geq W_k\quad {\rm in}\ \, B_1\setminus\{0\}.$$
  Passing to the limit as $k\to+\infty$, we have that 
 $$v_{\Gw,\infty}\geq v_{1,\infty}\geq U_p-H_s\quad {\rm in}\ \, B_1\setminus\{0\},$$
and from $(\ref{ABC})$ letting $k\to+\infty$, we obtain that  $v_{\Gw,\infty}\leq v_{0,\infty}=U_p$ in   $B_1\setminus\{0\}$.
   
   As a consequence, we have that 
$$U_p-\sup_{\BBR^N} H_s\leq v_{\Gw,\infty}\leq U_p\quad {\rm in}\ \, B_1\setminus\{0\}, $$
 which also implies $(\ref{DD7})$. 
\hfill$\Box$ \medskip

The next result leads to an integral criteria for characterizing positive solutions of Emden-Fowler equations with strong singularity without  {\it a priori} estimate.

\bprop{StrSing-EF} Assume $p\in(0,\frac{N}{N-2s})$ and $B_1\subset \Gw$. If $v$ is a positive solution of $(\ref{DD6})$ such that 
\bel{DD11}\BA{lll}
\displaystyle
\myint{B_1}{}v^p dx=+\infty,
\EA\ee
then   $v\geq v_{1,k}$ in $B_1\setminus\{0\}$ for any $k>0$,  where $v_{1,k}$ is the solution of $(\ref{DD6})$. As a consequence, $v\geq v_{1,\infty}$.
\es
\Proof Let $\ga>1$, then $w_k=\ga^{-\frac1{p-1}}v_{1,k}$ satisfies
\bel{W0}\left.\BA{lll}
(-\Gd)^sw_k+\ga w_k^p=\ga^{-\frac1{p-1}}k\gd_0\ \ &\mbox{in }\, B_1\\[1.5mm]
\phantom{(-\Gd)^s+\gth w_k^p}
w_k=0&\mbox{in }\, B_1^c.
\EA\right.\ee
The equation satisfied by $v$ can be written under the form 
$$\left.\BA{lll}
(-\Gd)^sv+\ga v^p=(\ga-1)v^p\ \ &\mbox{in }\, \Gw\setminus\{0\}\\[1mm]
\phantom{(-\Gd)^s+\ga v^p}
v=0&\mbox{in }\, \Gw^c.
\EA\right.$$
From $(\ref{DD11})$ and regularity of $v$ in $\Gw\setminus\{0\}$,  for any $\vge\in (0,1)$ there holds
$$
\int_{B_\vge}v^pdx=+\infty.
$$
Hence for any $k>0$, there exists $\ell=\ell(\vge)>0$ such that 
$$
(\ga-1)\int_{B_\vge}(\min\{v,\ell\})^pdx= k.
$$
 Let $\tilde v_\vge$ be the unique solution  of 
$$\left.\BA{lll}
(-\Gd)^sv+\ga v^p=(\ga-1)\big(\min\{v,\ell\}\big)^p\chi_{_{B_\vge}}\ \ &\mbox{in }\, B_1\\[1.5mm]
\phantom{(-\Gd)^s+\gth v^p}
v=0&\mbox{in }B^c_1.
\EA\right.$$
 Note that $\tilde v_\vge$ is bounded in $B_1$. 

The comparison principle (see  \rcor{teo 2.1-com}) shows that for any $l>0$,
$$v\geq  \tilde v_\vge\quad{\rm in}\ \,  \BBR^{N}\setminus\{0\}. $$
Since $\ell_\vge\to\infty$ when $\vge\to 0$ and 
$$(\ga-1)(\min\{v,\ell\})^p\chi_{_{B_\vge}}\rightharpoonup \ga^{-\frac1{p-1}} k\gd_0,$$
  then it follows by \rprop{Dirac-1} that $\tilde v_\vge $ converges to the solution $ w_{k}$ of $(\ref{W0})$. 
Thus, we deduce that $v\geq \ga^{-\frac1{p-1}} v_{1,k}$.   Since $k>0$ is arbitrary,   we infer the claim.\hfill$\Box$ \medskip


In the following we establish a local integral {\it a priori} estimate 
\bprop{LAE} Let $p>1+\frac{2s}{N}$ and $v$ be a nonnegative solution of $(\ref{DD6})$.
 Then for any $1\leq q<\frac{p}{1+\frac{2s}{N}}$, there holds
$$
\Big(\int_{B_\frac12}v^qdx\Big)^{\frac1q}\leq c_{23}+c_{24}\int_{\Gw\setminus B_\frac12}vdx,
$$
where $c_{23}$ and $c_{24}$ are positive constants depending on $p$, $q$, $s$ and $N$.
\es
\Proof We use the same test function $\gz$ as in the proof of \rprop{estLE-2} and we set $\tilde\gz_n(x)=\gz(2^nx)$ and $\Gg_n=B_{2^{-n-1}}\setminus B_{2^{-n-2}}$,
$$\myint{\Gg_n}{}v^pdx\leq \myint{\Gg_n}{}v^p\tilde\gz_ndx=\left |\myint{\Gw}{}(-\Gd)^s\tilde\gz_n v dx\right|\leq 2^{2ns}\norm{(-\Gd)^s\gz}_{L^\infty(\BBR^N)}\norm v_{L^1(\Gw)}.
$$
Then
$$\myint{\Gg_n}{}v^qdx\leq \Big(\myint{\Gg_n}{}v^pdx\Big)^{\frac{q}{p}}\Big(\myint{\Gg_n}{}dx\Big)^{1-\frac{q}{p}}
\leq c_{25}2^{-\frac{ N}{p}\big(p-q(1+\frac{2s}{N})\big)n}\norm{(-\Gd)^s\gz}^\frac{q}{p}_{L^\infty(\BBR^N)}\norm v^\frac{q}{p}_{L^1(\Gw)}.
$$
Since $1\leq q<\frac{p}{1+\frac{2s}{N}}$, we set $\ga_0=\frac{N}{p}\big(p-q(1+\frac{2s}{N})\big)>0$ by our assumption that $q<\frac{p}{1+\frac{2s}{N}}$. Summing the above inequalities from 
$n=0$ to $+\infty$, we derive
$$
\int_{B_\frac12}v^qdx\leq \frac{2^{\ga_0} c_{25}}{2^{\ga_0}-1}\norm{(-\Gd)^s\gz}^\frac{q}{p}_{L^\infty(\BBR^N)}\norm v^\frac{q}{p}_{L^1(\Gw)},
$$
therefore 
\bel{W7}\BA {lll}\displaystyle
\Big(\int_{B_\frac12}v^qdx\Big)^{\frac pq}\leq \big(\frac{2^{\ga_0} c_{25}}{2^{\ga_0}-1}\big)^{\frac pq}\norm{(-\Gd)^s\gz}_{L^\infty(\BBR^N)}\norm v_{L^1(\Gw)}\\
\phantom{\displaystyle \left(\int_{B_\frac12}v^qdx\right)^{\frac pq}}
\leq  c_{26}\norm v_{L^1(\Gw\setminus B_{\frac12})}
+c_{27}\Big(\int_{B_\frac12}v^qdx\Big)^{\frac 1q},
\EA\ee
where $c_{26}$ and  $c_{27}$ are positive constants depending on $N$, $p$, $q$ and $s$. If we set 
$X=\norm v^\frac{q}{p}_{L^q(B_\frac 12)}$, $(\ref{W7})$ becomes
$$
X^p\leq c_{26}\norm v_{L^1(\Gw\setminus B_{12})}+c_{27}X.
$$
Standard algebraic computations combined with a convexity argument lead to
$$
\norm v_{L^q(B_\frac 12)}\leq c_{27}^\frac{1}{p-1}+\frac{c_{26}}{(p-1)c_{27}}\norm v_{L^1(\Gw\setminus B_{\frac12})}.
$$
We complete the proof. \qeda

\subsection{A priori estimate for Emden-Fowler equation}

Without loss of generality, we can  assume that $\overline B_1\subset\Gw\subset\overline \Gw\subset B_\ell$ for some $\ell>1$. The main estimate is the following
\bprop{estLE} Assume that $s\in(0,1)$, $p>1+\frac{N}{2s}$ and $v$ is a nonnegative function satisfying $(\ref{DD6})$.
Then there exists a constant $c_{28}=c_{28}(N,p,s)>0$ such that 
\bel{YY-2}\BA{lll}
v(x)\leq c_{28} |x|^{-\frac{2s}{p-1}}\quad\text{for all }\;x\in B_\frac{1}{2}\setminus\{0\}.
\EA\ee
\es

\nind\Proof The proof follows the technique introduced by Gidas and Spruck in \cite {GS} dealing with the classical Lane-Emden equation and extended to the fractional Lane-Emden equation in  \cite{JLX} and \cite{YZ}. Only the final argument is different and specific to the Emden-Fowler equation. 
We assume by contradiction that there exists a sequence $\{x_k\}\subset B_1\setminus\{0\}$ converging to $0$ such that
$$
|x_k|^{\frac{2s}{p-1}}v(x_k)=\max_{|x_k|\leq |x|\leq 1}|x|^{\frac{2s}{p-1}}v(x)\to+\infty\quad\text{as }\,k\to+\infty.
$$
We set 
$$\displaystyle
\phi_k(x)=\Big(\frac{|x_k|}{2}-|x-x_k|\Big)^{\frac{2s}{p-1}}v(x),\;\ \phi_k(\bar x_k)=\max_{|x-x_k|\leq \frac{|x_k|}{2}}\phi_k(x).
$$
Since 
$$(2\gn_k)^{\frac{2s}{p-1}}v(\bar x_k)=\gf_k(\bar x_k)\geq \gf_k(x)\geq \gn^{\frac{2s}{p-1}}_k v(x)\quad\text{for all }\,x\in B_{\gn_k}(\bar x_k),
$$
which implies
$$
2^{\frac{2s}{p-1}}v(\bar x_k)\geq v(x)\quad\text{if }|x-\bar x_k|\leq\gn_k.
$$
Since $\frac 12|\bar x_k|\geq \frac 12|x_k|-|\bar x_k-x_k|$, we obtain 
$$
|\bar x_k|^{\frac{2s}{p-1}}v(\bar x_k)\geq |x_k|^{\frac{2s}{p-1}}v(x_k)\geq |x|^{\frac{2s}{p-1}}v(x)\quad\text{for all  }\,x\in B_1\setminus B_{|x_k|}
$$
and
 $$
(2\gn_k)^{\frac{2s}{p-1}}v(\bar x_k)=\gf(\bar x_k)\geq \gf( x_k)\geq \Big(\myfrac{|x_k|}{2}\Big)^{\frac{2s}{p-1}}v(x_k)\to+\infty\quad\text{as }k\to+\infty.
$$
Therefore,
$$
w_k(y)\leq 2^{\frac{2s}{p-1}} \quad\text{for all }|y|\leq r_k:=(v(\bar x_k))^{\frac{p-1}{2s}}\gn_k.
$$
As in \cite{YZ} we define the function $W_k$ by
$$
W_k(y,z)=\frac{1}{v(\bar x_k)}v\Big(\frac{y}{(v(\bar x_k))^{\frac{p-1}{2s}}}+\bar x_k,\frac{z}{(v(\bar x_k))^{\frac{p-1}{2s}}}\Big)\quad\text{for all }(y,z)\in{\bf \Gw}_k\setminus\{(\xi_k,0)\},
$$
where 
$${\bf \Gw}_k=\Big\{(y,z)\in\BBR^{N+1}_+:\Big(\frac{y}{(v(\bar x_k))^{\frac{p-1}{2s}}}+\bar x_k,\frac{z}{(v(\bar x_k))^{\frac{p-1}{2s}}}\Big)\in\Gw\ti(0,1)\Big\} 
$$
and 
$$\xi_k= -(v(\bar x_k))^{\frac{p-1}{2s}}\bar x_k.
$$
We set $w_k =W_k(\cdot,0)$ and $\CO_k=\overline{\bf \Gw}_k\cap\prt\BBR^{N+1}_+$. The function $W_k$ satisfies  
$$
\left.\BA{lll}\displaystyle
\qquad\qquad\qquad \quad\ \CD_s W_k=0\ \ &\text{in }{\bf \Gw}_k\setminus\{(\xi_k,0)\}\\[2mm]
\prt_{\gn^s}W_k(\cdot,0)+W_k(\cdot,0)^p =0&\text{in }\CO_k\setminus\{(\xi_k,0)\}\\[2mm]
\qquad\qquad\qquad\ \, W_k(\cdot,0)=0&\text{in }\CO^c_k,
\EA\right.$$
where $\CO^c_k$ denotes the complement of $\CO_k$ in $\prt\BBR^{N+1}_+$. Since $w_k^{p-1}\leq 2^{2s}$ in $B_{r_k}$, there exists a constant 
$c_{29}>0$ depending on $N,p$ and $s$ such that
\bel{WZ2}\BA{lll}\displaystyle
0\leq W_k(y,z)\leq c_{29}\qquad\text{for all }(y,z)\in B_{\frac{r_k}{2}}\ti [0,\tfrac{r_k}{2}].
\EA\ee
It is a consequence of the  Harnack inequality \cite[Proposition 2.2]{YZ} and the invariance of the problem by scaling. From this inequality follows a local regularity estimate of $W_k$, independent of $k$ already obtained in \cite{YZ} that we recall 
$$
\norm {W_k}_{W^{1,2}(B_R\ti (0,R),z^{1-2s})}+\norm {W_k}_{C^\ga(\overline B_R\ti [0,R])}+\norm {w_k}_{C^{2,\ga}(\overline B_R)}\leq C(R) 
$$
for some $\ga>0$ and any $R>0$, where  $W^{1,2}(B_R\ti (0,R),z^{1-2s})$ denotes the weighted Sobolev space with weight $z^{1-2s}dxdz$. Therefore 
there exists a subsequence of $\{W_k\}$ still denoted by $\{W_k\}$ and a nonnegative bounded function $W$ such that $W_k\to W$ weakly in 
$W^{1,2}(\BBR^{N+1}_+,z^{1-2s})$ and in $C_{loc}^s(\overline{\BBR^{N+1}_+})$ and such that $w_k\to w:=W(\cdot,0)\in C_{loc}^2(\BBR^N)$ . Therefore the functions $W$, $w$ satisfy $w(0)=1$ and, because of $(\ref{WZ2})$, it is a nonnegative and bounded solution  of 
$$
\left.\BA{lll}\displaystyle
\qquad\qquad \qquad \ \ \  \CD_s W=0\ \ &\text{in }\,\BBR^{N+1}_+\\[2mm]
 \prt_{\gn^s}W(\cdot,0)+W(\cdot,0)^p=0&\text{in }\, \BBR^{N}.
\EA\right.$$
By \rth{trace} $w$ belongs to $L^1_{\gm_s}(\BBR^N)$, thus $\prt_{\gn^s}W(\cdot,0)=(-\Gd)^sw$ and $w$  satisfies $(\ref{Li1})$. By \rcor{Liouv} such a function cannot exist. This implies that $(\ref{YY-2})$ holds.\qeda\medskip

The next statement is a surprising consequence of \rprop{estLE}.
\bcor{estCOR-2} 
Assume that $1+\frac{2s}{N}< p<\frac{N}{N-2s}$. If $v$ is a positive function satisfying $(\ref{DD6})$ and
$$
\int_{\Omega}v^pdx=+\infty.
$$
Then there holds
$$
\liminf_{x\to 0}|x|^{\frac{2s}{p-1}}v(x)=c_p.
 $$
 \es
\Proof By \rprop{Dirac-2} and \rprop{StrSing-EF} we have
$$
\liminf_{x\to 0}|x|^{\frac{2s}{p-1}}v(x)\geq c_p.
$$
Therefore we proceed by contradiction assuming that there exists a nonnegative solution  $v$  and $\bar\gm>0$
such that 
 \begin{equation}\label{ll 3.0-x2}
c_p<\bar\gm^{1/(p-1)}<\liminf_{x\to 0}|x|^{\frac{2s}{p-1}}v(x)<+\infty.
 \end{equation}
We can write the equation satisfied by $v$ under the form 
$$
\CL_{\bar \mu}^s v=g(x):=\frac{\bar\mu}{|x|^{2s}}v-v^p.
 $$
 The mapping $\gm\mapsto\gt_-(\gm)$ is continuous and monotone decreasing from $[\gm_0,+\infty)$ onto $(-\infty,\frac{2s-N}{2}]$ and 
 $\gt_-(c_p^{p-1})=-\frac{2s}{p-1}$. Therefore, there exists $\bar\gt>\frac{2s}{p-1}$ such that $\gt_-(\bar\gm)=-\bar\gt$. By $(\ref{YY-2})$ we have 
 $$\liminf_{x\to 0}|x|^{\gt}v(x)=0
 $$
 for any $\gt>\frac{2s}{p-1}$. Choosing $\gt=\bar\gt$ we see that there exist $r,c_{30}>0$ such that
 $$
    0<-g(x) \leq c_{30} |x|^{-p\bar\gt}\qquad \text{for any }\, 0<|x|<r.
$$
Set $\displaystyle m=r^{-\gt_+(\bar\gm)}v(r)$. For any $\vge>0$ the function $w_{\vge,m}:=\vge\Phi_{\bar\gm}+m\Gamma_{\bar\gm}$ satisfies 
$\CL_{\bar \mu}^sw_{\vge,m}=0$. Therefore the function $w_{\vge,m}-u$ satisfies in turn 
$$\CL_{\bar \mu}^s(w_{\vge,m}-v)=-g\geq 0\qquad\text{in }\; B_r\setminus\{0\}.$$
This function is nonnegative in 
$\BBR^N\setminus B_r$ and there holds
$$\liminf_{x\to 0}\frac{(w_{\vge,m}-u)(x)}{\Phi_{\bar\gm}(x)}=\vge>0.
$$
By Lemma \ref{L-comp}, $w_{\vge,m}-v\geq 0$. Letting $\vge\to 0$ yields $v(x)\leq m\Gamma_{\bar\gm}(x)$. Thus $(\ref{ll 3.0-x2})$ does not hold and no such  $\bar\gm$ exists, which implies the claim.\qeda\medskip

\mysection{Singularity of solutions }

\subsection{Energy method and limit set}
The estimates $(\ref{WW-2})$ and $(\ref{YY-2})$ combined with the lifting estimate $(\ref{S17})$ allow to use an energy method to describe the behaviour of signed solutions of $(\ref{X-2})$. Furthermore, this energy method does not depend of the signs of solutions, even if the bounds do. If we set
\bel{I17}\BA{lll}\displaystyle
u(\xi)=u(r,\gs)=r^{-\frac{2s}{p-1}}w(t,\gs)\quad\text{with }\,t=\ln r
\EA\ee
and $w(t,\gs)=w(t,\gs',\gf)$. Then $w$ satisfies 
$$
\left.\BA{lll}\displaystyle
w_{tt}+\Gth_{s,p,{\scriptscriptstyle N}}w_t+\Gl_{s,p,{\scriptscriptstyle N}}w+\CA_s[w]=0\ \, &\text{in }\,\BBR\ti \BBS^N_+\\[2mm]
\phantom{------,,-\ }
\myfrac{\prt w}{\prt\gn^s}+\ge |w|^{p-1}w=0\ \, &\text{in }\,\BBR\ti \BBS^{N-1},
\EA\right.$$
where $\Gth_{s,p,{\scriptscriptstyle N}}$ and $\Gl_{s,p,{\scriptscriptstyle N}}$ are defined in $(\ref{X-7})$ and 
$$\displaystyle\myfrac{\prt w}{\prt\gn^s}(t,\gs')=-\lim_{\phi\to 0}\big(\sin\gf)^{1-2s} w_\gf(t,\gs',\gf)\big).
$$
We define the negative (resp. positive) trajectory of $w$ by 
$$
\CT_-[w]=\bigcup \big\{w(t,\cdot), t\leq 0\big\},
$$
respectively
$$
\CT_+[w]=\bigcup \big\{w(t,\cdot), t\geq 0\big\}.
$$
The trajectory of $w$ denoted by $\CT[w]=\CT_-[w]\cup\CT_+[w]$.
We also define the  $\ga$-limit set $\Gg_-[w]$ and the $\omega$-limit set $\Gg_+[w]$ of a trajectory by
$$
\BA{lll}\displaystyle
\Gg_-[w]=\bigcap_{t\leq 0}\left(\overline{\bigcup \left\{w(\gt,\cdot), \gt\leq t\right\}}^{\bf X}\right)\quad\text{and }\;
\Gg_+[w]=\bigcap_{t\geq 0}\left(\overline{\bigcup \left\{w(\gt,\cdot),\gt\geq t\right\}}^{\bf X}\right),
\EA$$
where ${\bf X}=\Big\{\gz\in C^s(\BBS^{N}_+):\gz(\cdot,0)\in C^2(\BBS^{N-1})\Big\}$. If we assume that the set $\{w(t,\cdot)\}$ is bounded in $\BBR_-\ti \BBS^N_+$ (resp. in $\BBR_+\ti \BBS^N_+$), its trajectory $\CT_-[w]$ (resp. $\CT_+[w]$) is uniformly bounded in  $\left\{\gz\in C^\ga(\BBS^{N}_+):\gz(\cdot,0)\in C^{2+\ga}(\BBS^{N-1})\right\}$. By Ascoli's theorem they are relatively compact in ${\bf X}$.
By the general theory of limit sets, $\Gg_\pm[w]$ are compact and connected subsets of ${\bf X}$. We first give a general result about these limit sets under the assumption that the solutions are signed and satisfy {\it a priori} estimate which is indeed only proved for nonnegative solutions.

\bth{3} Assume that $s\in (0,1)$, $p\in (1,+\infty)\setminus\{\frac{N+2s}{N-2s}\}$, $\ge=\pm 1$ and $u\in C(\overline{\BBR^{N+1}_+}\setminus\{(0,0)\})\cap C^2(\BBR_+^{N+1})$
is a solution of $(\ref{X-2})$. If $u$ satisfies 
\bel{I21}\BA{lll}\displaystyle
|u(x,z)|\leq c_{31}\left(|x|^2+z^2\right)^{-\frac{s}{p-1}}\quad\text{for }\,0<|x|<1,\quad(\text{resp. for  }|x|>1)
\EA\ee
for some $c_{31}>0$,
then $\Gg_-[w]$ (resp. $\Gg_+[w]$) is a nonempty compact connected subset of the set $\CE_\ge$, which is the set of solutions of $(\ref{X-5})$. 
\es
\Proof We give the proof for $\Gg_-[w]$, the case  $\Gg_+[w]$ being the same up to changing $t$ in $-t$. We define the energy function
$$
\CI_\ge[w](t)=\frac12\int_{\BBS^{N}_+}\left(|\nabla' w|^2-\Gl_{s,p,{\scriptscriptstyle N}}w^2-w_t^2\right)d\gm_s-\frac\ge{p+1}\int_{\BBS^{N-1}}|w|^p dS'.
$$
The function $w$ is bounded in $C^2(\overline{\BBS^N_+}\ti (-\infty,-1])$ and there holds
 \bel{AA2}\displaystyle
\frac{d}{dt}\CI_\ge[w](t)=\left(p-\frac{2s(p+1)}{p-1}\right)\int_{\BBS^{N}_+}w^2_td\gm_s=\Gth_{s,p,{\scriptscriptstyle N}}\int_{\BBS^{N}_+}w^2_td\gm_s.
 \ee
Since $p\neq\frac{N+2s}{N-2s}$, $\Gth_{s,p,{\scriptscriptstyle N}}\neq 0$, and we have the damping estimate
$$
\int_{\BBS^{N}_+}w^2_td\gm_s<+\infty.
$$
We conclude by using the uniform continuity of $w_t$ that $w_t(t,\cdot)\to 0$ uniformly on $\overline{\BBS^N_+}$ when $t\to-\infty$. Let $\gw\in\Gg_-$. There exists 
a sequence $\{t_n\}$ converging to $-\infty$ such that $w(t_n,\cdot)\to w$ uniformly on $\overline{\BBS^N_+}$, and since $w_t(t,\cdot)\to 0$ uniformly on $\overline{\BBS^N_+}$, for any $T>0$, $w(t,\cdot)\to w$ uniformly on $[t_n-T,t_n+T]\ti \overline{\BBS^N_+}$. Let $\gz\in C^{2}(\overline{\BBS^N_+})$ such that $\prt_{\gn^s}\gz=0$, then 
$$
\BA{lll}\displaystyle
\int_{\BBS^N_+}\left(w_t(T+t_n,\cdot)-w_t(T-t_n,\cdot)\right)\gz d\gm_s +\Gth_{s,p,{\scriptscriptstyle N}}\int_{\BBS^N_+}\left(w(T+t_n,\cdot)-w(T-t_n,\cdot)\right)\gz d\gm_s\\[4mm]\phantom{-------\ }
\displaystyle
-\Gl_{s,p,{\scriptscriptstyle N}}\int_{-T}^{T}\int_{\BBS^N_+}w(t_n+t,\cdot)\gz d\gm_sdt
-\int_{-T}^{T}\int_{\BBS^N_+}w(t_n+t,\cdot)\CA_s[\gz] d\gm_sdt\\[4mm]\phantom{--------------------}
\displaystyle
+\ge\int_{-T}^{T}\int_{\BBS^{N-1}}(|w|^{p-1}w)(t_n+t,\cdot)\gz dS'dt =0.
\EA$$
Letting $t_n\to-\infty$ yields
$$
\int_{-T}^{T}\int_{\BBS^N_+}\left(\Gl_{s,p,{\scriptscriptstyle N}}\gw\gz+\gw\CA_s[\gz]\right)d\gm_sdt-\ge\int_{-T}^{T}\int_{\BBS^{N-1}}|\gw|^{p-1}\gw\gz dS'dt=0.
$$
This implies that $\gw$ is a weak solution of $(\ref{X-5})$ and therefore a strong one. \qeda\medskip

\nind\BBRemark The case $p=\frac{N+2s}{N-2s}$ is not covered since the damping coefficient $\Gth_{s,p,{\scriptscriptstyle N}}$ in $(\ref{AA2})$ vanishes. 

\subsection{Singularities of    Emden-Fowler equation}

The next result concerns the singularities of solutions to the   Emden-Fowler equations.

\bth{3.1} Assume that  $s\in (0,1)$, $\ge=1$ and $u\in C^1\big(\overline{\BBR^{N+1}_+}\setminus\{(0,0)\}\big)\cap C^2(\BBR_+^{N+1})$
is a solution of $(\ref{X-2})$ satisfying $(\ref{I21})$, $v=u(\cdot,0)$ and $w$ is defined by $(\ref{I17})$. The following convergences hold when $t\to-\infty$.\smallskip

\nind 1-a If $p\geq \frac{N}{N-2s}$, then $w(t,\cdot)$ converges uniformly in $\BBS^{N}_+$  to $0$, $v\equiv0$, hence $u$ is a smooth solution;\smallskip

\nind 1-b If $1+\frac{2s}{N}<p<\frac{N}{N-2s}$ and $u$ is nonnegative, then either $w(t,\cdot)$ converges uniformly in $\BBS^{N}_+$  to $\gw_1$ or it converges to $0$;\smallskip

\nind 1-c If $p^*\leq p<\frac{N}{N-2s}$, then $w(t,\cdot)$ converges uniformly in $\BBS^{N}_+$  to some $\ell\in\{0,\gw_1,-\gw_1\}$;\smallskip 

\nind 1-d If $1<p\leq 1+\frac{2s}{N}$, and $u$ is nonnegative,  then $w(t,\cdot)$ converges uniformly in $\BBS^{N}_+$  to $0$.
\smallskip

\nind 2. Let $1<p< \frac{N}{N-2s}$,  $u$ is nonnegative,  if $w(t,\cdot)$ converges uniformly in $\BBS^{N}_+$ to $0$, then there exists $k\in \BBR_+$ such that 
$$
e^{(N-\frac{2s}{p-1})t}w(t,\gs',\gf)\to k\sin^{2s}(\gf)\quad\text{uniformly in }\,\BBS^{N}_+.
$$
If $k=0$, then $w\equiv0$   and $u$ is smooth.
\es
\Proof {\it 1-a.} If $p\geq \frac{N}{N-2s}$, the set $\CE_1$ is reduced to the zero function. Therefore $w(t,\cdot)$ converges to $0$ when $t\to -\infty$. By $(\ref{X-3})$ the function $\tilde u(t,\gs)=u(r,\gs)$ with $t=\ln r$ satisfies 
 \bel{AA6}
\left.\BA{lll}
\tilde u_{tt}+(N-2s)\tilde u_t+\CA_s[\tilde u]=0\ \, &\text{in }(-\infty,0)\ti \BBS^N_+\\[2mm]
\phantom{----,,\ }
\myfrac{\prt \tilde u}{\prt\gn^s}+|\tilde u|^{p-1}\tilde u=0&\text{in }(-\infty,0)\ti \BBS^{N-1}.
\EA\right. \ee
Since $\frac{2s}{p-1}\leq N-2s$, we have that $\tilde u(t,\cdot)=o(e^{(2s-N)t})$ when $t\to-\infty$. Let $m=\max |\tilde u(0,\cdot)|$. For any $\vge>0$, 
$t\mapsto \vge e^{(2s-N)t}+m$ is a supersolution of $(\ref{AA6})$, larger than $\tilde u$ at $t=0$ and $t\to-\infty$. Hence $\tilde u(t,\gs)\leq \vge e^{(2s-N)t}+m$. 
Letting $\vge\to 0^+$ yields $\tilde u(t,\gs)\leq m$. Similarly $\tilde u(t,\gs)\geq -m$, hence $\tilde u$ is bounded and so is $u$. Therefore 
$v\ \big(=u(\cdot,0)\big)\in L^1_{\gm_s}(\BBR^N)$,
$$\myfrac{\prt \tilde u}{\prt\gn^s}=(-\Gd)^sv,
$$
and $v$ is a bounded solution of the fractional Emden-Fowler equation in $\Gw\setminus\{0\}$ vanishing on $\Gw^c$. It is therefore identically $0$. By the even extension of $u$ to $\BBR^{N+1}$, it follows that $u$ is smooth. \smallskip

\nind {\it 1-b} Since $\CE_1^+=\{\gw_1,0\}$ and is disconnected, the result follows. \smallskip

\nind {\it 1-c} Since $\CE_1=\{\gw_1,-\gw_1,0\}$ and is disconnected, the result follows.\smallskip

\nind {\it 1-d}  Since $\CE_1^+=\{0\}$ the result follows. \smallskip

\nind {\it 2.} If $u$ is nonnegative, then $v\in L^1_{\gm_s}(\BBR^N)$, thus $(\ref{DD6})$ holds.  If $(\ref{DD11})$ holds, then 
$$  v\geq v_\infty=\lim_{k\to+\infty }v_k$$
 by \rprop{StrSing-EF}, which is impossible if $w(t,\cdot)$ converges to $0$. Therefore $v\in L^p(B_1)$ and by \rlemma{lm 4.1-s}
 there exists $k\geq 0$ such that $(\ref{DD1})$ holds and therefore $v=v_k$. Furthermore, $(\ref{DD3})$ is satisfied. If $k=0$, then $v$ is identically $0$ and $u$ is smooth as in the previous case.\qeda
\subsection{Singularities for  Lane-Emden equation}
In the case of Lane-Emden equation we have the following.

\bth{3.2} Assume that  $s\in (0,1)$, $\ge=-1$, $p>1$ and $p\neq \frac{N+2s}{N-2s}$ and $u\in C(\overline{\BBR^{N+1}_+}\setminus\{(0,0)\})\cap C^2(\BBR_+^{N+1})$
is a nonnegative solution of $(\ref{I-5})$ satisfying $(\ref{I21})$. The following convergences hold when $t\to\infty$.\smallskip

\nind 1-a If $p> \frac{N}{N-2s}$, then $w(t,\cdot)$ converges uniformly in $\BBS^{N}_+$ either to $\gw_1$ or to $0$.\smallskip

\nind 1-b If $1<p\leq\frac{N}{N-2s}$,  then $w(t,\cdot)$ converges uniformly to $0$ in $\BBS^{N}_+$.

\nind 2-a  If $p> \frac{N}{N-2s}$ and  $w(t,\cdot)$ converges to $0$, then $u$ is a smooth function in $\overline{\BBR^{N+1}_+}$ and so is $v$.\smallskip

\nind 2-b  If $1<p<\frac{N}{N-2s}$, there exists $k\geq 0$ such that 
$$
e^{(N-\frac{2s}{p-1})t}w(t,\gs',\gf)\to k(\sin\gf)^{2s}\quad\text{uniformly in }\,\BBS^{N}_+.
$$
\es
\Proof Assertions {\it 1-a} and {\it 1-b} follow from \rth{3} and the description of $\CE^+_{-}$ given in \rth{Th2}. Furthermore, by \rlemma{lm 4.1-s} there exists $k\geq 0$ such that $(\ref{2.y4})$ holds with $\ge=-1$. Observe that $v(x)\geq kG_{s,\Gw}(x,0)$, where $G_{s,\Gw}$ is the Green function of $(-\Gd)^s$ in $\Gw$ subject to zero Dirichlet condition in $\Gw^c$. If $p\geq \frac{N}{N-2s}$, then $G_{s,\Gw}\notin L^p(B_1)$, therefore $k=0$ in that case.\smallskip

 \nind  {\it 2-a. Step 1:}  If $w(t,\cdot)$ converges to $0$, we first prove that there exist $\vge>0$ and $c_{32}>0$ such that 
 \bel{BB2}
0\leq  w(t,\gs)\leq c_{32}e^{\vge t}\qquad\text{for all }(t,\gs)\in (-\infty]\ti \overline {\BBS^{N}_+}.
\ee
Set $\gr(t)=\norm {w(t,.)}_{C^0(\BBS^N_+)}$. If $(\ref{BB2})$ does not hold we have that
$$
\limsup_{t\to-\infty}e^{-\vge t}\gr(t)=+\infty\quad\text{for all }\vge>0.
$$
 By \cite[Lemma 2.1]{CMV} there exists a function $\eta\in C^\infty\big((-\infty,0]\big)$ such that 
  \bel{BB4}\BA{lll}
&(i)\qquad\displaystyle &\displaystyle\eta>0,\,\eta'>0,\,\lim_{t\to-\infty}\eta(t)=0;\qquad\qquad\qquad\qquad\qquad\qquad\qquad\qquad\\[3mm]
&(ii)\qquad\displaystyle &\displaystyle0<\limsup_{t\to-\infty}\frac{\gr(t)}{\eta(t)}<+\infty;\\[3mm]
&(iii)\qquad\displaystyle &\displaystyle\lim_{t\to-\infty}e^{-\vge t}\eta(t)=+\infty\quad\text{for all }\vge>0;\\[3mm]
&(iv)\qquad\displaystyle &\displaystyle\left(\frac{\eta'}{\eta}\right)',\, \left(\frac{\eta''}{\eta}\right)'\in L^1((-\infty,0));\\[3mm]
&(v)\qquad\displaystyle &\displaystyle\lim_{t\to-\infty}\frac{\eta'(t)}{\eta(t)}=\lim_{t\to-\infty}\frac{\eta''(t)}{\eta(t)}=0.
\EA\ee
We define $\psi$ by $w(t,\cdot)=\eta(t)\psi(t,\cdot)$, then 
\bel{BB5}\left.\BA{lll}\displaystyle
 \psi_{tt}+\varrho_1\psi_t+\varrho_2\psi+\CA_s[\psi]=0\ \, &\text{in }(-\infty,0)\ti \BBS^N_+\\[2.5mm]
\phantom{----- }
\myfrac{\prt \psi}{\prt\gn^s}- \eta^{p-1}\psi^p=0\ \, &\text{in }(-\infty,0)\ti \BBS^{N-1},
\EA\right.\ee
where 
$$\varrho_1= \Gth_{s,p,{\scriptscriptstyle N}}+2\frac{\eta'}{\eta}\quad{\rm and}\quad \varrho_2= \Gl_{s,p,{\scriptscriptstyle N}}+\frac{\eta'}{\eta}\Gth_{s,p,{\scriptscriptstyle N}}+\frac{\eta''}{\eta}.$$

The function $\psi$ is uniformly bounded. Then by standard regularity theory for degenerate elliptic operators (see e.g. \cite{FKS}), $\psi$ is bounded in 
$C^2(-\infty,-1]\ti \overline{\BBS^N_+})$. Therefore, the trajectory of $\CT_-[\psi]:=\{\psi(t,\cdot): t\leq 0\}$ is relatively compact in ${\bf X}\ \, \big(\!=\big\{\gf\in C^s(\BBS^{N}_+):\gf(\cdot,0)\in C^2(\BBS^{N-1})\big\}\big)$. By $(\ref{BB4})$-(ii) there exist a nonzero nonnegative function $\gth\in C^{2}(\overline{\BBS^N_+})$ and a sequence $\{t_n\}$ tending to $-\infty$ such that $\psi(t_n,\cdot)\to \gth$ in $C^2(\overline {\BBS^N_+})$. If we multiply by $\psi_t$ the equation $(\ref{BB5})$ which is satisfied by $\psi$ and integrate over $\BBS^N_+$, we obtain the following energy formula
\bel{BB6}\BA{lll}\displaystyle
\myfrac{d}{dt}\CF[\psi](t)=\varrho_1 \int_{\BBS^N_+}\psi_t^2d\gm_s
-\frac{1}{2}\varrho_2'\int_{\BBS^N_+}\psi^2d\gm_s 
 -\frac{p-1}{p+1}\eta'\eta^{p-2}\int_{\BBS^{N-1}}\psi^{p+1} dS',
\EA\ee
where 
$$\varrho_2'=\big(\frac{\eta'}{\eta}\big)'\Gth_{s,p,{\scriptscriptstyle N}}+\big(\frac{\eta''}{\eta}\big)' $$
and
$$
\CF[\psi](t):=\frac12\int_{\BBS^{N}_+}\Big(|\nabla' w|^2- \varrho_2\Big) \psi^2d\gm_s
-\frac{\eta^{p-1}}{p+1}\int_{\BBS^{N-1}}\psi^{p+1} dS'.
$$
By properties $(i)$, $(iv)$ and $(v)$ in $(\ref{BB4})$, we have that 
$$\eta'\eta^{p-2}\ \Big(=\frac{\eta'}{\eta} \eta^{p-1}\Big), \  \varrho_2',\ \frac{\eta'(t)}{\eta(t)}\in L^1((-\infty,0)) \quad\text{as }t\to\infty.
$$
Since $\CF[\psi](t)$ is bounded independently of $t$, it follows that 
$$
\int_{-\infty}^0\int_{\BBS^N_+}\psi_t^2d\gm_sdt<+\infty
$$
using the regularity estimates it implies that
$$
\lim_{t\to-\infty}\int_{\BBS^N_+}\psi_t^2d\gm_s=0.
$$
This implies that  the convergence of $\psi(t_n,\cdot)$ to $\gth$ holds in the stronger following way: for any $T>0$, 
$$\psi(t,\cdot)\to \gth \quad\text{uniformly in } [t_n-T,t_n+T]\ti \overline{\BBS^N_+}.
$$
Let $\gz\in C^2( \overline{\BBS^N_+})$ such that $\prt_{\gn^s}\gz=0$, then from $(\ref{BB5})$, 
\bel{BB9}\BA{lll}\displaystyle
\quad\ \ \int_{\BBS^N_+}\left(w_t(T+t_n,\cdot)-w_t(T-t_n,\cdot)\right)\gz d\gm_s +\Gth_{s,p,{\scriptscriptstyle N}}\int_{\BBS^N_+}\left(w(T+t_n,\cdot)-w(T-t_n,\cdot)\right)\gz d\gm_s\\[4mm]\phantom{---}
\displaystyle
-\Gl_{s,p,{\scriptscriptstyle N}}\int_{-T}^{T}\int_{S^N_+}\psi(t_n+t,\cdot)\gz d\gm_sdt
-\int_{-T}^{T}\int_{\BBS^N_+}\psi(t_n+t,\cdot)\CA_s[\gz] d\gm_sdt\\[4mm]\phantom{ }
\displaystyle
=2\int_{-T}^{T}\int_{\BBS^N_+}\left(\frac{\eta'}{\eta}\psi_t(t_n+t,\cdot)+\varrho_2'\psi(t_n+t,\cdot)\right)\gz d\gm_s dt +
\int_{-T}^{T}\int_{\BBS^{N-1}}\eta^{p-1}\psi^p(t_n+t,\cdot)\gz dS'dt.
\EA\ee
When $t_n\to-\infty$, the right-hand side of $(\ref{BB9})$ tends to $0$. Since $\prt_{\gn^s}\gth=0$ on $\BBS^{N-1}$ from the equation satisfied by $\psi$ and the convergence of $\psi$, it follows, as in the proof of \rth{3}, that $\gth$ satisfies 
$$
\left.\BA{lll}\displaystyle
\CA_s[\gth]+\Gl_{s,p,{\scriptscriptstyle N}}\gth=0\ \ &\text{in }\, \BBS^N_+\\[2mm]
\phantom{---;;-\ }
\myfrac{\prt \gth}{\prt\gn^s}=0\ \  &\text{in }\, \BBS^{N-1}.
\EA\right.$$
As $\Gl_{s,p,{\scriptscriptstyle N}}<0 $ we derive $\gth=0$ by integration, contradiction. Hence $(\ref{BB2})$ holds for some $\vge>0$. \smallskip

\nind {\it 2-a.  Step 2.}   Next we prove that there exist $c_{33}>0$ such that
$$
w(t,\gs)\leq c_{33} e^{\frac{2s}{p-1}t}\qquad\text{for all }(t,\gs)\in (-\infty,0]\ti \overline {\BBS^N_+}.
$$
The set $\CQ$ of $\vge>0$ such that $(\ref{BB2})$ holds is a non-empty interval and it admits an upper bound $\vge^*$ which is finite if $v$ is not identically $0$. 
If $\vge\in\CQ$, we set  $\psi_{\vge}(t,\cdot)=e^{-\vge t}w(t,\cdot)$\smallskip

If $\vge^*\in\CQ$, we take $\vge=\vge^*$ and denote $\psi_{\vge^*}=\psi^*$. Then $\psi^*$ is bounded and satisfies
$$
\left.\BA{lll}\displaystyle
\!\!\!\psi^*_{tt}+\varrho_3^* \psi^*_t+\varrho_4^*\psi^*+\CA_s[\psi^*]=0\ \, &\text{in }(-\infty,0)\ti \BBS^N_+\\[3mm]
\phantom{ ---\ \,  }
\myfrac{\prt \psi^*}{\prt\gn^s}- e^{(p-1)\vge^*t}\psi^{*p}=0\ \, &\text{in }(-\infty,0)\ti \BBS^{N-1},
\EA\right.$$
where 
$$\varrho_3^*= \Gth_{s,p,{\scriptscriptstyle N}}+2\vge^* \quad{\rm and}\quad 
\varrho_4^*= \Gl_{s,p,{\scriptscriptstyle N}}+\vge^*\Gth_{s,p,{\scriptscriptstyle N}}+\vge^{*2}.$$

Because $\psi^*$ is bounded, the same analysis as in the previous cases shows that any element $\gth^*$ in the limit set of the trajectory of $\{\psi^*\}$ is a nonnegative function satisfying 
 \bel{BB13}
\left.\BA{lll}\displaystyle
\varrho_4^* \gth^*+\CA_s[\gth^*]=0\ \,  \text{in }\, \BBS^N_+,\qquad \myfrac{\prt \gth^*}{\prt\gn^s}=0\ \,  \text{in }\, \BBS^{N-1}.
\EA\right. \ee
The polynomial $P(z)=z^2+\Gth_{s,p,{\scriptscriptstyle N}}z+\Gl_{s,p,{\scriptscriptstyle N}}$ admits two roots $z_2<0<z_1$, where
$$z_1=\frac{2s}{p-1}\quad\text{and }\quad z_2=2s-N+\frac{2s}{p-1}.
$$
Notice that $0<\vge^*\leq z_1$ since $v(0)>0$. If $\vge^*<z_1$, the limit set is reduced to $\{0\}$, hence we proceed as is Step 1, by introducing a function $\eta^*$ which has the same properties as the function $\eta$ shown in $(\ref{BB4})$, except that $(\ref{BB4})$-(ii) is replaced by
$$
0<\limsup_{t\to-\infty}\frac{\gr^*(t)}{\eta^*(t)}<+\infty\quad\text{with }\,\gr^*(t)=\norm{\psi^*(t,.)}_{L^\infty(S^N_+)}.
$$
Defining the function $\gu^*=\frac{w^*}{\eta^*}$ which satisfies
$$
\left.\BA{lll}\displaystyle
\gu^*_{tt}+\varrho_5^*\gu^*_t+\varrho_6^*\gu^*+\CA[\gu^*]=0 \ \, &\text{in }(-\infty,0)\ti \BBS^N_+\\[3mm]
\phantom{\ \ }\myfrac{\prt\gu^*}{\prt\gn^s}- e^{(p-1)\vge^*t}\eta{^*}^{p-1}\gu^{*p}=0\ \, &\text{in }(-\infty,0)\ti \BBS^{N-1},
\EA\right.$$
where 
$$\varrho_5^*= \Gth_{s,p,{\scriptscriptstyle N}}+2\vge^*+2\frac{\eta{^*}'}{\eta{^*}}\quad
{\rm and}\quad \varrho_6^*= \Gl_{s,p,{\scriptscriptstyle N}}+\Big(\vge^*+\frac{\eta{^*}'}{\eta{^*}}\Big)\Gth_{s,p,{\scriptscriptstyle N}}+\vge^{*2}+\frac{\eta{^*}''}{\eta{^*}}.$$

The same analysis carried on in Step 1 shows that there exists a nonnegative nonzero element in the limit set of the trajectory of $\gu^*$ which satisfies the same equation as $\gth^*$ in $(\ref{BB13})$, contradiction. Therefore $\vge^*=z_1$ which implies that $v$ is bounded and therefore regular.  

\smallskip

If $\vge^*\notin\CQ$, then for any $\vge<\vge^*$ the function $\psi_\vge(t,.)=e^{-\vge t}w(t,\cdot)$ is bounded, tends to $0$ when $t\to-\infty$ and satisfies 
$$
\left.\BA{lll}\displaystyle
\!\!\!\psi_{\vge\,tt}+\varrho_3\psi_{\vge\,t}+\varrho_4\psi_\vge+\CA_s[\psi_\vge]=0\ \ &\text{in }(-\infty,0)\ti \BBS^N_+\\[3mm]
\phantom{ -----    }
\myfrac{\prt \psi_\vge}{\prt\gn^s}- e^{(p-1)\vge t}\psi^p_\vge=0\ \ &\text{in }(-\infty,0)\ti \BBS^{N-1},
\EA\right.$$
where
$$\varrho_3= \Gth_{s,p,{\scriptscriptstyle N}}+2\vge  \quad{\rm and}\quad 
\varrho_4= \Gl_{s,p,{\scriptscriptstyle N}}+\vge \Gth_{s,p,{\scriptscriptstyle N}}+\vge^{2}.$$

Put $X_\vge(t)=\myint{\BBS^N_+}{}\psi_\vge(t,.)d\gm_s$ and $\displaystyle F_\vge(t)=\int_{\BBS^{N-1}}\!\!\!\psi^p_\vge dS'$, then 
$$
X''_\vge+\varrho_3 X'_\vge+ \varrho_4 X_\vge
+e^{(p-1)\vge t}F_\vge=0.
$$
The equation
$$Y''+\varrho_3 Y'+ \varrho_4Y=0
$$
admits the two linearly independent solutions
$$
Y_1(t)=e^{\left(\frac{2s}{p-1}-\vge\right)t}\quad\text{and }\;Y_2(t)=e^{\left(2s-N+\frac{2s}{p-1}-\vge\right)t}.
$$
The function $X_\vge$ satisfies 
$$-X''_\vge-\varrho_3 X'_\vge-\varrho_4X_\vge
\leq c_{34}e^{(p-1)\vge t}\quad\text{on }\,(-\infty,0].
$$
Since it tends to $0$ when $t\to-\infty$, $X_\vge$ is bounded from above by the solution $\tilde X_\vge$ of 
$$\left.\BA{lll}\displaystyle-\tilde X''_\ge-\varrho_3\tilde X'_\vge-\varrho_4\tilde X_\vge
= c_{34}e^{(p-1)\vge t}\ \ \text{on }\,(-\infty,0]\\[3mm]
\phantom{------ \ \ \, }\displaystyle\tilde X_\vge(0)=X_\ge(0), \\[3mm]
\phantom{----\,   }\displaystyle\lim_{t\to-\infty}\tilde X_\vge(t)=0.
\EA\right.
$$
The solution of this last differential equation is easily computable as we can always assume that $(p-1)\vge\neq \frac{2s}{p-1}-\vge$ (the equality is called the resonance phenomenon) and we have
$$
\tilde X_\vge(t)=X_\vge(0)e^{\left(\frac{2s}{p-1}-\vge\right)t}+c_{35}e^{(p-1)\vge t}
$$
for some explicit constant $c_{35}$ depending on the coefficients. This implies 
$$
\myint{\BBS^N_+}{}w(t,\cdot)d\gm_s\leq X_\vge(0)e^{\frac{2s}{p-1}t}+c_{33}e^{p\vge t}.
$$
Applying Harnack inequality \cite[Proposition 3.1]{YZ} we obtain that 
$$
w(t,\cdot)\leq c_{36}e^{\frac{2s}{p-1}t}+c_{37}e^{p\vge t}\quad\text{for all }t\leq T_0.
$$
Because $\vge$ can be taken arbitrarily close to $\vge^*$, we have that $\min\{p\vge,\frac{2s}{p-1}\}>\vge^*$, which contradicts the definition of $\vge^*$. Hence 
$\vge^*\in\CQ$, which ends the proof of {\it 2-a}.\smallskip

\nind {\it 2-b.} If $1<p<\frac{N}{N-2s}$, there exists $k\geq 0$ such that $(\ref{2.y4})$ holds, and 
$$kG_{s,\Gw}(x,0)\leq v(x)\leq c_{38}|x|^{2s-N},
$$
where $c_{38}>0$. \qeda\medskip

\nind{\bf Proofs of Theorem E--F. } Theorem E is \rth{3}, Theorem  F follows \rth{3.1} part {\it 2} and Theorem  E does \rth{3.2} part {\it 1-a, 2-a}.
 \qeda\medskip

\subsection{Classification of fractional Emden-Fowler equation}
\nind{\bf Proof of Theorem G. }  Let $v$ be a nontrivial nonnegative solution of  $(\ref{IN17})$.  It suffices to prove that \\{\it

\nind 1.  for $p \in\big(0, 1+\frac{2s}{N}\big]$,   
 $$v=v_k\quad\ {\rm in}\ \  B_1\setminus\{0\}\ \ \text{for some }\, k\geq0;$$
 \nind 2.   for $p\in \big(1+\frac{2s}{N}, \frac{N}{N-2s}\big)$,  either 
  $$v=v_k\quad\ {\rm in}\ \  B_1\setminus\{0\}\ \ \text{for some }\, k\geq 0  $$
  or 
  $$v=v_\infty\quad\ {\rm in}\ \  B_1\setminus\{0\}.$$
  
   \nind 3.   for $p\geq  \frac{N}{N-2s}$, $v\equiv0$.} \medskip
 
{\it Part 1. }  For $p\in\big(1,\frac{N}{N-2s}\big)$,   if $\int_{\Omega} v^pdx<+\infty$, it follows by \rlemma{lm 4.1-s} that 
$v=v_k$ for some $k\geq0$. Moreover, if $k=0$, $v=v_0\equiv0$ in $\Omega\setminus\{0\}$. 

If $\int_{\Omega} v^pdx=+\infty$,  from \rprop{StrSing-EF}, we have that 
$$v\geq v_\infty\quad{\rm in}\ \, \BBR^N\setminus\{0\}.$$

If  $p \in\big(0, 1+\frac{2s}{N}\big]$, then $v_\infty=+\infty$ in $\Omega$,
which implies that $v$ can be taken in this case.

{\it Part 2. } If $p\in \big(1+\frac{2s}{N}, \frac{N}{N-2s}\big)$, $v_\infty$ is a solution of $(\ref{IN12})$ satisfying
$$\lim_{x\to0} v_\infty(x)  |x|^{\frac{2s}{p-1}}=c_p $$
by \rcor{estCOR-2}. 
Moreover, by the strong maximum principle, we have that 
$$v\equiv v_\infty\quad{\rm or}\quad v>v_\infty \quad{\rm in}\ \, \Omega\setminus\{0\}.$$

If $v>v_\infty$,  we let 
$$V=\frac{v^p-v_\infty^p}{v-v_\infty}$$
and then 
$$v_\infty^{p-1} \leq V\leq v^{p-1}\quad {\rm in}\ \, \Omega\setminus\{0\}.$$
Let  $\bar v=v-v_\infty$ in $\BBR^N\setminus\{0\}$, then $\bar v$ satisfies that 
$$\lim_{x\to0} \bar v(x)|x|^{\frac{2s}{p-1}}=0$$
and
\bel{BB22}
(-\Delta)^s \bar v+V\bar v=0\quad {\rm in}\ \, \Omega\setminus\{0\},\qquad  \bar v=0 
\quad {\rm in}\ \, \Omega^c.
\ee
Note that $\vge v_\infty$ is a super solution of $(\ref{BB22})$ for any $\vge>0$ and
by the comparison principle, we have that 
$$\bar v\leq \vge v_\infty \quad {\rm in}\ \, \Omega.$$
By the arbitrary of $\vge>0$, we have that $\bar v\leq 0$. 
This is impossible from our setting that $v>v_\infty$. \\
Therefore, $v\equiv v_\infty$.

{\it Part 3. } From \rprop{estLE} we have that 
$$v(x)\leq c_{28}|x|^{-\frac{2s}{p-1}}\quad {\rm in}\ \, \Gw\setminus\{0\}.$$

 When $p>\frac{N}{N-2s}$, we have that $-\frac{2s}{p-1}>2s-N$ and $\vge |\cdot|^{2s-N}$ is a super solution of $(\ref{BB22})$ for any $\vge>0$, then comparison principle implies that $v\leq \vge |\cdot|^{2s-N}$, which by the arbitrary of $\vge>0$, implies that 
 $v=0$. 

 When $p=\frac{N}{N-2s}$, there holds $-\frac{2s}{p-1}=2s-N$.  For $m\in\BBR$  let  $ \nu_{m}$ be a smooth,  radially symmetric   function such that
 \begin{equation}\label{cl-6-1}
 \nu_{ m}(x)= (-\ln|x|)^m\quad {\rm for}\ \ 0<|x|<\frac1{e^2}\quad {\rm and}\quad v_{ m}(x) =0\quad {\rm for }\ |x|>1.
 \end{equation}
 Moreover,  $\nu_m$ is non-increasing  in $|x|$ if $m>0$.
Denote
 \begin{equation}\label{cl-6}
 w_{m}(x)=|x|^{2s-N} \nu_m(x) \ \ {\rm for}\ \,  x\in \BBR^N\setminus\{0\}.
 \end{equation}

 From \cite[Proposition 2.1]{CH} for $m\not=0$, we denote
 \begin{equation}\label{cl-6-2}
 \CI_{m}=\CC_s'(0)m,\quad \CJ_{m}=\CC_s''(0)\frac{m(m-1)}2,
 \end{equation}
 where $\CC_s'(0)<0$ and $\CC_s''(0)<0$.
Then $\CI_{m}>0$ for $m<0$, and there exist $r_0\in(0,\frac1{e^2}]$ and $c_5>0$ such that for $0<|x|<r_0$
 $$
  \Big| (-\Delta)^s w_{m}(x)  -\CI_{m}  |x|^{-N}(-\ln |x|)^{m-1}-\CJ_{m}|x|^{-N}(-\ln|x|)^{m-2}\Big| 
 \leq   c_5 |x|^{-N}(-\ln |x|)^{m-3}. 
$$
 Now we take $m_0=-\frac{N-2s}{2s}<0$ and for some $r_1>0$  there holds
 $$(-\Delta)^s w_{m_0}(x)\geq 0\quad {\rm in}\ \,  B_{r_1}(0)\setminus\{0\}.$$
Now we assume that  $\Omega\subset B_{r_1}$, then $\vge w_{m_0}$ is a super solution of $(\ref{BB22})$ for any $\vge>0$ in $B_{r_1}$, then comparison principle implies that $v\leq \vge w_{m_0}$, which by the arbitrary of $\vge>0$, implies that 
 $v=0$. 
 
When $\Omega$ is  general domain, if $v$ is a nontrivial singular solution of $(\ref{BB22})$ in $\Omega$, then by the super and sub solution method, $(\ref{BB22})$ in $B_{r_1}$ has a solution $\tilde v$ such that 
$$\lim_{x\to0}\frac{\tilde v}{v}=1.$$
 Then a contradiction could be obtained.   \qeda\medskip

\subsection*{Open problems}

\nind{\bf Problem 1.} Non-radial elements of $\CE_1$ can be formally constructed in the following way: set 
$$\BBS^N_{+,1}=\BBS^N\cap\{(x_1,...,x_N,z):x_1>0,z>0\}.$$
Construct by minimization a positive function $\gw$ satisfying 
\bel{CC1}\left.\BA{lll}\displaystyle
 \CA_s[\gw]+\Gl_{s,p,{\scriptscriptstyle N}}\gw=0\quad&\text{in }\,\BBS^N_{+,1}\\[2mm]
\displaystyle\phantom{\,   \Gl_{s,p,{\scriptscriptstyle N}}}
\frac{\prt \gw}{\prt \gn^s}+ \gw^p=0\quad&\text{in }\, \BBS^{N-1}\cap\{x_1>0\}\\[2mm]
\displaystyle\phantom{\Gl_{s,p,{\scriptscriptstyle N}}\gw+\CA[]\ }
\gw=0\quad&\text{in }\,\BBS^{N}_+\cap\{x_1=0\}.
\EA\right.\ee
Then the function $\tilde w$ defined in  $\BBS^N_+$ by
$$
\tilde \gw=\left\{\BA{lll}\gw(\gs)\quad&\text {in } \,\BBS^N_{+,1}\\[2mm]
-\gw(\tilde\gs)\quad&\text {in } \, \BBS^N_{+}\cap\{x_1<0\},
\EA\right.
$$
where $\tilde\gs$ is the symmetric of $\gs$ with respect to the plane $x_1=0$, is a signed solution of $(\ref{I-5})$ with $\ge=1$. It vanishes on the plane $\{(x,z):x_1=0\}$. It could be interesting to see if the condition for existence of such a solution coincides with the condition $p<p^*$. 
\medskip

\nind{\bf Problem 2.} It could be interesting to investigate the type of operators in $\BBR^N_+$ which could lead by a trace process somehow similar to the one of \cite{CS} to perturbations of $(-\Gd)^s$ such as $(-\Gd+mI)^s$ or to Hardy type operators $(-\Gd)^s+\gm|x|^{-2s}$. This could be usefull for studying a wider class of nonlinear semilinear equations involving fractional type-Laplacians. \medskip

 \nind{\bf Problem 3.} The description of the singular behaviour in the case $p=\frac{N}{N-2s}$ is more complicated due to the fact that $\Gl_{s,p,{\scriptscriptstyle N}}=0$ and the linear blow-up rate $2s-N$ coincides with the nonlinear blow-up rate $\frac{2s}{p-1}$. This phenomenon has been described by Aviles 
 \cite{Av} in the case $s=1$. In a very recent article \cite{WW} Wei and Wu proved that the function $v$ either is smooth at $0$ or satisfies
 $$c_1|x|^{2s-N}(-\ln |x|)^{\frac{2s-N}{2s}}\leq v(x)\leq c_2|x|^{2s-N}(-\ln |x|)^{\frac{2s-N}{2s}}\quad\text{for all }x\in B_1\setminus\{0\}.
 $$
 for some constants $c_1,c_2>0$. 
 We conjecture that the following result holds: either the function $v$ can be extended as a smooth solution at $0$, or
 $$
 \lim_{x\to 0} |x|^{N-2s}(-\ln |x|)^{\frac{N-2s}{2s}}v(x)=C(N,s)
 $$
 for some positive explicit constant $C(N,s)$. \medskip

\nind{\bf Problem 4.} In an impotant article, Simon \cite{Si} studied the unique limit problem for bounded solutions of  analytic functional. His results were used in 
\cite{BV-V} to obtained the precise asymptotics of solutions of  
$$
-\Gd u=\gl e^u\quad \text{in }B_1\setminus\{0\}\subset\BBR^3
$$
 with $\gl>0$ satisfying $|x|^2e^{u(x)}\in L^{\infty}(B_\frac 12)$, or positive solutions of 
$$
-\Gd u=u^p+\frac{\ell}{|x|^2}u\quad \text{in }B_1\setminus\{0\}\subset\BBR^N
$$
with $p>1$ satisfying $|x|^{\frac{2}{p-1}}\in L^{\infty}(B_\frac 12)$. The technique was to reduce the problem to a semilinear ellictic equations in an infinite cylinder and to prove that the limit set of the renormalized corresponding solutions was reduced to a single element whatever is the structure of the set of possible limits (or self-similar solutions). It is natural to investigate similar questions concerning 
$$
\left.\BA{lll}\displaystyle
u_{zz}+\frac{1-2s}{z}u_z+\Gd_xu+\frac{\ell}{|x|^2+z^2}u=0\quad\text{in }\, \BBR^{N+1}_+ \\[3.5mm]\displaystyle
\phantom{--;;}-\lim_{z\to 0}z^{1-2s}u_z(x,z)-u^p(x,0)=0\quad\text{in } B_1\setminus\{0\}\subset \BBR^N.
\EA\right.$$
Assuming that a solution $u$ satisfies $(|x|^2+z^2)^{\frac{s}{p-1}}\in L^{\infty}(B_\frac 12)$ and using the variable $(\ref{I17})$ the new function 
$w(t,\gs)=r^{\frac {2s}{p-1}}u(r,\gs)$, $t=\ln r$, satisfies 
 \bel{I18}
\left.\BA{lll}\displaystyle
w_{tt}+\Gth_{s,p,{\scriptscriptstyle N}}w_t+\left(\Gl_{s,p,{\scriptscriptstyle N}}+\ell\right)w+\CA_s[w]=0\ \ &\text{in }\BBR\ti \BBS^N_+\\[2mm]
\phantom{-------------}
\myfrac{\prt w}{\prt\gn^s}+w^p=0\ \ &\text{in }\BBR\ti \BBS^{N-1},
\EA\right. \ee
where $\Gth_{s,p,{\scriptscriptstyle N}}$ and $\Gl_{s,p,{\scriptscriptstyle N}}$ are unchanged. 

 \bigskip\bigskip
 
\nind{\bf Aknowledgements:} {\small H. Chen is supported by the Natural Science Foundation of China, No. 12071189, by Jiangxi Province Science Fund, No. 20202ACBL201001 and 20212ACB211005.} The author are grateful to the referee for his careful reading of the article and for drawing to their attention reference 
\cite{WW} which partially solved Problem 3 above.



\begin {thebibliography}{99}

\bibitem{AH}  D. R. Adams, L. I. Hedberg,  Function spaces and potential theory. {\it Grundlehren  Math. Wissen.  314  Springer}  (1996).

\bibitem{Av} P. Aviles, Local behavior of solutions of some elliptic equations. {\it Comm. Math. Phys. 108}, 177--192 (1987).

\bibitem{BV-V} M.F. Bidaut-V\'eron, L. V\'eron, Nonlinear elliptic equations on compact Riemannian manifolds and asymptotics of Emden equations. {\it Invent. Math. 106}, 489--539 (1991).

\bibitem{BV} O. Boukarabila, L. V\'eron, Nonlinear boundary value problems relative to harmonic functions. {\it Nonlinear Analysis  201}, 112090 (2020).

\bibitem{BLb}  H. Brezis, E. H. Lieb, Long range atomic potentials in Thomas-Fermi theory. {\it Comm. Math. Phys.  65}, 231--246 (1979).

\bibitem{BV}  H. Brezis, L. V\'eron,  Removable Singularities for Some Nonlinear Elliptic Equations. {\it Arch. Rat. Mech. Anal.  75}, 1--6 (1981).

\bibitem {BL} H. Brezis, P. L. Lions, A note  on isolated singularities for linear elliptic equations.  {\it Math. Anal. Appl. Part A. Adv. in Math. Supplementary studies  Vol 7A}, 363--366 (1981).

\bibitem{CS} L. Caffarelli, L. Silvestre, An extension problem related to the fractional Laplacian. {\it Comm. Part.  Diff.  Eq. 32}, 1245--1260 (2007).

\bibitem{CGS} L. Caffarelli, B. Gidas, J. Spruck, Asymptotic symmetry and local behavior of semilinear elliptic equations with critical Sobolev growth. {\it  Comm. Pure Appl. Math.   42}, 271--297 (1989).

\bibitem{CJSX} L. Caffarelli, T. Jin, Y. Sire, J. Xiong, Local analysis of solutions of fractional semi-linear elliptic equations with isolated singularities. {\it Arch. Rat. Mech. Anal.  213}, 245--268 (2014).

\bibitem{Ch} S. Chandrasekhar, Introduction to Stellar Structure. {\it University of Chicago Press}  (1939).

\bibitem {CFQ} H. Chen, P. Felmer, A. Quaas, Large solutions to elliptic  equations involving fractional Laplacian.
{\it Ann. Inst. H.  Poincar\'e--AN  32}, 1199--1228 (2015). 

\bibitem{CH} H. Chen,  H. Hajaiej, On isolated singularities for fractional Lane-Emden equation in the Serrin's critical case, {\it arXiv:2109.07085} (2022).
 

\bibitem{CMV} X. Chen, H. Matano, L. V\'eron,  Anisotropic Singularities of Solutions of Nonlinear Elliptic Equations in $\BBR^2$.   {\it J. Funct. Anal.  83}, 50--97 (1989).

\bibitem{CV} H. Chen, L. V\'eron,  Semilinear fractional elliptic equations involving measures. {\it J. Diff. Eq.  257}, 1457--1486 (2014).

\bibitem{CV1} H. Chen, L. V\'eron,  Weakly and strongly singular solutions of semilinear fractional elliptic equations. {\it Asymp. Anal.  88}, 165--184 (2014).

\bibitem{CW} H. Chen, T.  Weth, The Poisson problem for the fractional Hardy operator: distributional identities and singular solutions. {\it Trans. Amer. Math. Soc.  374}, 6881--6925 (2021).

\bibitem{CY} H. Chen,  J. Yang,  Semilinear fractional elliptic equations with measures in unbounded domain. {\it Nonlinear Anal. 145},   118--142 (2016). 

\bibitem{Fa} M. Fall, Semilinear elliptic equations for the fractional Laplacian with Hardy potential. {\it Nonlinear Analysis  193}, 111311 (2020).

\bibitem{FT} Y. Fang, D. Tang, Method of sub-super solutions for fractional elliptic equations.  {\it Discr. Cont. Dyn. Sys.  23}, 3153--3165 (2018).

  \bibitem {FQT} P.  Felmer, A. Quass, J. Tan, 
Positive solutions of the nonlinear Schr\"odinger equation with the fractional Laplacian.  {\it Proceed.  Royal Society of Edinburgh, 142A,} 1237--1262 (2012).



\bibitem{Ga}  N. Garofalo, Fractional thoughts in  New Developments in the Analysis of Nonlocal Operators. {\it Contemp. Math.  723}, 1--136 (2019).

\bibitem{GNN} B. Gidas, W. Ni, L. Nirenberg, Symmetry and related properties via the maximum principle. {\it Comm. Math. Phys.  68}, 209--243 (1979).

\bibitem{GS} B. Gidas, J. Spruck, A priori bounds for positive solutions of nonlinear elliptic equations. {\it  Comm. Part. Diff. Eq.  6}, 883--901 (1981).

\bibitem{GS1} B. Gidas, J. Spruck, Global and local behaviour of positive solutions of nonlinear elliptic equations. {\it  Comm. Pure Appl. Math.  34}, 525--598 (1981).

\bibitem{FKS} E. Fabes, C. Kenig, R. Serapioni, The local regularity of solutions of degenerate elliptic equations. {\it  Comm. Part. Diff. Eq.  7}, 77--116  (1982).

\bibitem{Fo} R. Fowler, Further studies on Emden's and similar diffferential equations. {\it Quart. J. Math. 2}, 259--288 (1931).

\bibitem{JLX} T. Jin, Y.Y. Li, J. Xiong, On a fractional Nirenberg problem, part I: blow up analysis and compactness of solutions.
{\it J. Eur. Math. Soc.  16}, 1111--1171 (2014).

\bibitem{LWX} C. Li, Z. Wu, H. Xu, Maximum principles and B\^ocher type theorems. {\it Proc Natl Acad Sci U S A. 115}, 6976--6979 (2018).

\bibitem{Li} P. L. Lions, Isolated singularities in semilinear problems. {\it J. Diff. Eq.  38}, 441--450  (1980).


\bibitem{MV} M. Marcus, L. V\'eron, Uniqueness and asymptotic behaviour of solutions with boundary blow-up for a class of nonlinear elliptic equations. {\it Ann. Inst. H. Poincar\'e--AN 14}, 237--274 (1997).

\bibitem{Me} R. Menon,   Integro-differential operators: connections to degenerate elliptic equations and some free boundary problems. {\it PhD Thesis, Michigan State University} (2020).

\bibitem {Si} L. Simon, Isolated singularities of extrema of geometric variational problems, In: Giusti E. (ed.)  Harmonic Mappings and Minimal Immersions. {\it Lect. Notes Math.    1161}, 206--277, Berlin, Heidelberg, New-York (1985).

\bibitem{So} A. Sommerfeld, Asymptotische integration der differential-gleichung des Thomas-Fermischen  atoms. {\it Zeitschrift f\"ur Physik 78}, 283--308 (1932).

\bibitem{Ve} L. V\'eron, Singular solutions of some nonlinear elliptic equations. {\it Nonlinear Analysis, Th. Meth. Appl.  5}, 225--242 (1981).

\bibitem {WW} J. Wei, K. Wu, Local behavior of solutions to a fractional equation with isolated singularity and critical Serrin exponent. {\it Dis. Cont. Dyn. Syst. 42}, 4031-4050 (2022). 

\bibitem{YZ} H. Yang, W. Zou, On isolated singularities of fractional semi-linear elliptic equations. {\it Ann. Inst. H. Poincar\'e--AN. 38}, 403--420 (2021).

\bibitem{YZ1} H. Yang, W. Zou,  Sharp blow up estimates and precise asymptotic behavior of singular positive solutions to fractional Hardy-H\'enon equations, {\it J. Diff. Eq. 278}, 393--429 (2021).
 
\end{thebibliography}
 \end {document}